\documentclass[11pt,reqno]{amsart}
\usepackage{amsmath,amssymb,amsthm,mathtools,microtype,mathrsfs,bm}
\usepackage{tikz}
\usepackage{needspace}
\usepackage{euscript}
\let\eu\EuScript
\DeclareFontFamily{U}{stabilizationcal}{\skewchar\font=45}
\DeclareFontShape{U}{stabilizationcal}{m}{n}{<-> BOONDOX-r-calo}{}
\DeclareMathAlphabet{\stabilizationcal}{U}{stabilizationcal}{m}{n}
\newcommand{\stratumZ}{\stabilizationcal{Z}}

\usepackage[colorlinks=true,linkcolor=blue,citecolor=red]{hyperref}
\usetikzlibrary{arrows.meta}
\numberwithin{equation}{section}
\numberwithin{figure}{section}
\newtheorem{theorem}{Theorem}[section]
\newtheorem{proposition}[theorem]{Proposition}
\newtheorem{lemma}[theorem]{Lemma}
\newtheorem{corollary}[theorem]{Corollary}
\newtheorem{maintheorem}{Theorem}

\theoremstyle{definition}
\newtheorem{definition}[theorem]{Definition}
\newtheorem{remark}[theorem]{Remark}
\newtheorem*{unnumberedremark}{Remark}
\newtheorem{example}[theorem]{Example}
\newcommand\xqed[1]{\leavevmode\unskip\penalty9999 \hbox{}\nobreak\hfill
  \quad\hbox{#1}}
\newcommand\qetr{\xqed{$\triangle$}}
\newcommand\qrem{\xqed{$\spadesuit$}}
\DeclareMathOperator{\Var}{Var}

\DeclareMathOperator{\Log}{Log}
\DeclareMathOperator{\Sym}{Sym}
\DeclareMathOperator{\des}{des}
\DeclareMathOperator{\Aff}{Aff}
\DeclareMathOperator{\Aut}{Aut}
\DeclareMathOperator{\Vertices}{Vert}
\DeclareMathOperator{\Edge}{Edge}
\newcommand{\C}{\mathbb C}
\newcommand{\Q}{\mathbb Q}
\newcommand{\Z}{\mathbb Z}
\newcommand{\LL}{\mathbb L}
\newcommand{\Mbar}{\overline{\mathcal M}}
\newcommand{\GG}{\mathbb G_m}
\newcommand{\Vspace}{\mathcal V}
\allowdisplaybreaks

\begin{document}
\title[Heavy--Light Wall Crossings and Motivic Asymptotics
of Moduli Spaces of Curves, II]{Heavy--Light Wall Crossings and Motivic Asymptotics
of Moduli Spaces of Curves\\[2mm]
Part II: Tameness of Motivic Measures}

\author[Giordano Cotti]{Giordano Cotti$\>^{\circ}$}

{\let\thefootnote\relax
\footnotetext{
\noindent
$^\circ\>$\textit{ E-mail}:  giordano.cotti@tecnico.ulisboa.pt}}

\subjclass[2020]{Primary 14H10; Secondary 14C15, 14C30, 05A16}
\keywords{Moduli spaces of curves, Hassett spaces, motivic measures,
wall crossing, stabilization graphs, residual quotients, tameness,
Hodge--Deligne polynomials, singularity analysis}
\maketitle
\begin{center}
\textit{ 
$^\circ\>$Grupo de F\'isica Matem\'atica \\
Departamento de Matem\'atica, Instituto Superior T\'ecnico\\
Av. Rovisco Pais, 1049-001 Lisboa, Portugal\/\\}
\end{center}
\vskip1cm
\begin{abstract}
This paper is concerned with the study of the fixed-genus asymptotics
of motivic measures of moduli
spaces of stable pointed curves and their minimal all-light Hassett
compactifications as the number of markings tends to infinity.
A motivic wall-crossing argument gives a universal change of variables
between the two generating series, depending only on the value of
the Lefschetz class. For an arbitrary motivic measure, we use
stabilization graphs and support patterns to express the all-light
series in terms of universal chain factors and finite residual quotient
families. This formal decomposition isolates the rational trivalent
contribution without requiring a fixed-locus formula for finite quotients.
To pass from formal identities to asymptotics, we formulate analytic
hypotheses (the \emph{$\mathrm H_g$-hypotheses}) for scalar realizations
with positive real Lefschetz value. We introduce \emph{vertex tameness}
and \emph{residual-quotient tameness}, whose conjunction,
\emph{graph tameness}, provides geometric sufficient conditions
for the $\mathrm H_g$-hypotheses. In genus at least two, this yields explicit
leading asymptotics in both chambers, depending only on the genus
and the Lefschetz value; genera zero and one are treated separately.
As an application, we prove graph tameness for Hodge--Deligne
evaluations and deduce fixed-product universality and asymptotic
concentration of Hodge numbers on the diagonal.
\end{abstract}

\tableofcontents

\section{Introduction}
\label{sec:introduction}

\medskip\noindent\textbf{1.1. The problem and the generating series.}
The comparison between ordinary and minimal-weight compactifications of
moduli spaces of curves raises a natural question:
\begin{center}
\begin{minipage}{.9\textwidth}
\centering\itshape
Which features of their fixed-genus asymptotics are determined by the
geometry of the marking weights, and which depend on the invariant
being computed?
\end{minipage}
\end{center}
For the ordinary Euler characteristic, this comparison was studied in
\cite{CottiPaperI} by combining the geometry of stabilization graphs
with a change of variables between the two chambers. Here we consider
the more general problem in which the Euler characteristic is replaced
by an arbitrary \emph{motivic measure}
\[
\mu:K_0(\Var_{\C})\longrightarrow\mathsf R,
\qquad \lambda:=\mu(\mathbb L).
\]
Here $\mathsf R$ is a commutative ring, $\mu$ is a unital ring
homomorphism, and $\mathbb L=[\mathbb A^1]$. The value $\lambda$
of the \emph{Lefschetz class} will govern the universal
part of the construction, while the remaining contributions retain
information about the chosen measure. Throughout, the invariants are
those of coarse moduli spaces.

Let $\overline M_{g,n}$ be the ordinary moduli space of stable curves
with $n$ labeled markings, and let $\overline M_{g,(0^+)^n}$ denote
the Hassett space in which all markings have sufficiently small positive
weight, so that every subset is allowed to coincide
\cite{DeligneMumford,Knudsen,Hassett}. For exponential generating
series we extend coefficients to
$\mathsf R_{\Q}:=\mathsf R\otimes_{\Z}\Q$; values of $\mu$ are
understood through the natural map $\mathsf R\to\mathsf R_{\Q}$.
This convention is used throughout, and for the universal measure we
write $\eu K:=K_0(\Var_{\C})\otimes_{\Z}\Q$. The series to be compared are
\begin{equation}
A_g^\mu(z):=\sum_{n\ge0}
\mu(\overline M_{g,(0^+)^n})\frac{z^n}{n!},\qquad
G_g^\mu(t):=\sum_{n\ge0}
\mu(\overline M_{g,n})\frac{t^n}{n!}.
\label{eq:intro-series}
\end{equation}
These series are defined for $g\ge1$, with zero constant term in
genus one. There is no stable minimal all-light chamber in genus zero;
its ordinary series starts at $n=3$ and will enter through the universal
change of coordinates.

\medskip\noindent\textbf{1.2. The geometric decomposition and the change of variables.}
For $g\ge2$, forgetting the light markings contracts the maximal
chains of two-valent rational components and produces a stable
unmarked curve. Its dual graph is the \emph{stabilization graph}.
Only finitely many such graphs occur in a fixed genus, even though
the lengths of the contracted chains are unbounded
(Section~\ref{sec:stabilizationgraphs}). As in
\cite[Section~4]{CottiPaperI}, we reconstruct the all-light spaces by
\emph{decorated expansions}: smooth component moduli give the
\emph{vertex factors}, and ordered rational chains give the
\emph{edge factors}. These constructions are introduced in
Sections~\ref{sec:decorated-expansions}--\ref{sec:parameter-and-quotient}.
The chain calculation of Sections~\ref{sec:single-chain-component}--\ref{sec:motivicpropagator}
gives the \emph{propagator}\footnote{As in \cite{CottiPaperI}, this
terminology is borrowed from the physics of Feynman diagrams: a
propagator encodes the contribution associated with an edge.}
\[
P_\lambda(z)=
\left(1-\sum_{m\ge1}(\lambda-1)^{m-1}\frac{z^m}{m!}\right)^{-1},
\qquad P_1(z)=(1-z)^{-1}.
\]
We prove in Section~\ref{sec:losevmaninpropagator} that this is also
the generating series of the measures of the \emph{Losev--Manin
spaces} $\overline L_n=\overline M_{0,2|n}$, with a formal constant
term for the empty chain. For their geometry, see \cite{LosevManin}
and \cite[Section~7.1]{CottiRHB}.

A second part of the strategy relates the two marking chambers.
We prove a single-wall identity for coarse classes in the Grothendieck
ring (Section~\ref{sec:single-wall}), then apply it successively while
lowering the weights of the heavy markings
(Section~\ref{sec:wall}). The resulting change of variables is
independent of the genus and depends on the measure only through
$\lambda$. Moreover, it is determined by the ordinary genus-zero
series, and its derivative recovers the edge propagator.
To state these links together, set
\[
G_0^\mu(t):=\sum_{n\ge3}\mu(\overline M_{0,n})\frac{t^n}{n!},
\qquad G_0^{\mu,\bullet}(t):=1+t+(G_0^\mu)'(t),
\]
where $1+t$ are formal unstable initial terms.

\begin{maintheorem}[Universal transformation and its genus-zero geometry]
\label{intro:universal}
There is a universal coordinate $\Phi_\lambda(z)=z+O(z^2)$ with
coefficients in $\Q[\lambda]$, and formal inverse $\Psi_\lambda$,
such that for every motivic measure $\mu$ and every $g\ge1$,
\begin{equation}
A_g^\mu=G_g^\mu\circ\Phi_\lambda,
\qquad G_g^\mu=A_g^\mu\circ\Psi_\lambda.
\label{eq:intro-chamber-change}
\end{equation}
The coordinate, the ordinary genus-zero series, and the propagator
satisfy the formal identities
\[
e^{\Psi_\lambda(t)}=G_0^{\mu,\bullet}(t),
\qquad\qquad
P_\lambda(z)=1+\sum_{n\ge1}\mu(\overline L_n)\frac{z^n}{n!}
=\frac{e^z}{\Phi_\lambda'(z)}.
\]
These statements hold already for the universal classes in $\eu K$
and hence after every motivic realization in $\mathsf R_{\Q}$.
\end{maintheorem}
The chamber identities are Theorem~\ref{thm:pde} and
Corollary~\ref{cor:realized-composition}; the Losev--Manin interpretation
is Proposition~\ref{prop:LM-interpretation}. The genus-zero and derivative
identities are Propositions~\ref{prop:Psi-genuszero}
and~\ref{prop:propagator-coordinate}, in
Sections~\ref{sec:genuszero-coordinate}--\ref{sec:propagator-coordinate};
Corollary~\ref{cor:genuszero-heavy-light} writes the chamber change
directly in terms of $G_0^{\mu,\bullet}$.

The universal coordinate has the explicit expression
\[
\Phi_\lambda(z)=
\frac{(e^{\lambda z}-1)/\lambda-\lambda(e^z-1)}{1-\lambda},
\qquad
\Phi_0(z)=z,\quad \Phi_1(z)=(2-z)e^z-2.
\]
The fraction is a shorthand for its
formal expansion: every coefficient is polynomial in $\lambda$.
Its apparent singularities at $\lambda=0,1$ are removable
coefficientwise, so no inversion of $\lambda$ or $1-\lambda$ in the
coefficient ring is required.

Thus the genus-zero spaces determine the coordinate in every positive
genus, while the Losev--Manin spaces supply the chain factors in the
graph decomposition. Their connection through $\Phi_\lambda'$ also
links the singularities governing the two chambers. The identification
with Losev--Manin spaces is an equality of Grothendieck classes and
series: the disjoint union of chain strata and their proper
compactification are not identified as varieties.

S.\,Kannan, S.\,Serpente, and C.\,H.\,Yun first obtained the change of variables for
symmetric-group-equivariant Hodge--Deligne generating series
\cite[Theorem~A]{KSY}. Their non-equivariant specialization is stated
in \cite[Corollary~B]{KSY}; the ordinary Euler transformation follows
by setting $u=v=1$, taking the removable limit in the explicit formula.
Equivariant invariants retain the marking-permutation actions and
require an equivariant framework. Our wall-crossing argument proves
the non-equivariant transformation in $\eu K$, before choosing a
realization. At $\lambda=1$ it gives the Euler coordinate
of \cite[Section~2.4]{CottiPaperI}; its
Hodge--Deligne specialization also gives the heavy/light recursion
of \cite[Theorem~3.1]{Kiem}.

\medskip\noindent\textbf{1.3. Graph symmetries and finite residual quotients.}
For $g\ge2$, the geometric decomposition according to stabilization
graphs gives a finite stratification of each all-light moduli space
by locally closed loci\footnote{For algebraic spaces or stacks,
$\bigsqcup$ denotes a decomposition into locally closed pieces, which
may be related by specialization; $\coprod$ denotes their disjoint
union as algebraic spaces or stacks, in which each piece is open and
closed.}
(Section~\ref{sec:decomposition-by-graph}).
By additivity of $\mu$, it induces a
decomposition of the generating series:
\[
A_g^\mu(z)=\sum_{\Gamma\in\mathfrak C_g}A_{g,\Gamma}^\mu(z),
\]
where $\mathfrak C_g$ is the finite set of stabilization graphs of
genus $g$, and $A_{g,\Gamma}^\mu$ collects the measures of the loci
with stabilization graph isomorphic to $\Gamma$.

To describe one such locus, we follow
\cite[Section~4.3]{CottiPaperI} and use \emph{decorations} of a fixed
representative $\Gamma$ (Section~\ref{sec:decorated-expansions}).
A decoration specifies the rational chains
inserted along its edges and the distribution of the labeled light
markings among the graph vertices and the inserted components; every
inserted component carries a nonempty set of light labels. Write
$\mathscr D_n(\Gamma)$ for the set of decorations with light labels
$[n]$. For each decoration $\delta$, the component moduli form a coarse
\emph{parameter space} $X_\delta$, whose points give the curves before
graph symmetries are divided out.

The group $\Aut\Gamma$ acts on $\mathscr D_n(\Gamma)$ by changing
the identification with the fixed graph, while keeping every light
label fixed. Decorations in the same orbit describe the same geometric
stratum. If $G_\delta\subseteq\Aut\Gamma$ is the stabilizer of
$\delta$, the stratum is the finite quotient $X_\delta/G_\delta$;
these quotient loci form a finite locally closed stratification of
the locus belonging to $\Gamma$, as in
\cite[Lemma~4.11]{CottiPaperI}. The quotient description is established
in Section~\ref{sec:parameter-and-quotient}; consequently,
\[
A_{g,\Gamma}^\mu(z)
=\sum_{n\ge0}\frac{z^n}{n!}
\sum_{[\delta]\in\mathscr D_n(\Gamma)/\Aut\Gamma}
\mu(X_\delta/G_\delta).
\]

We next separate a contribution that can be computed directly from
the vertex and edge factors. The \emph{parameter series} sums over all
decorations, before passage to orbits:
\[
A_{\Gamma,\mathrm{par}}^\mu(z)
:=\sum_{n\ge0}\frac{z^n}{n!}
\sum_{\delta\in\mathscr D_n(\Gamma)}\mu(X_\delta)
=\left(\prod_{v\in\Vertices(\Gamma)}V_{h_v,d_v}^\mu(z)\right)
P_\lambda(z)^{E(\Gamma)}.
\]
Here $V_{h,d}^\mu$ is the \emph{vertex series} for a smooth component
of genus $h$ with $d$ nodal branches and arbitrary light markings (see \eqref{eq:vertexseries}), and
$E(\Gamma)$ is the number of edges. The product follows from the
independent allocation of the labels among vertices and edge chains.

An orbit represented by $\delta$ has
$|\Aut\Gamma|/|G_\delta|$ elements, so its contribution to
$A_{\Gamma,\mathrm{par}}^\mu/|\Aut\Gamma|$ is
$\mu(X_\delta)/|G_\delta|$. Comparing this with the actual quotient
contribution gives
\eqref{eq:graphorbitdecomposition}--\eqref{eq:nonfreecorrection}:
\begin{align*}
A_{g,\Gamma}^\mu(z)
&=\frac1{|\Aut\Gamma|}A_{\Gamma,\mathrm{par}}^\mu(z)
+\mathcal R_\Gamma^\mu(z),\\
n![z^n]\mathcal R_\Gamma^\mu(z)
&=\sum_{\substack{[\delta]\in\mathscr D_n(\Gamma)/\Aut\Gamma\\
G_\delta\ne\{1\}}}
\left(\mu(X_\delta/G_\delta)
-\frac{\mu(X_\delta)}{|G_\delta|}\right).
\end{align*}
Thus the remainder records precisely the differences between the
measures of the quotient spaces and the normalized measures of their
parameter spaces. Only decorations with nontrivial stabilizer can
contribute. For the Euler characteristic, Burnside--Lefschetz rewrites
these differences as sums of nonidentity fixed-locus contributions,
which can then be studied analytically
\cite[Sections~4.4 and~6.1]{CottiPaperI}. For an arbitrary motivic
measure this formula is unavailable, and we must reorganize
$\mathcal R_\Gamma^\mu$ while retaining the quotient geometry.

A \emph{support pattern} records which vertices carry light labels
and sorts edges by the number of inserted components: zero, one,
or at least two. This determines the discrete stabilizer and identifies
the chain factors on which it acts trivially
(Section~\ref{sec:motivic-residual-types}). Separating these factors
and expanding their series leaves a finite collection of families
of spaces with prescribed group actions. For each choice $\tau$ in
this expansion, write $Y_{\Gamma,\tau}^\mu$ for the series of the
remaining spaces and $Q_{\Gamma,\tau}^\mu$ for the series of their
quotients by the inherited group $G_{\Gamma,\tau}$. These
\emph{finite residual quotient families} contribute a quotient
difference times one factor of $P_\lambda$ for each extracted edge
(Section~\ref{sec:finite-residual-families}). The resulting identity is
\[
\mathcal R_\Gamma^\mu(z)
=\sum_{\tau\in\mathscr T_\Gamma}c_{\Gamma,\tau}
\left(Q_{\Gamma,\tau}^\mu(z)
-\frac{Y_{\Gamma,\tau}^\mu(z)}{|G_{\Gamma,\tau}|}\right)
P_\lambda(z)^{k_{\Gamma,\tau}},
\qquad k_{\Gamma,\tau}\le E(\Gamma)-1,
\]
where $\mathscr T_\Gamma$ indexes the finitely many choices,
$k_{\Gamma,\tau}$ counts the extracted edges, and $c_{\Gamma,\tau}$
are the explicit rational coefficients of
Lemma~\ref{lem:residual-decomposition}.

Finally, we sum over stabilization graphs. A rational trivalent graph
has $2g-2$ vertices and $3g-3$ edges, and each vertex contributes
$V_{0,3}^\mu(z)=e^{(\lambda-2)z}$. Its normalized parameter contribution
therefore has the maximal propagator exponent $3g-3$. All other graph
parameter contributions, as well as all residual correction terms,
have propagator exponents at most $3g-4$. These observations yield the following formal result.

\begin{maintheorem}[Formal stabilization-graph expansion]
\label{intro:formal}
For every motivic measure $\mu$ and every fixed $g\ge2$, one has
\begin{equation}
A_g^\mu(z)=\gamma_g e^{2(g-1)(\lambda-2)z}
P_\lambda(z)^{3g-3}+\mathcal R_g^\mu(z),
\qquad
\gamma_g=\sum_{\Gamma\in\mathfrak T_g}\frac1{|\Aut\Gamma|},
\label{eq:intro-graph-expansion}
\end{equation}
as an identity of formal power series in $\mathsf R_{\Q}[[z]]$.
Here $\mathfrak T_g$ is the set of rational trivalent stabilization
graphs.\footnote{A \emph{rational trivalent graph} is a connected
multigraph whose vertices all have genus decoration zero and valence
three. Loops and multiple edges are allowed, with each loop
contributing two to the valence of its vertex. The genus of such a
graph is its first Betti number
$b_1(\Gamma)=|E(\Gamma)|-|V(\Gamma)|+1$.} The remainder admits a finite representation by products of
vertex series and by residual quotient differences of the form
$Q_{\Gamma,\tau}^\mu-Y_{\Gamma,\tau}^\mu/|G_{\Gamma,\tau}|$,
multiplied by powers of $P_\lambda$ of exponent at most $3g-4$.
\end{maintheorem}
This is Theorem~\ref{thm:motivicgraphexpansion}, obtained by summing
Proposition~\ref{prop:motivicgraph}; its residual representation is
Lemma~\ref{lem:residual-decomposition}. These results are proved in
Sections~\ref{sec:finite-residual-families}
and~\ref{sec:leadingmotivicpropagator}.

Identity~\eqref{eq:intro-graph-expansion} makes no convergence
assertion. The smaller propagator exponents in its remainder become
analytically subdominant only under additional control of the vertex
and residual series.

\medskip\noindent\textbf{1.4. From formal identities to asymptotics.}
Theorem~\ref{intro:formal} is an identity of formal series. The bound
on propagator exponents alone gives no convergence information about
the vertex and residual series. To draw asymptotic conclusions, we
pass to scalar realizations $\mu:K_0(\Var_{\C})\to\C$ and specify
the analytic information needed in addition to the formal decomposition.

For $\lambda>0$, the universal geometry of
Section~\ref{sec:universal-analytic-geometry} gives a unique dominant
propagator pole and a square-root singularity of the inverse coordinate.
Their locations and the inverse's leading local coefficient are
\[
\sigma_\lambda=\frac{\log\lambda}{\lambda-1},\quad \sigma_1=1,
\qquad \rho_\lambda=\Phi_\lambda(\sigma_\lambda),\qquad
a_\lambda=\sqrt{\frac{2\rho_\lambda}{\lambda e^{\sigma_\lambda}}},
\]
with $0<\rho_\lambda<\sigma_\lambda$
(Proposition~\ref{prop:universal-singularity-pair}). Hypotheses
$(\mathrm H_1)$ and $(\mathrm H_g)$, stated in
Section~\ref{sec:minimal-analytic-hypothesis}, concern convergence near
zero and continuation to suitable $\Delta$-domains.\footnote{We follow
the terminology of P.\,Flajolet and R.\,Sedgewick \cite[Chapter~VI]{FS}. The domains used here are
defined in \eqref{eq:Delta-domain}.} For $g\ge2$,
$(\mathrm H_g)$ requires the total remainder in
\eqref{eq:intro-graph-expansion} to have singular growth of order at
most $3g-4$. These conditions allow the standard
\emph{singularity-transfer theorems} of P.\,Flajolet and A.\,Odlyzko
\cite{FO90}, as developed in \cite[Chapter~VI]{FS}, to turn local
singular expansions into coefficient asymptotics.

We also give geometric sufficient conditions for these analytic
hypotheses. A realization is \emph{vertex-tame} if for every fixed
stable vertex type $(h,d)$ there are constants $B,C>0$ such that
\[
|\mu(\mathcal V_{h,d,n})|\le BC^n\qquad(n\ge0).
\]
It is \emph{residual-quotient tame} if for each prescribed residual
family and its group action there are constants $B,C>0$ such that
\[
|\mu(Y_{\Gamma,\tau,n})|
+|\mu(Y_{\Gamma,\tau,n}/G_{\Gamma,\tau})|
\le BC^n\qquad(n\ge0).
\]
The constants may depend on the realization and the indicated type,
and may differ between the two bounds, but are independent of $n$.
Their conjunction is \emph{graph tameness}
(Definition~\ref{def:graph-tame}, Section~\ref{sec:vertex-tameness}).
These estimates control precisely the coefficients collected in the
vertex, residual, and quotient exponential generating series, making
them entire of exponential type.\footnote{An entire function $F$ is
of \emph{exponential type} if there exist $A,T>0$ with
$|F(z)|\le Ae^{T|z|}$ for every $z\in\C$. The coefficient estimate
implying this bound is written out after
Definition~\ref{def:graph-tame}.}

\begin{maintheorem}[Tameness and fixed-genus asymptotics]
\label{intro:asymptotics}
Let $\mu:K_0(\Var_{\C})\to\C$ have Lefschetz value $\lambda>0$.
Vertex tameness implies $(\mathrm H_1)$, and graph tameness implies
$(\mathrm H_g)$ for every $g\ge2$.
For a fixed $g\ge2$, assume $(\mathrm H_g)$ and put $m=3g-3$.
Then, as $n\to\infty$,
\begin{align*}
\mu(\overline M_{g,(0^+)^n})
&\sim\frac{\gamma_g e^{2(g-1)(\lambda-2)\sigma_\lambda}}
{\Gamma(m)(\lambda\sigma_\lambda)^m}
\,n!\,\sigma_\lambda^{-n}n^{m-1},\\
\mu(\overline M_{g,n})
&\sim\frac{\gamma_g e^{2(g-1)(\lambda-2)\sigma_\lambda}}
{\Gamma(m/2)(\lambda a_\lambda)^m}
\,n!\,\rho_\lambda^{-n}n^{m/2-1}.
\end{align*}
In particular, these equivalents hold for every graph-tame realization
with $\lambda>0$, and their leading constants depend only on $g$ and
$\lambda$.
\end{maintheorem}
The implications from tameness are Proposition~\ref{prop:vertex-tame-Hg};
the coefficient equivalents are
Theorems~\ref{thm:motivic-alllight-asymptotics}
and~\ref{thm:motivic-ordinary-asymptotics} in
Section~\ref{sec:motivic-fixed-genus}. The analytic hypotheses constrain
only the total remainder; graph tameness is sufficient, not necessary.

Genus one is treated separately. The elliptic and rational-necklace
loci of Sections~\ref{sec:elliptic-motivic-collisions}--\ref{sec:motivic-genus-one}
give $A_1^\mu=\tfrac12\log P_\lambda+S_1^\mu$
(Proposition~\ref{prop:motivic-genus-one}). Under $(\mathrm H_1)$,
$S_1^\mu$ extends holomorphically across $\sigma_\lambda$, and
logarithmic transfer gives
\[
\mu(\overline M_{1,(0^+)^n})\sim\tfrac12 n!\,\sigma_\lambda^{-n}n^{-1},
\qquad
\mu(\overline M_{1,n})\sim\tfrac14 n!\,\rho_\lambda^{-n}n^{-1}
\]
(Theorem~\ref{thm:motivic-genus-one-asymptotics}). In genus zero there
is no all-light family. The ordinary series satisfies
\[
\mu(\overline M_{0,n})\sim
\frac{\rho_\lambda e^{\sigma_\lambda}a_\lambda}{2\sqrt\pi}
\,n!\,\rho_\lambda^{-n}n^{-5/2},
\]
without an additional analytic hypothesis. In fact, the inverse
coordinate gives a complete expansion in powers of $n^{-1}$
(Theorem~\ref{thm:motivic-genus-zero-asymptotics}).

\medskip\noindent\textbf{1.5. Hodge--Deligne realizations.}
The cohomological estimates of Section~\ref{sec:hodge-vertex-tame}
provide a family of realizations to which the preceding theory applies.
The argument combines exponential Betti bounds for vertex spaces with
the invariant-cohomology description of finite quotients, without
requiring explicit quotient Hodge polynomials. For the Hodge
consequences, let $X_n$ denote either $\overline M_{g,(0^+)^n}$ or
$\overline M_{g,n}$, and write $h^{p,r}(X_n)$ for its Hodge numbers.

\begin{maintheorem}[Hodge tameness, universality, and diagonal concentration]
\label{intro:hodge}
For every fixed $u,v\in\C$, the measure
$\mu_{u,v}(X)=E_c(X;u,v)$ is graph-tame.
If $q=uv>0$, the leading fixed-genus equivalents in both chambers,
for every $g\ge1$, depend on $u$ and $v$ only through $q$.
For each fixed $g\ge1$ and $q>0$, one also has
\[
\frac{\sum_{p\ne r}h^{p,r}(X_n)q^{(p+r)/2}}
{\sum_{p,r}h^{p,r}(X_n)q^{(p+r)/2}}\longrightarrow0
\qquad(n\to\infty).
\]
All parameters are fixed in these limits; no uniformity over pairs
$(u,v)$ or over increasing genus is asserted.
\end{maintheorem}
Graph tameness is Proposition~\ref{prop:hodge-vertex-tame}.
Universality and diagonal concentration are
Theorems~\ref{thm:fixed-q-hodge-universality}
and~\ref{thm:asymptotic-hodge-tate}, in
Sections~\ref{sec:fixed-product-universality}--\ref{sec:diagonal-concentration}.

At $q=1$, the Euler characteristic is asymptotic to the total Betti
number, and odd-degree cohomology has negligible total dimension
(Corollary~\ref{cor:euler-total-betti},
Section~\ref{sec:euler-betti-consequence}). These are statements about
relative growth; off-diagonal Hodge components need not vanish.
Purity gives $P(X_n;s)=E_c(X_n;-s,-s)$, so the fixed-product results
also apply to Poincar\'e evaluations with $q=s^2>0$.
The all-light to ordinary ratio tends to zero with explicit exponential
rates (Corollary~\ref{cor:hodge-heavy-light-ratio},
Section~\ref{sec:hodge-chamber-comparison}). In genus one we additionally
obtain closed Hodge--Deligne and Poincar\'e series
(Theorems~\ref{thm:A1explicit} and~\ref{thm:poincare-resummation},
Sections~\ref{sec:bessel}--\ref{sec:poincare1}), using the local-system
formulas of \cite[Section~4]{Kiem} and Bessel resummation.

\begin{unnumberedremark}[Related work on Hodge asymptotics]
Some of the Hodge--Deligne applications overlap with the recent
preprint of J.\,Choi and Y.-H.\,Kiem \cite{ChoiKiem}.
They obtain coefficient asymptotics for the even Betti polynomials
in both chambers, uniformly in a neighborhood of the parameter value
$1$ \cite[Theorem~4.1]{ChoiKiem}. For $g\ge2$, they also prove
concentration on the Hodge diagonal, with off-diagonal proportions of
exact orders $n^{-1/2}$ in the ordinary chamber and $n^{-1}$ in the
all-light chamber \cite[Theorem~6.5]{ChoiKiem}. In particular, this
refines the qualitative concentration at $q=1$ in
Theorem~\ref{thm:asymptotic-hodge-tate}. Their Euler-characteristic
comparison \cite[Corollary~4.6]{ChoiKiem} also agrees with the
specialization $u=v=1$ of
Corollary~\ref{cor:hodge-heavy-light-ratio}.
In this paper, the Hodge--Deligne results follow from a general framework for arbitrary
motivic measures, based on finite residual quotient families and
sufficient tameness conditions for the analytic hypotheses.
\qrem
\end{unnumberedremark}

\medskip\noindent\textbf{1.6. Structure of the paper.}
Section~\ref{sec:all-light-decomposition} sets up the motivic
conventions and the geometry of the minimal all-light chamber. It
introduces stabilization graphs, their decorated expansions, and the
parameter spaces whose finite quotients describe the coarse strata.
Section~\ref{sec:edge-propagator} computes the rational-chain series,
identifies the resulting propagator with the Losev--Manin generating
series, and records its Eulerian-polynomial coefficients.

Section~\ref{sec:motivicgraphs} separates the normalized parameter
contribution from the correction due to graph symmetries. Support
patterns determine the discrete stabilizers, and the finite residual
quotient construction extracts the unrestricted chain factors. Summing
over graphs gives the formal trivalent expansion of
Theorem~\ref{intro:formal}. Section~\ref{sec:motivic-heavy-light}
proves the single-wall identity and the heavy/light recursion in the
Grothendieck ring, derives the universal coordinate, and establishes
its connections with genus zero and the propagator, collected in
Theorem~\ref{intro:universal}.

Section~\ref{sec:g1} treats the exceptional genus-one geometry through
elliptic collision strata and rational necklaces, and derives explicit
Hodge--Deligne and Poincar\'e series. Section~\ref{sec:motivic-asymptotics}
studies the universal singularities, formulates the minimal analytic
hypotheses, and obtains the coefficient asymptotics in both chambers.
It then introduces tameness and proves the cohomological bounds needed
for Hodge--Deligne realizations. Section~\ref{sec:asymptotic-hodge-tate}
deduces fixed-product universality, asymptotic diagonal concentration,
the Euler-to-total-Betti comparison, and the heavy/light ratios.

Appendix~\ref{app:triangular} develops the coefficient form of the
universal transformation, including its triangular and derangement
recurrences and a forest interpretation.
Appendix~\ref{app:poincare-resummation} proves the Poincar\'e
resummation and gives an integral representation of the unilateral
Bessel series. Appendix~\ref{app:nonpositive-lambda} discusses zero,
negative, and complex Lefschetz values, where the dominant singularity
geometry can differ from the positive-real case.

\medskip\noindent\textbf{Notation and conventions.}
We write $[n]=\{1,\ldots,n\}$ for $n\ge1$ and $[0]=\varnothing$,
and denote by $\mathfrak S_n$ the symmetric group on $[n]$.
For a formal power series $F(z)$, the notation $[z^n]F(z)$ denotes
its coefficient of $z^n$; thus $n![z^n]F(z)$ extracts the labeled
coefficient when $F$ is written as an exponential generating series.
All markings are labeled, and graph automorphisms preserve every
light label individually. The symbols $\square$, $\triangle$, and
$\spadesuit$ mark the end of a proof, an example, and a remark,
respectively.

\medskip\noindent\textbf{Acknowledgements.}
This research was supported by Portuguese national funds through the
Funda\c{c}\~ao para a Ci\^encia e a Tecnologia (FCT),
under project UID/00208/2025.
\newline
\noindent\textbf{Use of AI tools.}
OpenAI’s ChatGPT 6.1 and Anthropic’s Claude were accessed both directly and through IAedu, the platform provided by FCT through its digital services unit FCCN, to assist with literature searches, English-language revision, the transcription of dictated passages, and the conversion of handwritten material into \LaTeX{}. The author takes full responsibility for the mathematical content, references, and final text.

\section{Stabilization graphs and the all-light decomposition}
\label{sec:all-light-decomposition}

We set up the geometric decomposition of the minimal all-light
spaces used throughout the paper. After fixing the motivic and
heavy/light conventions, we describe the curves that occur in this
chamber and the stabilization graphs that index their loci
\cite[Section~4.2]{CottiPaperI}. Decorations specify how light labels
are distributed among the surviving components and the rational
chains contracted by stabilization. We then construct the corresponding
parameter spaces and their finite coarse quotients. The chain series
and the reorganization of the quotient corrections are developed in
Sections~\ref{sec:edge-propagator} and~\ref{sec:motivicgraphs}, respectively.

\subsection{Motivic conventions}
\label{sec:motivic}

All varieties and stacks are defined over $\C$. The Grothendieck ring
$K_0(\Var_{\C})$ is generated by isomorphism classes $[X]$ of complex algebraic
varieties, subject to the \emph{scissors relations}
$[X]=[Z]+[X\setminus Z]$ for closed
subvarieties $Z\subseteq X$, and multiplication $[X][Y]=[X\times Y]$.
The class $\LL=[\mathbb A^1]$ is called the \emph{Lefschetz class}.
A motivic measure is a unital ring homomorphism
\[
\mu:K_0(\Var_{\C})\longrightarrow\mathsf{R},
\qquad \lambda:=\mu(\LL),
\]
where $\mathsf{R}$ is a commutative ring. We call $\lambda$ the
\emph{Lefschetz value} of $\mu$.

For a Deligne--Mumford stack $\mathcal X$ with coarse moduli variety $X$,
we use $[\mathcal X]_{\mathrm c}:=[X]$ and $\mu(\mathcal X):=\mu(X)$.
In particular, a finite quotient contributes the class of the quotient
variety. These are not classes in a Grothendieck ring of stacks.
Our main examples are the compactly supported Hodge--Deligne measure
$E_c(-;u,v)$, with Lefschetz value $uv$, and its Euler specialization
$E_c(-;1,1)$. The cohomological interpretation of the Poincar\'e
specialization will be recalled when we discuss Hodge realizations.

All constructions below may also be performed universally by taking
$\mu$ to be the identity of $K_0(\Var_{\C})$.

\subsection{Hassett spaces and heavy/light chambers}
\label{sec:hassett}

Let $g,n\ge0$ and let $\bm a=(a_1,\ldots,a_n)$ be an admissible
vector of rational weights, meaning
\[
0<a_i\le1,\qquad 2g-2+\sum_{i=1}^n a_i>0.
\]
An $\bm a$-stable curve is a connected projective nodal curve $C$
of arithmetic genus $g$, with labeled markings $p_1,\ldots,p_n$ in its
smooth locus, such that coincident markings have total weight at most
one and the $\Q$-line bundle
\[
\omega_C\Bigl(\sum_i a_i p_i\Bigr)
\]
is ample. On a component whose normalization has genus $h$ and $d$
branches of nodes, the latter condition is
\begin{equation}
2h-2+d+\sum_{p_i\text{ on the component}}a_i>0.
\label{eq:weighted-component-stability}
\end{equation}
A self-node contributes two branches. B.\,Hassett constructs the smooth
proper Deligne--Mumford stack $\Mbar_{g,\bm a}$ of these curves;
its coarse moduli space $\overline M_{g,\bm a}$ is projective
\cite[Theorem~2.1]{Hassett}. We retain this distinction between stack
and coarse notation throughout. When all weights are one, the markings
are distinct and one recovers the classical moduli stack $\Mbar_{g,n}$
and its coarse space $\overline M_{g,n}$, constructed by
Deligne--Mumford and F.\,F.\,Knudsen~\cite{DeligneMumford,Knudsen}.

The hyperplanes $\sum_{i\in I}a_i=1$, with $|I|\ge2$, divide the
admissible weight domain into chambers. Within a chamber, the collision
and stability conditions are unchanged, and the moduli spaces are
canonically identified. The corresponding changes when weights are
lowered will be studied through reduction morphisms in Section~\ref{sec:motivic-heavy-light}.

We write $\Mbar_{g,r|s}$ and $\overline M_{g,r|s}$ for the heavy/light
spaces with $r$ heavy markings of weight one and $s$ light markings of
equal weight $\epsilon$, where $0<\epsilon<1/s$ if $s\ge1$. Heavy
markings are pairwise distinct and cannot meet any light marking;
every subset of light markings is allowed to coincide. The notation
$(1^r,(0^+)^s)$ records this chamber. The representative $\epsilon$
is chosen separately for each $s$; the symbol $0^+$ never denotes
weight zero. For $s=0$, this is the ordinary unweighted space whenever
it is stable.

In genus zero the heavy/light space exists precisely when $r\ge3$,
or when $r=2$ and $s\ge1$. In genus one at least one marking is
required, whereas for $g\ge2$ the unmarked space is already stable.
In particular, the minimal all-light spaces are obtained by setting
$r=0$, and occur only in positive genus.

\subsection{Geometry of the minimal all-light chamber}
\label{sec:all-light-stability}

For $n\ge1$, choose $0<\epsilon<1/n$ and write
$\overline M_{g,(0^+)^n}=\overline M_{g,0|n}$.
We study the series
\begin{equation}
A_g^\mu(z):=\sum_{n\ge0}
\mu\bigl(\overline M_{g,(0^+)^n}\bigr)\frac{z^n}{n!}\in\mathsf R_{\Q}[[z]].
\label{eq:alllightseries}
\end{equation}
For $g\ge2$, the term with $n=0$ is $\mu(\overline M_g)$; for $g=1$
we set this unstable term to zero and start the geometric family at $n=1$.
There is no stable minimal all-light chamber in genus zero: the total
marking weight is less than one, whereas global stability would require
it to exceed two.

On a component whose normalization has genus $h$, with $d$ branches of
nodes and $m$ light markings, \eqref{eq:weighted-component-stability}
becomes
\begin{equation}
2h-2+d+m\epsilon>0.
\label{eq:alllight-component-stability}
\end{equation}
Since $m\epsilon<1$, a rational
component cannot have fewer than two branches of nodes. A rational
component with exactly two branches is stable precisely when $m\ge1$;
one with at least three branches is stable even when $m=0$.
Thus no rational tail occurs, and every rational two-valent component
carries a nonempty set of light labels.

For $g\ge2$, forgetting all markings and stabilizing therefore contracts
exactly the maximal chains of rational two-valent components. All other
components remain stable as unmarked components. Genus one is exceptional:
a smooth elliptic curve may be stabilized by its light markings alone,
and a cycle consisting entirely of rational two-valent components may
have arbitrary length. We treat that case separately in Section~\ref{sec:g1}.

\subsection{Stabilization graphs and decomposition by graph}
\label{sec:stabilizationgraphs}
\label{sec:decomposition-by-graph}

We use the graph conventions of \cite[Section~4.1]{CottiPaperI}.
A finite graph with legs is a quadruple $\Gamma=(V,H,\iota,\partial)$,
where $V,H$ are finite sets, $\iota$ is an involution of $H$, and
$\partial:H\to V$ is the incidence map. The elements of $H$ are
\emph{half-edges}; the two-element orbits of $\iota$ are the
\emph{edges}, and its fixed points are the \emph{legs}. A loop has
both half-edges incident to the same vertex, and multiple edges are
allowed. Legs record markings in the dual graph of a pointed curve.
For a graph without legs, the involution $\iota$ is fixed-point-free,
and we write
\[
V(\Gamma)=|V|,\qquad E(\Gamma)=|H|/2,\qquad
d_v=|\partial^{-1}(v)|.
\]
A genus decoration assigns an integer $h_v\ge0$ to each vertex.
For a connected graph its genus is
\[
g(\Gamma)=\sum_vh_v+E(\Gamma)-V(\Gamma)+1.
\]
An automorphism preserves incidence, the involution, and the genus
decoration. In particular, reversing a loop exchanges its two half-edges
and is part of the automorphism data.

\begin{definition}
A \emph{stabilization graph of genus $g\ge2$} is a connected genus-decorated graph
without legs, of genus $g$, such that
\begin{equation}
2h_v-2+d_v>0\qquad\text{at every vertex }v.
\label{eq:graphstability}
\end{equation}
The stabilization graph of an all-light curve is obtained from its dual graph by
forgetting all marking legs and contracting each maximal chain of
genus-zero two-valent vertices.
\end{definition}

Equivalently, the graph is the dual graph of the stable unmarked curve
obtained by forgetting the markings. It is a stable graph of type $(g,0)$
in the usual terminology~\cite{GetzlerKapranov}. The automorphism group
$\Aut\Gamma$ has no marking labels to preserve: these reappear in the
decorations considered below.

\begin{lemma}[The maximum number of graph edges]
\label{lem:maxedges}
For a stabilization graph of genus $g\ge2$,
\[
E(\Gamma)\le3g-3.
\]
Equality holds precisely when all vertices have genus zero and valence
three. In that case $V(\Gamma)=2g-2$.
\end{lemma}

\begin{proof}
The genus formula and $\sum_vd_v=2E(\Gamma)$ give
\[
\sum_v(2h_v-2+d_v)=2g-2.
\]
Each summand is a positive integer, so $V(\Gamma)\le2g-2$. Hence
\[
E(\Gamma)=g-1+V(\Gamma)-\sum_vh_v\le3g-3.
\]
Equality forces $\sum_vh_v=0$ and $V(\Gamma)=2g-2$, so every summand
above equals one and every vertex has valence three. The converse follows
from the same identities.
\end{proof}

The bounds on vertices and edges, together with $\sum_vh_v\le g$, show
that there are finitely many isomorphism classes of stabilization graphs in a
fixed genus. Let $\mathfrak C_g$ be a set of representatives and let
$\mathfrak T_g\subseteq\mathfrak C_g$ consist of the rational trivalent
graphs. We shall use the finite symmetry sum
\begin{equation}
\gamma_g:=\sum_{\Gamma\in\mathfrak T_g}\frac1{|\Aut\Gamma|}.
\label{eq:trivalent-constant-definition}
\end{equation}
Its enumeration and large-genus behavior are studied in
\cite[Section~7]{CottiPaperI}.

\smallskip\noindent\emph{Decomposition by stabilization graph.}

Fix $g\ge2$. Forgetting the light markings and stabilizing gives a
morphism $f_n:\overline M_{g,(0^+)^n}\to\overline M_g$.
As in \cite[Section~4.2]{CottiPaperI}, write
$\stratumZ_\Gamma\subseteq\overline M_g$ for the stratum with dual
graph $\Gamma\in\mathfrak C_g$. The inverse image
$f_n^{-1}(\stratumZ_\Gamma)$ is the locally closed coarse locus of
all-light curves whose stabilization graph is isomorphic to $\Gamma$.
The resulting finite decomposition is
\[
\overline M_{g,(0^+)^n}
=\bigsqcup_{\Gamma\in\mathfrak C_g}f_n^{-1}(\stratumZ_\Gamma).
\]
Define
\[
A_{g,\Gamma}^\mu(z)
:=\sum_{n\ge0}\mu\bigl(f_n^{-1}(\stratumZ_\Gamma)\bigr)
\frac{z^n}{n!}.
\]
Additivity gives the exact identity
\begin{equation}
A_g^\mu(z)=\sum_{\Gamma\in\mathfrak C_g}A_{g,\Gamma}^\mu(z).
\label{eq:finite-graph-sum}
\end{equation}
The index set does not grow with $n$. What grows above a fixed graph is
the number of light markings and the lengths of the chains inserted
along its edges. We next isolate these variable parts.

\subsection{Decorations and geometric building blocks}
\label{sec:decorated-expansions}

Fix a representative stabilization graph
$\Gamma=(V,H,\iota,\partial)$. Its vertices and half-edges will index
the component factors and their attaching markings in the constructions
below.
A \emph{discrete decoration} with label set $[n]=\{1,\ldots,n\}$
is a triple $\delta=(\widetilde\Gamma,c,\ell)$. Here
$\widetilde\Gamma$ is obtained by subdividing the edges of $\Gamma$
with genus-zero two-valent vertices, and
\[
c:\widetilde\Gamma^{\mathrm{st}}\xrightarrow{\sim}\Gamma
\]
is a genus-preserving graph isomorphism, where
$\widetilde\Gamma^{\mathrm{st}}$ is obtained by suppressing the inserted
vertices. The assignment
\[
\ell:[n]\longrightarrow\Vertices(\widetilde\Gamma)
\]
is required to give every inserted vertex at least one label; graph
vertices may receive none. The isomorphism $c$ identifies the surviving
vertices and the end half-edges of each subdividing chain with the
vertices and half-edges of the fixed graph, including the two branches
of a loop. Thus suppression and identification with $\Gamma$ are
separate operations, as in \cite[Section~4.3]{CottiPaperI}.

Two decorations $(\widetilde\Gamma,c,\ell)$ and
$(\widetilde\Gamma',c',\ell')$ are identified if there is a
genus-preserving graph isomorphism
$\varphi:\widetilde\Gamma\to\widetilde\Gamma'$ such that
$c'\circ\varphi^{\mathrm{st}}=c$ and
$\varphi_V\circ\ell=\ell'$. Their finite set of isomorphism classes
is denoted $\mathscr D_n(\Gamma)$; at most $n$ inserted vertices can
occur. Write $B_w=\ell^{-1}(w)$.
\footnote{We record marking labels by the assignment $\ell$ rather than
by legs attached to the expanded graph. These decorated expansions have only
the half-edges corresponding to nodal branches; the positions and
collisions of markings are encoded by the component moduli spaces.}

The decoration records which component carries each label, but not its
position or its collisions with other markings on that component. Those
continuous moduli are supplied by the following two building blocks.

\smallskip\noindent\emph{Vertex spaces.}
For a stable pair $(h,d)$, meaning $2h-2+d>0$, let
$\mathscr V_{h,S,B}$ be the moduli stack of smooth genus-$h$ curves with
distinct heavy markings indexed by a set $S$ of cardinality $d$, and
light markings indexed by $B$. Light markings may coincide with one
another but must avoid the heavy markings. Thus $\mathscr V_{h,S,B}$
is the smooth-curve locus of the heavy/light Hassett stack with weights
$(1^S,(0^+)^B)$. Let $\Vspace_{h,S,B}$ be its
coarse moduli space. Up to relabeling, write $\Vspace_{h,d,m}$ when
$|B|=m$, and define
\begin{equation}
V_{h,d}^\mu(z):=\sum_{m\ge0}\mu(\Vspace_{h,d,m})\frac{z^m}{m!}.
\label{eq:vertexseries}
\end{equation}
At a graph vertex $v$, the set $S_v$ of incident half-edges indexes the
heavy markings that will become branches of nodes.

\begin{example}[The rational trivalent vertex]
For $(h,d)=(0,3)$, normalize the heavy markings to $0,1,\infty$.
The $m$ light markings vary independently in
$\mathbb P^1\setminus\{0,1,\infty\}$, including when some of them
coincide. Consequently,
\[
\Vspace_{0,3,m}\cong
\bigl(\mathbb P^1\setminus\{0,1,\infty\}\bigr)^m,
\qquad
\mu(\Vspace_{0,3,m})=(\lambda+1-3)^m=(\lambda-2)^m,
\]
and the vertex series is
\[
V_{0,3}^\mu(z)=\sum_{m\ge0}(\lambda-2)^m\frac{z^m}{m!}
=e^{(\lambda-2)z}.
\]
\qetr
\end{example}

\smallskip\noindent\emph{Chain spaces.}
An inserted component has genus zero, two distinct distinguished nodal
branches, and a nonempty light label set $B$. Denote its moduli space by
$\mathcal C_B$. Normalizing the branches to $0,\infty$ and one selected
light marking to $1$ identifies
\[
\mathcal C_B\cong\GG^{|B|-1}.
\]
The selection is only a choice of coordinates. In particular, an
inserted component has no nontrivial automorphism fixing its nodal
branches and every light marking. A chain is an ordered sequence of
these components, ordered from one distinguished end to the other.

Section~\ref{sec:edge-propagator} assembles these component spaces into spaces of chains and
introduces their generating series, the edge propagator. For now, we
retain the individual component factors when constructing the space
attached to a decoration.

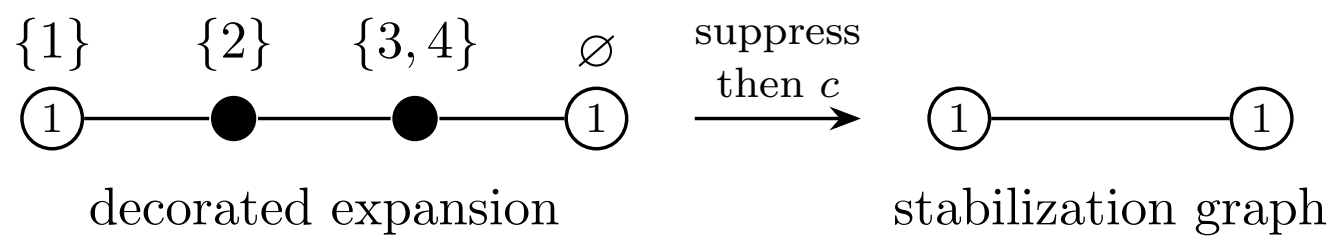
\begin{figure}[htbp]
\centering
\begin{tikzpicture}[>=Stealth,line width=.65pt,every node/.style={font=\small}]
\tikzset{graph/.style={circle,draw,line width=.65pt,minimum size=4mm,inner sep=0pt,fill=white,font=\scriptsize},
 inserted/.style={circle,fill,minimum size=3mm,inner sep=0pt}}
\node[graph] (a) at (0,0) {$1$};
\node[inserted] (u) at (1.2,0) {};
\node[inserted] (w) at (2.4,0) {};
\node[graph] (b) at (3.6,0) {$1$};
\draw (a)--(u)--(w)--(b);
\node[above=5pt] at (a) {$\{1\}$};
\node[above=5pt] at (u) {$\{2\}$};
\node[above=5pt] at (w) {$\{3,4\}$};
\node[above=5pt] at (b) {$\varnothing$};
\draw[->] (4.25,0)--(5.35,0) node[midway,above,align=center,font=\scriptsize] {suppress\\then $c$};
\node[graph] (aa) at (6,0) {$1$};
\node[graph] (bb) at (8,0) {$1$};
\draw (aa)--(bb);
\node[below=9pt] at (1.8,0) {decorated expansion};
\node[below=9pt] at (7,0) {stabilization graph};
\end{tikzpicture}
\caption{A genus-two example. The inserted rational vertices carry
nonempty label sets and are contracted to the edge of the graph. Their
component spaces are respectively a point and $\mathbb G_m$; the two
genus-one graph vertices contribute vertex spaces with one heavy marking.}
\label{fig:graph-strategy}
\end{figure}

\subsection{Parameter contributions and graph symmetries}
\label{sec:parameter-and-quotient}

For $\delta=(\widetilde\Gamma,c,\ell)$, the parameter stack is
\[
\mathscr X_\delta=
\prod_{v\in\Vertices(\Gamma)}\mathscr V_{h_v,S_v,B_v}
\times
\prod_{u\in\Vertices(\widetilde\Gamma)\setminus\Vertices(\Gamma)}
\mathcal C_{B_u}.
\]
Its coarse space is the corresponding product
\begin{equation}
X_\delta=
\prod_{v\in\Vertices(\Gamma)}\Vspace_{h_v,S_v,B_v}
\times
\prod_{u\in\Vertices(\widetilde\Gamma)\setminus\Vertices(\Gamma)}
\mathcal C_{B_u}.
\label{eq:decoration-parameter-space}
\end{equation}
Gluing the markings indexed by paired half-edges reconstructs the nodal
curve. No graph symmetries have yet been divided out.

The group $\Aut\Gamma$ acts on $\mathscr D_n(\Gamma)$ by changing
the identification with the fixed graph: for $\alpha\in\Aut\Gamma$,
\[
\alpha\cdot(\widetilde\Gamma,c,\ell)
=(\widetilde\Gamma,\alpha\circ c,\ell).
\]
In particular, the light labels remain fixed. The stabilizer
$G_\delta:=\operatorname{Stab}_{\Aut\Gamma}(\delta)$ acts on
$\mathscr X_\delta$ and on $X_\delta$ by permuting the component
factors and relabeling their attaching points; on an inserted component,
it may exchange the two nodal branches. The discrete
stabilizer $G_\delta$ is distinct from the automorphism group of an
individual curve, whose component automorphisms are already encoded
by the parameter stack.

Normalization at the nodes identifies the stack stratum indexed by the
orbit $[\delta]$ with $[\mathscr X_\delta/G_\delta]$. The coarse moduli
map is equivariant, and in characteristic zero we have
\[
\operatorname{coarse}\bigl([\mathscr X_\delta/G_\delta]\bigr)
\cong X_\delta/G_\delta;
\]
see \cite[Theorem~3.2 and Corollary~3.3(a)]{AOV08} and
\cite[Section~4.3, the coarse-quotient identification following
(4.4)]{CottiPaperI}. The quotient on the right is the quotient variety,
whose class is the one evaluated by $\mu$.
Consequently,
\begin{equation}
A_{g,\Gamma}^\mu(z)
=\sum_{n\ge0}\frac{z^n}{n!}
\sum_{[\delta]\in\mathscr D_n(\Gamma)/\Aut\Gamma}
\mu(X_\delta/G_\delta).
\label{eq:actual-graph-contribution}
\end{equation}

Before passing to graph-automorphism orbits, define the parameter series
by summing over all decorations of the fixed representative graph:
\begin{equation}
A_{\Gamma,\mathrm{par}}^\mu(z)
:=\sum_{n\ge0}\frac{z^n}{n!}
\sum_{\delta\in\mathscr D_n(\Gamma)}\mu(X_\delta).
\label{eq:parameter-graph-series}
\end{equation}
This sum includes decorations with nontrivial stabilizers. The independent
distribution of labels among vertices and edge chains will give its
product expression once the chain series has been introduced in Section~\ref{sec:edge-propagator}.

To compare \eqref{eq:actual-graph-contribution} with the normalized
parameter series, consider one orbit of decorations represented by
$\delta$. The orbit--stabilizer theorem gives
\[
|\Aut\Gamma\cdot\delta|=\frac{|\Aut\Gamma|}{|G_\delta|}.
\]
All parameter spaces indexed by this orbit are isomorphic, so its
contribution to the coefficient of $z^n/n!$ in
$A_{\Gamma,\mathrm{par}}^\mu/|\Aut\Gamma|$ is
\[
\frac1{|\Aut\Gamma|}
\sum_{\delta'\in\Aut\Gamma\cdot\delta}\mu(X_{\delta'})
=\frac{\mu(X_\delta)}{|G_\delta|},
\]
whereas its actual contribution is $\mu(X_\delta/G_\delta)$.
Subtracting these contributions and summing over the orbits gives
\begin{equation}
A_{g,\Gamma}^\mu(z)
=\frac1{|\Aut\Gamma|}A_{\Gamma,\mathrm{par}}^\mu(z)
+\mathcal R_\Gamma^\mu(z),
\label{eq:graphorbitdecomposition}
\end{equation}
where
\begin{equation}
n![z^n]\mathcal R_\Gamma^\mu(z)
=\sum_{\substack{[\delta]\in\mathscr D_n(\Gamma)/\Aut\Gamma\\
G_\delta\ne\{1\}}}
\left(\mu(X_\delta/G_\delta)
-\frac{\mu(X_\delta)}{|G_\delta|}\right).
\label{eq:nonfreecorrection}
\end{equation}
For a general motivic measure, the quotient value in this expression
cannot be replaced by formal division by the group order, or by an
average of the measures of fixed loci. Section~\ref{sec:motivicgraphs} retains these coarse
quotients and separates a finite residual geometry from the unrestricted
chain factors. Analytic bounds on the residual series will subsequently
allow the maximal trivalent term to determine the leading asymptotics.

\section{Rational chains and the Losev--Manin propagator}
\label{sec:edge-propagator}

We compute the series of the rational chains introduced in
Section~\ref{sec:decorated-expansions}, with their two ends distinguished.
Starting with a single component, we sum over ordered chains of all
lengths to obtain the edge propagator. We then identify this series
with the series of Losev--Manin spaces and describe its coefficients
through Eulerian polynomials.

The chain formulas are first defined with an indeterminate $\lambda$
over $\Q[\lambda]$ and then specialized to the Lefschetz value
$\mu(\LL)$. In particular, the universal class formulas are obtained
by setting $\lambda=\LL$.

\subsection{A single inserted component}
\label{sec:single-chain-component}

Let $B$ be a nonempty light label set, with $|B|=m$. As in Section~\ref{sec:all-light-decomposition},
$\mathcal C_B$ parametrizes a smooth rational component with two
distinct distinguished end points and light markings indexed by $B$.
Light markings lie away from the end points and may coincide with one
another. Normalizing the end points to $0,\infty$ and one chosen light
marking to $1$ gives
\[
\mathcal C_B\cong\GG^{m-1},\qquad
[\mathcal C_B]=(\LL-1)^{m-1}.
\]
Thus the exponential generating series of one inserted component is
\begin{equation}
R_\lambda(z)
:=\sum_{m\ge1}(\lambda-1)^{m-1}\frac{z^m}{m!}.
\label{eq:Rlambda}
\end{equation}
It has no constant term, since a two-valent rational component without
light markings is unstable. Its linear coefficient is one: one light
marking removes the remaining scaling automorphism and gives a point.

The series also admits the expression
\begin{equation}
R_\lambda(z)=\frac{e^{(\lambda-1)z}-1}{\lambda-1}.
\label{eq:Rlambda-closed}
\end{equation}
This quotient denotes the coefficientwise polynomial series
\eqref{eq:Rlambda}; it requires no inversion of $\lambda-1$ in $\mathsf{R}$.
In particular, $R_1(z)=z$.

\subsection{Ordered chains and the edge propagator}
\label{sec:motivicpropagator}

For a nonempty finite set $B$, let $\mathcal P_B$ be the disjoint union
\begin{equation}
\mathcal P_B
=\coprod_{k=1}^{|B|}
\ \coprod_{\substack{B=B_1\sqcup\cdots\sqcup B_k\\B_i\ne\varnothing}}
\prod_{i=1}^k\mathcal C_{B_i},
\label{eq:chain-parameter-union}
\end{equation}
where the inner sum is over ordered partitions. Set
$\mathcal P_\varnothing=\operatorname{pt}$, representing a graph edge
without inserted components. These spaces assemble the component
factors $\mathcal C_B$ introduced in Section~\ref{sec:decorated-expansions}. The components of a
chain are ordered from one
distinguished end to the other, so there is no division by $k!$.
We call their exponential generating series
\begin{equation}
P_\lambda(z):=\sum_{m\ge0}\mu(\mathcal P_{[m]})\frac{z^m}{m!}
\label{eq:Plambda-geometric}
\end{equation}
the \emph{edge propagator}. Since the component factors are products
of algebraic tori, its coefficients depend on $\mu$ only through
the Lefschetz value $\lambda$.

\begin{proposition}[The chain series]
The edge propagator satisfies
\begin{equation}
P_\lambda(z)=\sum_{k\ge0}R_\lambda(z)^k
=\frac1{1-R_\lambda(z)}.
\label{eq:Plambda}
\end{equation}
\end{proposition}

\begin{proof}
By additivity, multiplicativity, and
\eqref{eq:chain-parameter-union},
\begin{align*}
P_\lambda(z)
&=1+\sum_{m\ge1}\frac{z^m}{m!}
\sum_{k=1}^m
\sum_{\substack{[m]=B_1\sqcup\cdots\sqcup B_k\\B_i\ne\varnothing}}
\prod_{i=1}^k\mu(\mathcal C_{B_i})\\
&=1+\sum_{k\ge1}\sum_{m_1,\ldots,m_k\ge1}
\prod_{i=1}^k
\left((\lambda-1)^{m_i-1}\frac{z^{m_i}}{m_i!}\right)\\
&=\sum_{k\ge0}
\left(\sum_{m\ge1}(\lambda-1)^{m-1}\frac{z^m}{m!}\right)^k
=\sum_{k\ge0}R_\lambda(z)^k
=\frac1{1-R_\lambda(z)}.
\end{align*}
\end{proof}

Substitution of \eqref{eq:Rlambda-closed} yields the shorthand
\begin{equation}
P_\lambda(z)=\frac{\lambda-1}{\lambda-e^{(\lambda-1)z}}.
\label{eq:PlambdaClosed}
\end{equation}
The defining expression is \eqref{eq:Plambda}; the apparent singularity
of \eqref{eq:PlambdaClosed} at $\lambda=1$ is removable coefficientwise.
In particular,
\begin{equation}
P_1(z)=\frac1{1-z},\qquad P_0(z)=e^z.
\label{eq:propagator-basic-values}
\end{equation}

\begin{remark}[Reversal of the ends]
Exchanging the two ends reverses the order of the blocks in
\eqref{eq:chain-parameter-union} and exchanges the two distinguished
points on every component. This is an involution of the parameter
construction, not an identification made in $P_\lambda$.
For example, on a one-component chain, the normalized light coordinates
are sent to their inverses. The graph automorphism group may induce
such reversals, including when the graph edge is a loop. Its action is
retained when forming the actual coarse quotients in Section~\ref{sec:motivicgraphs}.
\qrem
\end{remark}

\subsection{Losev--Manin spaces}
\label{sec:losevmaninpropagator}

For $n\ge1$, let
\[
\overline L_n:=\overline M_{0,2|n}
\]
be the coarse Losev--Manin space, with two labeled heavy markings
$x_-,x_+$ and $n$ labeled light markings in the minimal positive-weight
chamber~\cite{LosevManin,Hassett}. Its curves are chains of rational
components from $x_-$ to $x_+$, each component carrying at least one
light marking. See also \cite[Section~7.1]{CottiRHB} for a discussion
of these spaces and their role in Losev--Manin cohomological field theories.

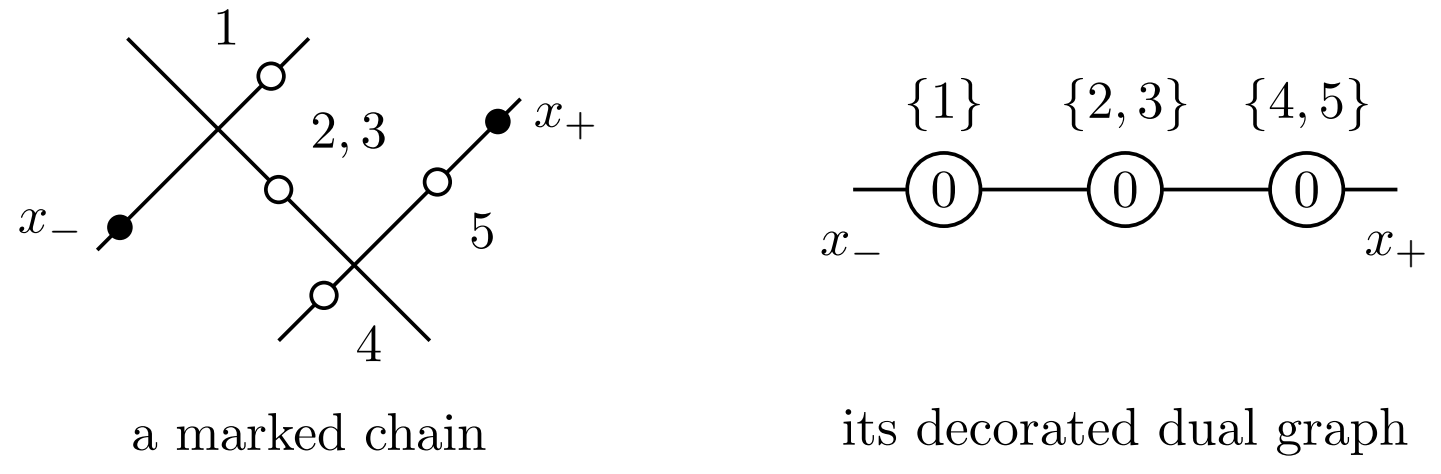
\begin{figure}[htbp]
\centering
\begin{tikzpicture}[line width=.7pt,every node/.style={font=\small},
 heavy/.style={circle,fill=black,inner sep=1.7pt},
 light/.style={circle,draw,fill=white,inner sep=1.7pt},
 dual/.style={circle,draw,fill=white,inner sep=2pt}]
\draw (0,-.6)--(1.4,.8);
\draw (.2,.8)--(2.2,-1.2);
\draw (1.2,-1.2)--(2.8,.4);
\node[heavy,label=left:$x_-$] at (.15,-.45) {};
\node[heavy,label=right:$x_+$] at (2.65,.25) {};
\node[light,label=above left:$1$] at (1.15,.55) {};
\node[light,label=above right:{$2,3$}] at (1.2,-.2) {};
\node[light,label=below right:$4$] at (1.5,-.9) {};
\node[light,label=below right:$5$] at (2.25,-.15) {};
\node at (1.4,-1.8) {a marked chain};
\draw (5,-.2)--(8.6,-.2);
\node[dual] (a) at (5.6,-.2) {$0$};
\node[dual] (b) at (6.8,-.2) {$0$};
\node[dual] (c) at (8,-.2) {$0$};
\node[above=7pt] at (a) {$\{1\}$};
\node[above=7pt] at (b) {$\{2,3\}$};
\node[above=7pt] at (c) {$\{4,5\}$};
\node[below=4pt] at (5,-.2) {$x_-$};
\node[below=4pt] at (8.6,-.2) {$x_+$};
\node at (6.8,-1.8) {its decorated dual graph};
\end{tikzpicture}
\caption{A point of $\overline L_5$. The segments represent three
copies of $\mathbb P^1$, meeting at two nodes. Black points are the
heavy markings and white points are the light markings. Labels $2$
and $3$ coincide, whereas $4$ and $5$ are distinct points on the same
component. On the dual graph, the vertices have
genus zero, the two legs record the heavy ends, and the sets above
the vertices record the light-label decoration. As in
Section~\ref{sec:decorated-expansions}, these sets record the component
carrying each label, without specifying coincidences.}
\label{fig:losev-manin-chain}
\end{figure}

Figure~\ref{fig:losev-manin-chain} shows a chain and the label
allocation on its dual graph.

Indeed, the dual graph is a tree. Each terminal component must carry
a heavy marking, since a terminal component with only light markings
would be unstable. The tree therefore has at most two leaves and is
a path, with the heavy markings at its ends. For a single component,
the two heavy markings play the same role. On every component the
number of nodal branches plus heavy markings is two; stability
requires at least one light marking.

The ordering from $x_-$ to $x_+$ determines an ordered partition
$[n]=B_1\sqcup\cdots\sqcup B_k$ into nonempty blocks. For a fixed
partition, the corresponding locally closed stratum is
\[
\prod_{i=1}^k\mathcal C_{B_i}\cong\GG^{n-k}.
\]
There is no graph quotient here: the two heavy ends are labeled and
fix the order of every component of the path. Thus the chain spaces
in \eqref{eq:chain-parameter-union} describe the strata of
$\overline L_n$.

\begin{proposition}[Losev--Manin interpretation]
\label{prop:LM-interpretation}
For every $n\ge1$,
\begin{equation}
[\overline L_n]
=\sum_{k=1}^n k!\,\genfrac{\{}{\}}{0pt}{}{n}{k}(\LL-1)^{n-k}
=[\mathcal P_{[n]}]
\quad\text{in }K_0(\Var_{\C}),
\label{eq:LM-strata-class}
\end{equation}
where $\genfrac{\{}{\}}{0pt}{}{n}{k}$ denotes the Stirling number of
the second kind, the number of partitions of $[n]$ into $k$ nonempty
blocks. Consequently,
\begin{equation}
P_\lambda(z)=1+\sum_{n\ge1}\mu(\overline L_n)\frac{z^n}{n!}.
\label{eq:LMpropagator}
\end{equation}
\end{proposition}

\begin{proof}
There are $k!\genfrac{\{}{\}}{0pt}{}{n}{k}$ ordered partitions into $k$ nonempty blocks, and
each corresponding stratum has class $(\LL-1)^{n-k}$. Additivity
gives \eqref{eq:LM-strata-class}; applying $\mu$ and summing proves
\eqref{eq:LMpropagator}.
\end{proof}

The equality in \eqref{eq:LM-strata-class} is an equality of classes.
The disjoint union $\mathcal P_{[n]}$ keeps the strata separate, whereas
in the proper space $\overline L_n$ they are glued along their
degenerations. We do not identify these two varieties. Also, the
constant term in \eqref{eq:LMpropagator} represents an empty chain;
there is no stable two-heavy, zero-light space $\overline M_{0,2|0}$.

\subsection{Eulerian polynomials and specializations}
\label{sec:eulerian-propagator}

The Losev--Manin space is the smooth projective toric variety associated
with the permutohedron~\cite{LosevManin}. Its $h$-polynomial\footnote{For
a simple $d$-dimensional polytope $P$, let $f_i(P)$ be the number of
$i$-dimensional faces, including $f_d(P)=1$. Its $h$-polynomial is
$h_P(t)=\sum_{i=0}^d f_i(P)(t-1)^i$; see
\cite[Section~2.2]{PostnikovReinerWilliams}. For the
$(n-1)$-dimensional permutohedron this is the Eulerian polynomial;
see \cite[Section~4.2, equation~(4.2)]{PostnikovReinerWilliams}.
For the cohomological interpretation on the associated smooth
projective toric variety, see \cite[Section~5.2]{FultonToric}.}
is the Eulerian polynomial
\[
\mathrm A_n(t):=\sum_{\sigma\in\mathfrak S_n}t^{\des(\sigma)},
\qquad n\ge1,
\]
where $\des(\sigma)=|\{i<n:\sigma(i)>\sigma(i+1)\}|$; see
\cite[\S1.4]{StanleyEC1}.
A generic one-parameter subgroup of the torus gives an affine-cell
decomposition, with cell counts encoded by the $h$-polynomial
\cite{BialynickiBirula,FultonToric}. Consequently,
\begin{equation}
[\overline L_n]=\mathrm A_n(\LL),\qquad
\mu(\overline L_n)=\mathrm A_n(\lambda).
\label{eq:LM-Eulerian}
\end{equation}
Combining the two descriptions gives
\begin{equation}
\mathrm A_n(\lambda)
=\sum_{k=1}^n k!\genfrac{\{}{\}}{0pt}{}{n}{k}(\lambda-1)^{n-k},
\qquad
P_\lambda(z)=1+\sum_{n\ge1}\mathrm A_n(\lambda)\frac{z^n}{n!}.
\label{eq:Eulerian-chain-expansion}
\end{equation}
The first terms are
\[
P_\lambda(z)=1+z+(1+\lambda)\frac{z^2}{2!}
+(1+4\lambda+\lambda^2)\frac{z^3}{3!}
+(1+11\lambda+11\lambda^2+\lambda^3)\frac{z^4}{4!}+O(z^5).
\]
Here $O(z^5)$ denotes terms of degree at least five in the formal
variable $z$. Since these spaces are smooth projective toric varieties,
their Hodge--Deligne and ordinary Poincar\'e polynomials are respectively
\[
E_c(\overline L_n;u,v)=\mathrm A_n(uv),\qquad
P(\overline L_n;t)=\mathrm A_n(t^2).
\]
In particular, their cohomology is concentrated in even degrees and
$\chi(\overline L_n)=\sum_j b_j(\overline L_n)=n!$, in agreement
with $P_1(z)=(1-z)^{-1}$.

We can now complete the product formula for the parameter series
introduced in Section~\ref{sec:parameter-and-quotient}. For a fixed representative graph $\Gamma$,
before passing to graph-automorphism orbits, label sets are distributed
independently among its vertices and edge chains. Multiplicativity
and the exponential product rule give
\begin{equation}
A_{\Gamma,\mathrm{par}}^\mu(z)
=\left(\prod_{v\in\operatorname{Vert}(\Gamma)}
V_{h_v,d_v}^\mu(z)\right)P_\lambda(z)^{E(\Gamma)}.
\label{eq:parameter-motivic-graph}
\end{equation}
An orientation may be chosen on each graph edge to describe the order
of its chain; reversing this choice gives an isomorphic parameter
space. There is one factor $P_\lambda$ for every edge, including loops.
For a rational trivalent graph of genus $g$, this gives the normalized
parameter contribution
\[
\frac{e^{2(g-1)(\lambda-2)z}}{|\Aut\Gamma|}
P_\lambda(z)^{3g-3}.
\]
Here $V_{0,3}^\mu=e^{(\lambda-2)z}$ is the vertex factor computed in
Section~\ref{sec:parameter-and-quotient}. The formula for the parameter series imposes no condition
that decoration stabilizers be trivial. The next section determines
how graph symmetries alter this contribution and isolates the term
supplied by rational trivalent graphs.

\section{Motivic contributions of stabilization graphs}
\label{sec:motivicgraphs}

The product of vertex series and edge propagators computed in the
preceding sections records the parameter contribution of a fixed
stabilization graph. To pass from this product to the contribution of
the corresponding coarse locus, we must account for the finite graph
symmetries and their actions on the component moduli. Our aim is to
express the resulting correction in a form that preserves the quotient
geometry while separating off the unrestricted chain factors.

For the Euler characteristic, the Burnside--Lefschetz argument of
\cite[Section~4.4]{CottiPaperI} expresses this correction through fixed
loci. Its vertex calculations use the fact that a finite covering has
Euler characteristic equal to its degree times that of the base,
even when monodromy permutes the points of the fibers. For general
motivic measures, however, the fixed-locus average and the covering
degree formula can both fail. The construction below therefore uses
the coarse quotient varieties themselves and keeps track of the
actions that define them.

We organize the correction into finitely many residual quotient
families multiplied by powers of $P_\lambda$, and show that each term
has fewer unrestricted propagators than the number of edges of the
fixed graph. When we sum over stabilization graphs, this separation
isolates the contribution of the rational trivalent graphs as a
universal leading term. The analytic conditions under which it governs
the coefficient asymptotics will be established in Section~\ref{sec:motivic-asymptotics}.

Throughout this section $g\ge2$, and all identities are formal.

\subsection{The correction from finite coarse quotients}
\label{sec:motivicgraphcontribution}

Fix a representative stabilization graph $\Gamma$.
We retain the notation of Section~\ref{sec:all-light-decomposition}: $\mathscr D_n(\Gamma)$ is the set
of decorations with light labels $[n]$, $X_\delta$ is the coarse
parameter space, and $G_\delta$ is the stabilizer of the discrete
decoration. The actual contribution and the parameter contribution before passing to
graph-automorphism orbits are
\begin{align*}
A_{g,\Gamma}^\mu(z)
&=\sum_{n\ge0}\frac{z^n}{n!}
\sum_{[\delta]\in\mathscr D_n(\Gamma)/\Aut\Gamma}\mu(X_\delta/G_\delta),\\
A_{\Gamma,\mathrm{par}}^\mu(z)
&=\left(\prod_{v\in\Vertices(\Gamma)}V_{h_v,d_v}^\mu(z)\right)
P_\lambda(z)^{E(\Gamma)}.
\end{align*}
By the orbit computation in Section~\ref{sec:parameter-and-quotient},
\begin{equation}
A_{g,\Gamma}^\mu=\frac1{|\Aut\Gamma|}A_{\Gamma,\mathrm{par}}^\mu
+\mathcal R_\Gamma^\mu.
\label{eq:section4-graph-correction}
\end{equation}
For the present argument it is useful to sum over decorations rather
than their orbits:
\begin{equation}
n![z^n]\mathcal R_\Gamma^\mu
=\frac1{|\Aut\Gamma|}\sum_{\substack{\delta\in\mathscr D_n(\Gamma)\\G_\delta\ne\{1\}}}
\left(|G_\delta|\mu(X_\delta/G_\delta)-\mu(X_\delta)\right).
\label{eq:all-decorations-correction}
\end{equation}
Indeed, an orbit has $|\Aut\Gamma|/|G_\delta|$ elements, so its contribution
to the right-hand side is exactly
$\mu(X_\delta/G_\delta)-\mu(X_\delta)/|G_\delta|$.

The quotient value in this formula remains a coarse-space invariant.
In general it is neither division of $\mu(X_\delta)$ by the group
order nor an average of the measures of fixed loci. Our task is to
separate the chain factors on which the prescribed stabilizer acts
trivially from the geometry on which its action must be retained.

\subsection{Support patterns and exact decoration stabilizers}
\label{sec:motivic-residual-types}

Which data of a decoration $\delta$ determine its stabilizer
$G_\delta$ and tell us on which factors it acts trivially? Each light
label is fixed individually, so a vertex carrying a light label must
be fixed. On an edge with inserted components, the labels force the
edge to be preserved. If there is just one inserted component, its
labels do not distinguish the two ends of the edge. If there are at
least two components, the terminal components carry disjoint nonempty
sets of labels, so reversing the chain would move those labels.
An empty edge imposes no such restriction. These are precisely the restrictions recorded by the following definition.

A \emph{support pattern} on the fixed representative graph is a pair
$p=(S_p,a_p)$, where $S_p\subseteq\Vertices(\Gamma)$ and
$a_p:\Edge(\Gamma)\to\{0,1,2\}$. For a decoration, $S_p$ is the set
of graph vertices carrying light labels and
\[
a_p(e)=
\begin{cases}
0,&\text{the chain on $e$ is empty},\\
1,&\text{the chain has one inserted component},\\
2,&\text{the chain has at least two inserted components}.
\end{cases}
\]
Put $L_p=\{e:a_p(e)=2\}$ and $m_p=|L_p|$. Define $H_p\subseteq \Aut\Gamma$
to consist of automorphisms that fix every vertex in $S_p$, preserve
each edge with $a_p(e)=1$ setwise, and fix both half-edges individually
on every edge in $L_p$.

\begin{lemma}[Exact decoration stabilizer]
\label{lem:exact-support-stabilizer}
Every decoration with support pattern $p$ has $G_\delta=H_p$.
\end{lemma}
\begin{proof}
First, we prove $G_\delta\subseteq H_p$. Every element of $G_\delta$
fixes all light labels individually. Thus an occupied graph vertex
must be fixed, and an occupied edge cannot move to another edge.
A one-component chain can be reversed without changing its discrete
label assignment. A chain with at least two components cannot be
reversed: the two terminal inserted components carry disjoint nonempty
label sets. Its two end half-edges must therefore be fixed individually.
These conditions are necessary, including for loops.

Second, we prove $H_p\subseteq G_\delta$. An element of $H_p$ lifts
along each subdivided edge.
On a long edge the lift fixes each component in its order; on an edge
with one inserted component it fixes that vertex, possibly exchanging
its two nodal branches. It fixes the labels on occupied graph vertices,
and empty edges and unoccupied vertices impose no further label
conditions. Hence the lift preserves the entire discrete decoration.
\end{proof}

Here $H_p$ is the stabilizer of a \emph{labeled decoration}, not in
general the setwise stabilizer of its support pattern. An automorphism
may preserve the pattern while permuting occupied positions carrying
different labels. All sums over $p$ below run over patterns on the
fixed representative graph, not over pattern orbits.

For $H_p\ne\{1\}$ we have
\begin{equation}
m_p\le E(\Gamma)-1.
\label{eq:residual-propagator-bound}
\end{equation}
If all edges were long, an element of $H_p$ would fix every half-edge
individually and hence every incident vertex. A connected graph with
an edge has no isolated vertex, so the automorphism would be the
identity. If $E(\Gamma)=0$, the graph has a single vertex and trivial
automorphism group; no nontrivial-stabilizer pattern occurs at all.

\subsection{Finite residual quotient families}
\label{sec:finite-residual-families}

Fix a support pattern $p$ with $H_p\ne\{1\}$ and omit the chains on
$L_p$. The remaining component positions are the vertices of $\Gamma$
and the edges with $a_p(e)=1$. For a finite label set $B$, assign a
subset $B_v$ to each vertex and a subset $B_e$ to each such edge,
subject to
\[
B=\bigsqcup_{v\in\Vertices(\Gamma)}B_v
\sqcup\bigsqcup_{e:a_p(e)=1}B_e,\qquad
B_v\ne\varnothing\ \Longleftrightarrow\ v\in S_p,\qquad
B_e\ne\varnothing.
\]
For each allocation, take the coarse decoration parameter space of
\eqref{eq:decoration-parameter-space}: its vertex factors are
$\Vspace_{h_v,S_v,B_v}$, and each edge with $a_p(e)=1$ contributes
$\mathcal C_{B_e}$. Thus
\[
Z_{p,B}:=
\coprod_{\substack{
B=\bigsqcup_v B_v\sqcup\bigsqcup_{e:a_p(e)=1}B_e\\
B_v\ne\varnothing\ \Longleftrightarrow\ v\in S_p,\quad B_e\ne\varnothing}}
\left(\prod_v\Vspace_{h_v,S_v,B_v}
\times\prod_{e:a_p(e)=1}\mathcal C_{B_e}\right).
\]
The incidence action of $H_p$ on these parameter spaces is retained.
Define
\[
Z_p^\mu(z)=\sum_{n\ge0}\mu(Z_{p,[n]})\frac{z^n}{n!},\qquad
Z_p^{\mu,\mathrm{quot}}(z)=
\sum_{n\ge0}\mu(Z_{p,[n]}/H_p)\frac{z^n}{n!}.
\]

Write $T_\lambda=P_\lambda-1-R_\lambda$ for the series of chains
with at least two components. Every element of $H_p$ fixes the two
ends of every edge of $L_p$ individually. Its lift therefore acts
trivially on each long-chain parameter factor, including its light
coordinates. For each label allocation the relevant quotient is
\[
\left(Z\times\prod_{e\in L_p}W_e\right)/H_p
\cong (Z/H_p)\times\prod_{e\in L_p}W_e,
\]
where $W_e$ is the corresponding chain factor. The label allocations
are themselves preserved by $H_p$.
Grouping \eqref{eq:all-decorations-correction} by $p$ consequently gives
\begin{equation}
\mathcal R_\Gamma^\mu
=\sum_{p:H_p\ne\{1\}}
\frac{|H_p|}{|\Aut\Gamma|}
\left(Z_p^{\mu,\mathrm{quot}}-\frac1{|H_p|}Z_p^\mu\right)
T_\lambda^{m_p}.
\label{eq:pattern-correction}
\end{equation}

To see how the extraction works, consider one long edge in
\eqref{eq:pattern-correction}. In its chain factor $T_\lambda$,
taking $P_\lambda$ extracts an unrestricted propagator. Taking the
subtracted term $-(1+R_\lambda)$ leaves an empty edge or a single inserted component
in the residual spaces and contributes a minus sign. The inherited
group acts trivially on these optional chain factors. Making this
choice independently on every long edge gives
\begin{equation}
T_\lambda^{m_p}
=\sum_{J\subseteq L_p}(-1)^{m_p-|J|}
P_\lambda^{|J|}(1+R_\lambda)^{m_p-|J|}.
\label{eq:signed-chain-extraction}
\end{equation}
A \emph{residual type} is a pair $\tau=(p,J)$ with $H_p\ne\{1\}$
and $J\subseteq L_p$. It records the pattern and the edges on which
we chose the unrestricted propagator. Let $\mathscr T_\Gamma$ be this
finite set and put
\[
G_{\Gamma,\tau}=H_p,\quad k_{\Gamma,\tau}=|J|,\quad
c_{\Gamma,\tau}=\frac{|H_p|}{|\Aut\Gamma|}(-1)^{m_p-|J|}.
\]

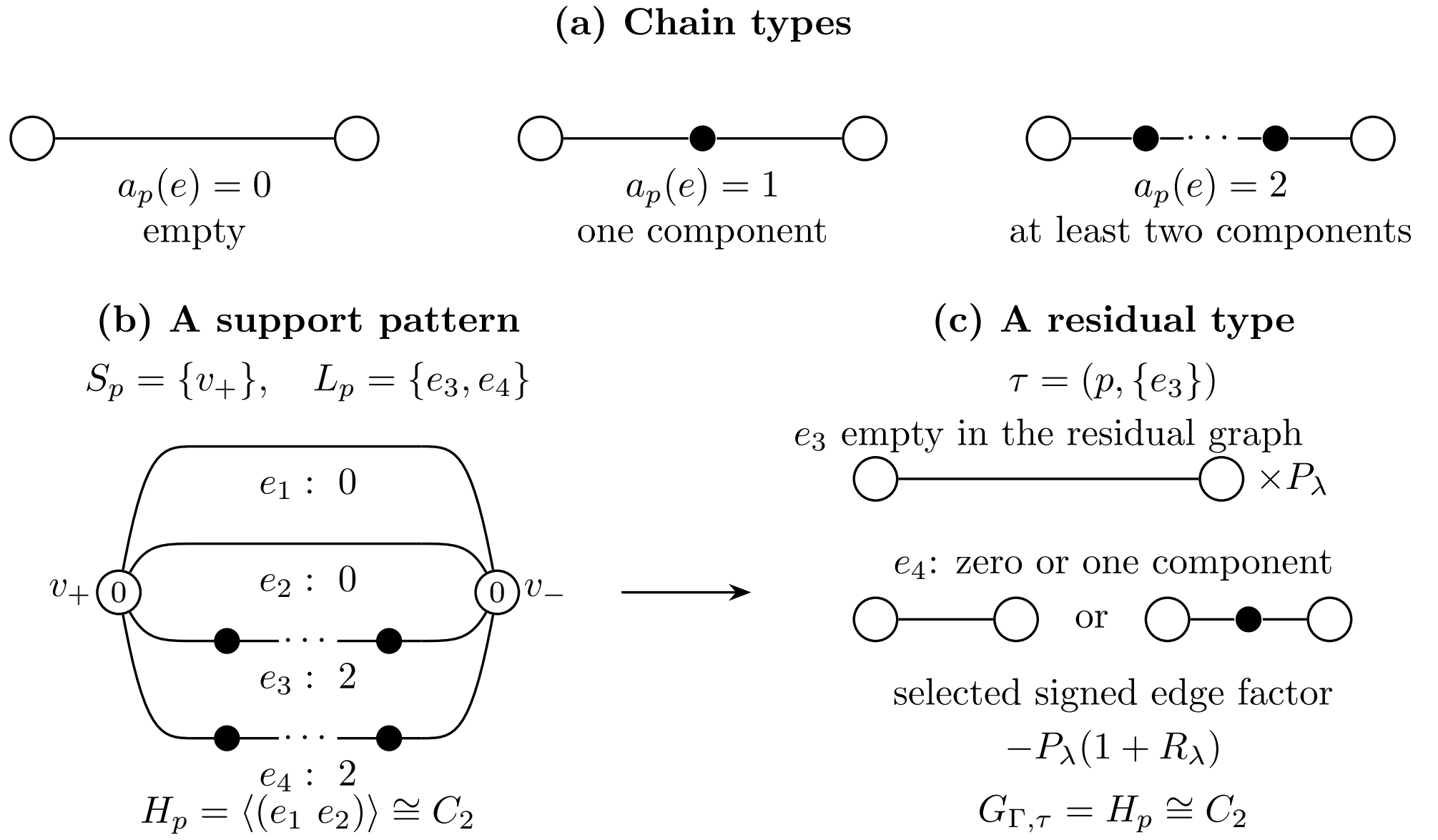
\begin{figure}[htbp]
\centering
\begin{tikzpicture}[line width=.65pt,>=Stealth,
 every node/.style={font=\small},
 graph/.style={circle,draw,line width=.65pt,minimum size=4mm,
 inner sep=0pt,fill=white,font=\scriptsize},
 inserted/.style={circle,fill,minimum size=2.4mm,inner sep=0pt}]
\node at (6.2,1.05) {\textbf{(a) Chain types}};
\foreach \x/\k in {0/0,4.7/1,9.4/2}{
 \node[graph] (a\k) at (\x,0) {};
 \node[graph] (b\k) at (\x+3,0) {};
}
\draw (a0)--(b0);
\node[inserted] (u1) at (6.2,0) {};
\draw (a1)--(u1)--(b1);
\node[inserted] (u2) at (10.3,0) {};
\node[inserted] (w2) at (11.5,0) {};
\draw (a2)--(u2)--(10.65,0);
\node at (10.9,0) {$\cdots$};
\draw (11.15,0)--(w2)--(b2);
\node[align=center] at (1.5,-.65) {$a_p(e)=0$\\empty};
\node[align=center] at (6.2,-.65) {$a_p(e)=1$\\one component};
\node[align=center] at (10.9,-.65) {$a_p(e)=2$\\at least two components};

\node at (2.55,-1.7) {\textbf{(b) A support pattern}};
\node at (2.55,-2.25) {$S_p=\{v_+\},\quad L_p=\{e_3,e_4\}$};
\node[graph] (v) at (.8,-4.2) {$0$};
\node[graph] (w) at (4.3,-4.2) {$0$};
\node[left=3pt] at (v) {$v_+$};
\node[right=3pt] at (w) {$v_-$};
\foreach \y in {-2.85,-3.75,-4.65,-5.55}{
 \draw (v) .. controls (1.05,\y) and (1.2,\y) .. (1.5,\y);
 \draw (3.6,\y) .. controls (3.9,\y) and (4.05,\y) .. (w);
}
\draw (1.5,-2.85)--(3.6,-2.85);
\draw (1.5,-3.75)--(3.6,-3.75);
\foreach \y in {-4.65,-5.55}{
 \node[inserted] (u\y) at (1.8,\y) {};
 \node[inserted] (w\y) at (3.3,\y) {};
 \draw (1.5,\y)--(2.28,\y);
 \draw (2.82,\y)--(3.6,\y);
 \node at (2.55,\y) {$\cdots$};
}
\node[below=2pt] at (2.55,-2.85) {$e_1:\ 0$};
\node[below=2pt] at (2.55,-3.75) {$e_2:\ 0$};
\node[below=2pt] at (2.55,-4.65) {$e_3:\ 2$};
\node[below=2pt] at (2.55,-5.55) {$e_4:\ 2$};
\node at (2.55,-6.25) {$H_p=\langle(e_1\ e_2)\rangle\cong C_2$};

\draw[->] (5.45,-4.2)--(6.65,-4.2);
\node at (10,-1.7) {\textbf{(c) A residual type}};
\node at (10,-2.25) {$\tau=(p,\{e_3\})$};
\node[graph] (r3a) at (7.8,-3.15) {};
\node[graph] (r3b) at (11,-3.15) {};
\draw (r3a)--(r3b);
\node[above=3pt] at (9.4,-3.15) {$e_3$ empty in the residual graph};
\node[right=5pt] at (r3b) {$\times P_\lambda$};
\node[graph] (r4a) at (7.8,-4.45) {};
\node[graph] (r4b) at (9.1,-4.45) {};
\draw (r4a)--(r4b);
\node at (9.8,-4.45) {or};
\node[graph] (r4c) at (10.5,-4.45) {};
\node[graph] (r4d) at (12,-4.45) {};
\node[inserted] (r4u) at (11.25,-4.45) {};
\draw (r4c)--(r4u)--(r4d);
\node at (10,-3.95) {$e_4$: zero or one component};
\node at (10,-5.15) {selected signed edge factor};
\node at (10,-5.65) {$-P_\lambda(1+R_\lambda)$};
\node at (10,-6.25) {$G_{\Gamma,\tau}=H_p\cong C_2$};
\end{tikzpicture}
\caption{From chain support to residual geometry. (a) The three values
of $a_p(e)$. (b) A genus-three graph with four parallel edges: only
$v_+$ carries light labels, $e_1,e_2$ are empty, and $e_3,e_4$ are
long. The displayed group fixes every label in any decoration with
this pattern. (c) Choosing $J=\{e_3\}$ extracts $P_\lambda$ from $e_3$
and leaves an optional empty or one-component chain on $e_4$.
The selected term in \eqref{eq:signed-chain-extraction} has sign minus.
The graph vertex factors and their prescribed $H_p$-action are retained;
the diagrams of the residual edges show only their allowed component
positions. This is an algebraic expansion, not a decomposition into
strata.}
\label{fig:support-pattern-residual-type}
\end{figure}

For a finite label set $B$ define $\mathcal F_\varnothing=\operatorname{pt}$
and $\mathcal F_B=\mathcal C_B$ if $B\ne\varnothing$. The \emph{residual family} for $\tau$ collects the spaces left after
these choices, as the number $n$ of available labels varies. Define
its members by
\begin{equation}
Y_{\Gamma,\tau,n}:=
\coprod_{[n]=B_0\sqcup\bigsqcup_{e\in L_p\setminus J}B_e}
Z_{p,B_0}\times\prod_{e\in L_p\setminus J}\mathcal F_{B_e}.
\label{eq:residual-varieties}
\end{equation}
Thus edges with $a_p(e)=1$ still carry exactly one component; edges
in $L_p\setminus J$ carry zero or one; all remaining edges, including
$J$, are empty. Graph vertex support is unchanged. The action of $H_p$
on $Z_{p,B_0}$ is extended trivially to the optional factors.
Let
\[
Y_{\Gamma,\tau}^\mu(z)=\sum_{n\ge0}\mu(Y_{\Gamma,\tau,n})\frac{z^n}{n!},
\qquad
Q_{\Gamma,\tau}^\mu(z)=\sum_{n\ge0}
\mu(Y_{\Gamma,\tau,n}/G_{\Gamma,\tau})\frac{z^n}{n!}.
\]
Residual geometry has a bounded number of component positions, but
the number of residual labels is unrestricted.
Figure~\ref{fig:support-pattern-residual-type} illustrates the passage
from a support pattern to one term of the residual expansion.

\begin{lemma}[Finite residual decomposition]
\label{lem:residual-decomposition}
For every stabilization graph $\Gamma$,
\begin{equation}
\mathcal R_\Gamma^\mu(z)
=\sum_{\tau\in\mathscr T_\Gamma}c_{\Gamma,\tau}
\left(Q_{\Gamma,\tau}^\mu(z)
-\frac{Y_{\Gamma,\tau}^\mu(z)}{|G_{\Gamma,\tau}|}\right)
P_\lambda(z)^{k_{\Gamma,\tau}}.
\label{eq:residual-decomposition}
\end{equation}
Every term has $k_{\Gamma,\tau}\le E(\Gamma)-1$. When $E(\Gamma)=0$
the correction is zero and the sum is empty.
\end{lemma}
\begin{proof}
The exponential product rule and the trivial action on the optional
factors give
\[
Y_{\Gamma,\tau}^\mu=Z_p^\mu(1+R_\lambda)^{m_p-|J|},\qquad
Q_{\Gamma,\tau}^\mu=Z_p^{\mu,\mathrm{quot}}
(1+R_\lambda)^{m_p-|J|}.
\]
Substitute these identities and \eqref{eq:signed-chain-extraction}
into \eqref{eq:pattern-correction}. The coefficient is precisely
$c_{\Gamma,\tau}$; the bound follows from
\eqref{eq:residual-propagator-bound}.
\end{proof}

\begin{remark}[The specified residual action]
This is a signed algebraic expansion, not a stratification of the
original graph locus by residual types. After extracting or shortening
a long chain, the residual decoration may have a larger stabilizer.
The quotient in \eqref{eq:residual-decomposition} is nevertheless by
the prescribed subgroup $H_p$, acting as constructed above, not by
that larger stabilizer. The signs belong to the coefficients
$c_{\Gamma,\tau}$; the residual families themselves are varieties.
\qrem
\end{remark}

\subsection{Example: the theta graph}
\label{ex:theta-graph}

\begin{figure}[htbp]
\centering
\begin{tikzpicture}[line width=.65pt,>=Stealth,
 every node/.style={font=\small},
 graph/.style={circle,draw,line width=.65pt,minimum size=4mm,
 inner sep=0pt,fill=white,font=\scriptsize},
 inserted/.style={circle,fill,minimum size=2.4mm,inner sep=0pt}]
\node at (1.65,1.65) {\textbf{(a) Original decoration}};
\node at (9.05,1.65) {\textbf{(b) Residual decoration}};
\node[graph] (a) at (0,0) {$0$};
\node[graph] (b) at (3.3,0) {$0$};
\node[inserted] (u) at (1.1,0) {};
\node[inserted] (w) at (2.2,0) {};
\draw (a)--(u)--(w)--(b);
\draw (a) to[bend left=60] node[above] {$e_2$} (b);
\draw (a) to[bend right=60] node[below] {$e_3$} (b);
\node[above=4pt] at (u) {$\{2\}$};
\node[above=4pt] at (w) {$\{3\}$};
\node[left=3pt,align=center] at (a) {$v_+$\\$\{1\}$};
\node[right=3pt] at (b) {$v_-$};
\node[below=4pt] at (1.65,0) {$e_1$};
\node at (1.65,-1.5) {$G_\delta=H_p=\langle(e_2\ e_3)\rangle$};
\node at (1.65,-2.05) {$H_p\cong C_2$};
\draw[->] (4.35,0)--(6.35,0);
\node[above=4pt] at (5.35,0) {$J=\{e_1\}$};
\node[below=4pt] at (5.35,0) {factor $P_\lambda$};
\node[graph] (c) at (7.4,0) {$0$};
\node[graph] (d) at (10.7,0) {$0$};
\draw (c)--(d);
\draw (c) to[bend left=60] node[above] {$e_2$} (d);
\draw (c) to[bend right=60] node[below] {$e_3$} (d);
\node[below=4pt] at (9.05,0) {$e_1$};
\node[left=3pt,align=center] at (c) {$v_+$\\$\{1\}$};
\node[right=3pt] at (d) {$v_-$};
\node at (9.05,-1.5) {full stabilizer: $\mathfrak S_3$};
\node at (9.05,-2.05) {quotient group: $G_{\Gamma,\tau_1}=H_p\cong C_2$};
\end{tikzpicture}
\caption{The prescribed action survives extraction. (a) A labeled
decoration of the theta graph with $L_p=\{e_1\}$ and
$H_p=\langle(e_2\ e_3)\rangle$. (b) A residual configuration for
$\tau_1=(p,\{e_1\})$, with $e_1$ empty and only label $1$ retained.
All three empty edges can now be permuted, so the full discrete
stabilizer is $\mathfrak S_3$. The residual quotient is nevertheless
by the original subgroup $H_p\cong C_2$. The arrow records the term
with factor $P_\lambda$ in the signed expansion; it is not a morphism
between the two strata.}
\label{fig:theta-residual-example}
\end{figure}
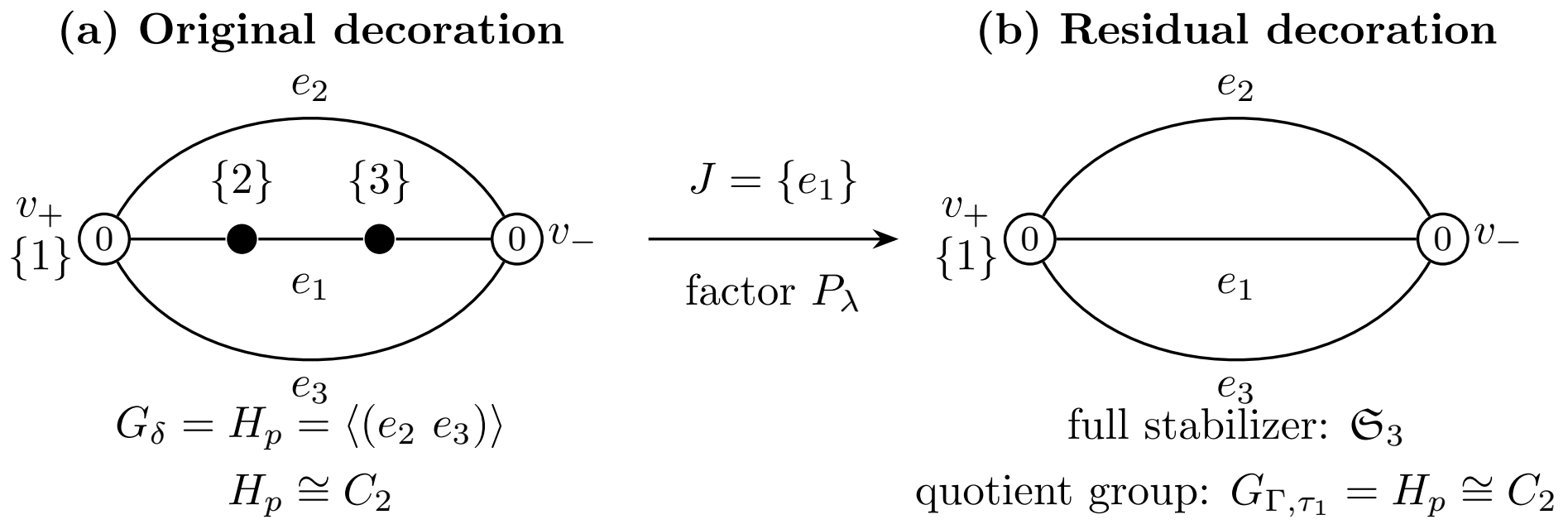

Let $\Gamma$ consist of two genus-zero vertices $v_+,v_-$ joined by
three edges $e_1,e_2,e_3$. Its automorphism group is
$\mathfrak S_3\times\mathbb Z/2\mathbb Z$, of order twelve. Choose
\[
S_p=\{v_+\},\qquad a_p(e_1)=2,\qquad a_p(e_2)=a_p(e_3)=0.
\]
The occupied vertex prevents exchange of $v_+$ and $v_-$; the long
chain fixes the two half-edges of $e_1$. The exact decoration
stabilizer is $H_p=\langle(e_2\ e_3)\rangle$, of order two.
Here $L_p=\{e_1\}$, and the coefficient $|H_p|/|\Aut\Gamma|$ in
\eqref{eq:pattern-correction} is $1/6$ for this pattern on the fixed graph.
Figure~\ref{fig:theta-residual-example} compares a decoration with
this pattern and a residual configuration for the type
$\tau_1=(p,\{e_1\})$ considered below.

Set $U=\mathbb P^1\setminus\{0,1,\infty\}$. Normalize the heavy
branches at $v_+$ indexed by $e_1,e_2,e_3$ to $\infty,0,1$.
The unoccupied trivalent vertex contributes a point, and
\[
Z_{p,B}=\begin{cases}\varnothing,&B=\varnothing,\\U^B,&B\ne\varnothing.
\end{cases}
\]
The generator of $H_p$ acts diagonally by $x_b\mapsto1-x_b$.
Thus the stabilizer is combinatorial, but its action on the parameter
space is nontrivial. In particular,
\[
Z_p^\mu(z)=e^{(\lambda-2)z}-1,
\]
whereas $Z_p^{\mu,\mathrm{quot}}$ retains the simultaneous quotient.

There are two residual types, $\tau_1=(p,\{e_1\})$ and
$\tau_0=(p,\varnothing)$, with coefficients
$c_{\Gamma,\tau_1}=1/6$ and $c_{\Gamma,\tau_0}=-1/6$.
Their series satisfy
\begin{align*}
Y_{\Gamma,\tau_1}^\mu&=Z_p^\mu,&
Q_{\Gamma,\tau_1}^\mu&=Z_p^{\mu,\mathrm{quot}},\\
Y_{\Gamma,\tau_0}^\mu&=Z_p^\mu(1+R_\lambda),&
Q_{\Gamma,\tau_0}^\mu&=Z_p^{\mu,\mathrm{quot}}(1+R_\lambda).
\end{align*}
The two terms combine to the required contribution
\[
\frac16\left(Z_p^{\mu,\mathrm{quot}}-\frac12Z_p^\mu\right)
\bigl(P_\lambda-(1+R_\lambda)\bigr).
\]
This is the contribution of the chosen pattern on the fixed graph, not the
entire theta graph. When the residual edge $e_1$ is empty, all three
empty edges can be permuted, so the residual decoration has stabilizer
$\mathfrak S_3$, larger than $H_p$. The residual quotient still uses
the original order-two action. This illustrates why that action must
be kept independently of the new decoration stabilizer.

\subsection{Fixed-graph formula and the trivalent term}
\label{sec:leadingmotivicpropagator}

\begin{proposition}[Motivic contribution of a stabilization graph]
\label{prop:motivicgraph}
For a stabilization graph $\Gamma$ of genus $g\ge2$,
\begin{equation}
A_{g,\Gamma}^\mu(z)=
\frac1{|\Aut\Gamma|}
\left(\prod_{v\in\Vertices(\Gamma)}V_{h_v,d_v}^\mu(z)\right)
P_\lambda(z)^{E(\Gamma)}+\mathcal R_\Gamma^\mu(z),
\label{eq:motivicgraph}
\end{equation}
where the correction has the finite residual representation
\eqref{eq:residual-decomposition}. If $E(\Gamma)\ge1$, every
correction term contains at most $E(\Gamma)-1$ unrestricted
propagators. If $E(\Gamma)=0$, the correction vanishes.
\end{proposition}
\begin{proof}
Combine \eqref{eq:section4-graph-correction} with the finite residual
decomposition.
\end{proof}

Let $\mathfrak C_g$ and $\mathfrak T_g$ be the sets of graphs and
rational trivalent graphs defined in Section~\ref{sec:all-light-decomposition}, and recall
$\gamma_g=\sum_{\Gamma\in\mathfrak T_g}|\Aut\Gamma|^{-1}$.

\begin{theorem}[Motivic stabilization-graph expansion]
\label{thm:motivicgraphexpansion}
For fixed $g\ge2$,
\begin{equation}
A_g^\mu(z)=\gamma_g e^{2(g-1)(\lambda-2)z}
P_\lambda(z)^{3g-3}+\mathcal R_g^\mu(z).
\label{eq:motivicgraphexpansion}
\end{equation}
More explicitly, writing
$D_{\Gamma,\tau}^\mu=Q_{\Gamma,\tau}^\mu-
Y_{\Gamma,\tau}^\mu/|G_{\Gamma,\tau}|$, the remainder is
\begin{equation}
\mathcal R_g^\mu=
\sum_{\Gamma\in\mathfrak C_g\setminus\mathfrak T_g}
\frac{\prod_v V_{h_v,d_v}^\mu}{|\Aut\Gamma|}
P_\lambda^{E(\Gamma)}
+\sum_{\Gamma\in\mathfrak C_g}\ \sum_{\tau\in\mathscr T_\Gamma}
c_{\Gamma,\tau}D_{\Gamma,\tau}^\mu
P_\lambda^{k_{\Gamma,\tau}}.
\label{eq:explicit-global-remainder}
\end{equation}
Both sums are finite and every displayed propagator exponent is at
most $3g-4$.
\end{theorem}
\begin{proof}
By the edge bound from Section~\ref{sec:all-light-decomposition}, a nontrivalent graph has at most
$3g-4$ edges. A rational trivalent graph has $2g-2$ vertices and
$3g-3$ edges. Each of its vertices contributes
$V_{0,3}^\mu=e^{(\lambda-2)z}$, so the sum of their normalized
parameter contributions is the first term of
\eqref{eq:motivicgraphexpansion}. Their corrections, like all other
graph corrections, have exponents at most $E(\Gamma)-1\le3g-4$.
Summing \eqref{eq:motivicgraph} over the finite set $\mathfrak C_g$
gives the explicit remainder.
\end{proof}

The bound on the powers of $P_\lambda$ refers to this constructed
finite representation, with its specified vertex and residual series.
It is not an intrinsic degree on arbitrary formal series. In particular,
it does not by itself bound the analytic singularities of the
coefficients $V_{h,d}^\mu$, $Y_{\Gamma,\tau}^\mu$ and
$Q_{\Gamma,\tau}^\mu$. Section~\ref{sec:motivic-asymptotics} supplies sufficient bounds through
tameness before drawing asymptotic conclusions.

For the Euler specialization, Burnside--Lefschetz identifies each
quotient defect with the nonidentity fixed sectors:
\[
\chi(X/H)-\frac{\chi(X)}{|H|}
=\frac1{|H|}\sum_{\substack{\sigma\in H\\\sigma\ne1}}
\chi(X^\sigma).
\]
This recovers the correction studied in \cite[Section~4.4]{CottiPaperI}.
Here the residual coarse quotients have been retained throughout,
so the construction applies before selecting a realization. The
universal trivalent term specializes, at $\lambda=1$, to
$\gamma_g e^{-(2g-2)z}(1-z)^{-(3g-3)}$, the leading Euler term of
\cite[Section~6.1]{CottiPaperI}.

\section{The universal motivic heavy/light transformation}
\label{sec:motivic-heavy-light}

Having described the all-light series through stabilization graphs,
we now compare it with the ordinary series by varying the marking
weights. The basic geometric step is the crossing of a single collision
wall, for which the contracted rational tail gives an identity of coarse
classes in the Grothendieck ring. Applying this identity while lowering
one heavy marking at a time yields a recursion, and its generating
series is governed by a change of coordinates independent of the genus.
We then identify the inverse coordinate through the ordinary genus-zero
spaces and relate its derivative to the edge propagator. This last
identity will connect the wall-crossing geometry with the singularities
used in the asymptotic analysis.

\subsection{Reduction morphisms and collision walls}
\label{sec:reduction-morphisms}

We retain the Hassett notation of Section~\ref{sec:all-light-decomposition}. Let $\bm a=(a_1,\ldots,a_n)$
and $\bm b=(b_1,\ldots,b_n)$ be admissible weight vectors for genus $g$,
with $\bm a\ge\bm b$ coordinatewise. Lowering the weights permits
additional collisions and may destabilize rational tails. Hassett's
reduction morphism
\[
\Mbar_{g,\bm a}\longrightarrow\Mbar_{g,\bm b}
\]
stabilizes the resulting weighted curve by contracting the components
that become unstable; it induces a morphism of projective coarse spaces
\cite[Theorem~4.1]{Hassett}. We will use the coarse morphism throughout
the wall calculation.

\begin{definition}
Admissible weight vectors $\bm a\ge\bm b$ are \emph{separated by a
single wall indexed by $I$} if $I\subseteq[n]$ is the unique subset
satisfying
\[
\sum_{i\in I}b_i\le1<\sum_{i\in I}a_i.
\]
\end{definition}
The markings are labeled, so the subset $I$ is fixed and no permutation
of its labels is divided out.

\subsection{The motivic single-wall identity}
\label{sec:single-wall}

\begin{lemma}
\label{lem:motivic-wall}
Let $\bm a\ge\bm b$ be admissible weight vectors separated by a single
collision wall indexed by $I$, with $k=|I|\ge2$. Let
$Z_I\subseteq\overline M_{g,\bm b}$ be the locus where the markings
indexed by $I$ coincide. Then
\begin{equation}
[\overline M_{g,\bm a}]-[\overline M_{g,\bm b}]
=\bigl([\mathbb P^{k-2}]-1\bigr)[Z_I]
=\left(\sum_{j=1}^{k-2}\LL^j\right)[Z_I]
\label{eq:motivic-single-wall}
\end{equation}
in $K_0(\Var_{\C})$. The sum is empty for $k=2$.
\end{lemma}

\begin{figure}[htbp]
\centering
\begin{tikzpicture}[>=Stealth,line width=.65pt,
 every node/.style={font=\small}]
\node at (1,2.15) {\textbf{(a) Before reduction}};
\node at (8.1,2.15) {\textbf{(b) After reduction}};
\node at (1,1.55) {$\sum_{i\in I}a_i>1$};
\node at (8.1,1.55) {$\sum_{i\in I}b_i\le1$};
\draw (0,0) ellipse [x radius=1.6,y radius=.9];
\node at (0,.15) {$C$};
\node[font=\footnotesize] at (0,-.3) {genus $g$};
\draw (2.4,0) circle (.8);
\node at (2.4,0) {$\mathbb P^1$};
\fill (1.6,0) circle (1.6pt);
\node[below left=2pt] at (1.6,0) {$q$};
\fill (2.674,.752) circle (1.3pt);
\node[above=3pt] at (2.674,.752) {$p_{i_1}$};
\fill (3.2,0) circle (1.3pt);
\node[right=3pt] at (3.2,0) {$p_{i_2}$};
\fill (2.674,-.752) circle (1.3pt);
\node[below=3pt] at (2.674,-.752) {$p_{i_k}$};
\node[font=\scriptsize] at (3.17,-.48) {$\vdots$};
\draw[->,thick] (4.3,0)--(6.2,0);
\node[above=4pt] at (5.25,0) {reduction};
\draw (8.1,0) ellipse [x radius=1.6,y radius=.9];
\node at (8.1,.15) {$C$};
\node[font=\footnotesize] at (8.1,-.3) {genus $g$};
\fill (9.7,0) circle (1.6pt);
\node[right=4pt,align=left] at (9.7,0) {labels $I$\\coincide};
\node at (1,-2) {$D_I\simeq\mathbb P^{k-2}\times Z_I$};
\node at (8.1,-2) {$Z_I$};
\draw[->] (3.6,-2)--(6.8,-2);
\node[above=3pt] at (5.2,-2) {projection};
\node at (5.2,-2.95)
 {$[\overline M_{g,\bm a}]-[\overline M_{g,\bm b}]
   =\bigl([\mathbb P^{k-2}]-1\bigr)[Z_I]$};
\end{tikzpicture}
\caption{A single collision wall with $k=|I|\ge3$. The reduction
contracts the rational tail carrying the labels of $I$ to their
common smooth marked point on $C$. The main curve is schematic;
its other markings and possible nodes are omitted. The tail moduli
form $\mathbb P^{k-2}$, and the restriction to $D_I$ is the displayed
projection. The complements of $D_I$ and $Z_I$ are isomorphic,
which gives the motivic difference below the diagrams.
For $k=2$ the reduction is an isomorphism and the difference vanishes.}
\label{fig:single-wall-contraction}
\end{figure}
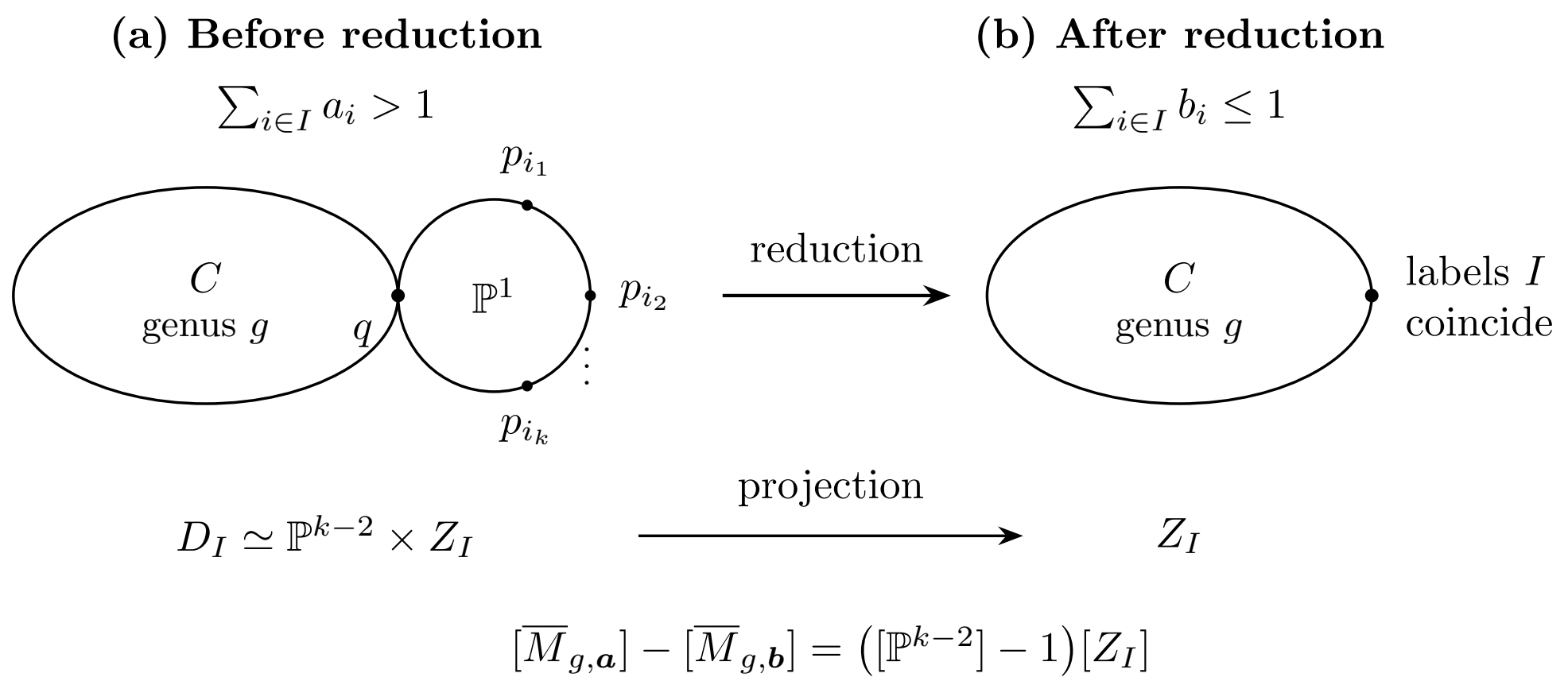

\begin{proof}
We use Hassett's description of reduction morphisms
\cite[Theorem~4.1, Proposition~4.5, Remark~4.6 and Corollary~4.7]{Hassett},
applied to coarse spaces. Put $J=[n]\setminus I$ and
$b_I=\sum_{i\in I}b_i$. The single-wall assumption gives
\[
\sum_{i\in S}a_i\le1\quad(S\subsetneq I),\qquad
b_I+b_j>1\quad(j\in J).
\]
Otherwise the wall indexed by $S$ or $I\cup\{j\}$ would also be
crossed. Merging the markings in $I$ identifies
\begin{equation}
\begin{aligned}
Z_I&\simeq\overline M_{g,((b_j)_{j\in J},b_I)}\\
&\simeq\overline M_{g,((b_j)_{j\in J},1)}
\simeq\overline M_{g,((a_j)_{j\in J},1)}.
\end{aligned}
\label{eq:wall-center-identifications}
\end{equation}
Indeed, the merged marking cannot collide with any other marking.
A rational tail carrying it alone is unstable for either weight, and
a tail carrying it together with any other marking has total weight
greater than one. No wall involving only labels in $J$ is crossed.

For $k\ge3$, the reduction is an isomorphism away from $Z_I$.
Its inverse image $D_I$ is the divisor with a rational tail carrying
precisely the labels in $I$. Hassett's gluing description gives the
tail factor and the main component with its attaching marking:
\[
D_I\simeq\overline M_{0,((a_i)_{i\in I},1)}\times Z_I
\simeq\mathbb P^{k-2}\times Z_I.
\]
For completeness, the tail is irreducible: any terminal component
not containing its attaching point would carry a proper subset of
$I$ of weight greater than one. Put the attaching point at infinity.
The remaining moduli are
\[
(\mathbb A^k\setminus\Delta)/\Aff_1\simeq\mathbb P^{k-2},
\qquad
\Delta=\{x\in\mathbb A^k:x_1=\cdots=x_k\}.
\]
Every proper subset of the tail markings may coincide, but all $k$
markings cannot coincide. Translation sets the last coordinate to
zero, and dilation projectivizes the nonzero vector of differences.
The attaching point and labeled markings leave no nontrivial tail
automorphisms; automorphisms of the main component act trivially on
this factor. Thus the product description also holds on coarse spaces,
and the restriction of the reduction to $D_I$ is the projection to
$Z_I$.
Figure~\ref{fig:single-wall-contraction} illustrates this contraction
and its restriction to the exceptional locus for $k\ge3$.

The complements are isomorphic:
\[
\overline M_{g,\bm a}\setminus D_I
\simeq\overline M_{g,\bm b}\setminus Z_I.
\]
The scissors relation and multiplicativity now give
\[
[\overline M_{g,\bm a}]-[\overline M_{g,\bm b}]
=[\mathbb P^{k-2}\times Z_I]-[Z_I]
=\bigl([\mathbb P^{k-2}]-1\bigr)[Z_I].
\]
Since $[\mathbb P^{k-2}]=1+\LL+\cdots+\LL^{k-2}$, this is
\eqref{eq:motivic-single-wall}. For $k=2$ the reduction is an
isomorphism by \cite[Corollary~4.7]{Hassett}, and both sides vanish.
\end{proof}

This proof follows the geometric argument for the Euler identity in
\cite[Section~2.3]{CottiPaperI}. Its product description on coarse spaces
allows the final step to be made in the Grothendieck ring. For any
motivic measure with Lefschetz value $\lambda$, it gives
\[
\mu(\overline M_{g,\bm a})-\mu(\overline M_{g,\bm b})
=\left(\sum_{j=1}^{k-2}\lambda^j\right)\mu(Z_I).
\]
For the Euler specialization, $\lambda=1$ and the wall coefficient
is $k-2$.

\subsection{Recursion and the universal change of coordinates}
\label{sec:wall}

Fix $g\ge1$. For $r,s\ge0$ write
$\mathfrak h_{r,s}=[\overline M_{g,r|s}]$ in $K_0(\Var_{\C})$.
Set $\mathfrak h_{0,0}=0$ when $g=1$; for $g\ge2$ it is
$[\overline M_g]$. For $m\ge2$, define the wall polynomial
\[
p_m(\lambda):=\lambda+\cdots+\lambda^{m-1}.
\]

\begin{proposition}[Motivic one-step recursion]
\label{prop:motivic-onestep}
For $r\ge1$ and $s\ge0$,
\begin{equation}
\mathfrak h_{r,s}-\mathfrak h_{r-1,s+1}
=\sum_{m=2}^s\binom{s}{m}p_m(\LL)\mathfrak h_{r,s-m}
\label{eq:motivic-onestep}
\end{equation}
in $K_0(\Var_{\C})$. The sum is empty for $s=0,1$.
\end{proposition}

\begin{proof}
We follow the weight-lowering argument for the Euler recursion in
\cite[Section~2.4]{CottiPaperI}, using the motivic
single-wall identity \eqref{eq:motivic-single-wall} at each crossing.
Choose positive rational weights $\epsilon_1,\ldots,\epsilon_s$
with total less than one and distinct subset sums, and choose
$\eta>0$ with $\eta+\sum_i\epsilon_i<1$. Lower the distinguished
heavy weight in
\[
\bm w(t)=(1^{r-1},t,\epsilon_1,\ldots,\epsilon_s),
\qquad \eta\le t\le1.
\]
The endpoints represent $\overline M_{g,r|s}$ and
$\overline M_{g,r-1|s+1}$. Although the light weights used along
the path are unequal, both endpoints have the collision and
stability conditions of the indicated minimal heavy/light chambers.
All weights along the path are admissible for $g\ge1$.

The walls crossed are precisely
$t+\sum_{i\in S}\epsilon_i=1$ for nonempty subsets $S\subseteq[s]$.
Their crossing values are distinct, so each crossing is a single
wall. Fix $|S|=m$. Immediately below this wall, the merged marking
has weight $1-\delta$, where $\delta>0$ can be chosen smaller than
every remaining light weight. It cannot collide with another marking,
and increasing its weight to one does not change stability. The
collision center is therefore canonically
$Z_S\simeq\overline M_{g,r|s-m}$.

For $m\ge2$, the wall involves $m+1$ markings, so the single-wall
identity gives the contribution
\[
\bigl([\mathbb P^{m-1}]-1\bigr)[Z_S]
=p_m(\LL)\mathfrak h_{r,s-m}.
\]
For $m=1$ the contribution is zero. Let $X_0,\ldots,X_\ell$ be the
successive coarse spaces along the path, with
$X_0=\overline M_{g,r|s}$ and
$X_\ell=\overline M_{g,r-1|s+1}$. Then
\begin{align*}
\mathfrak h_{r,s}-\mathfrak h_{r-1,s+1}
&=[X_0]-[X_\ell]
=\sum_{i=1}^{\ell}\bigl([X_{i-1}]-[X_i]\bigr)\\
&=\sum_{m=2}^s\ \sum_{\substack{S\subseteq[s]\\|S|=m}}
p_m(\LL)\mathfrak h_{r,s-m}
=\sum_{m=2}^s\binom{s}{m}p_m(\LL)\mathfrak h_{r,s-m}.
\end{align*}
For $s=0$ there is no crossing and both sums are empty.
\end{proof}

\begin{corollary}[Motivic form of Y.-H.\,Kiem's recursion]
\label{cor:motivic-kiem}
For $r,s\ge0$ and $N=r+s>0$,
\begin{equation}
\mathfrak h_{r,s}=\mathfrak h_{0,N}
+\sum_{j=1}^{\min(r,N-2)}\ \sum_{k=0}^{N-j-2}
\binom{N-j}{k}\mathfrak h_{j,k}p_{N-j-k}(\LL).
\label{eq:motivic-kiem}
\end{equation}
Empty sums vanish.
\end{corollary}
\begin{proof}
At fixed $N=r+s$, summing \eqref{eq:motivic-onestep} over
$j=1,\ldots,r$ gives
\begin{multline*}
\mathfrak h_{r,s}-\mathfrak h_{0,N}
=\sum_{j=1}^r
\bigl(\mathfrak h_{j,N-j}-\mathfrak h_{j-1,N-j+1}\bigr)\\
=\sum_{j=1}^r\sum_{m=2}^{N-j}
\binom{N-j}{m}p_m(\LL)\mathfrak h_{j,N-j-m}
=\sum_{j=1}^{\min(r,N-2)}\sum_{k=0}^{N-j-2}
\binom{N-j}{k}\mathfrak h_{j,k}p_{N-j-k}(\LL),
\end{multline*}
where the last equality uses $k=N-j-m$.
\end{proof}

Hodge--Deligne realization of \eqref{eq:motivic-kiem} gives the
recursion of \cite[Theorem~3.1]{Kiem}.

Define
\begin{align*}
\mathcal H_{g,r}(z)&=\sum_{s\ge0}\mathfrak h_{r,s}\frac{z^s}{s!},&
\mathcal H_g(y,z)&=\sum_{r\ge0}\mathcal H_{g,r}(z)\frac{y^r}{r!},\\
\mathcal A_g(z)&=\mathcal H_g(0,z),&
\mathcal G_g(t)&=\mathcal H_g(t,0).
\end{align*}
Thus $\mathcal A_g$ is the universal class version of the all-light
series of Section~\ref{sec:all-light-decomposition}, and $\mathcal G_g$ encodes ordinary stable
pointed curves. The preceding finite recursions hold before this
extension of coefficients.

Define the universal coordinate over $\Q[\lambda]$ by
\begin{equation}
\Phi_\lambda(z):=z-\sum_{m\ge2}p_m(\lambda)\frac{z^{m+1}}{(m+1)!},
\qquad
\Psi_\lambda:=\Phi_\lambda^{\langle-1\rangle}.
\label{eq:phi-formal}
\end{equation}
Its linear coefficient is one, so the inverse exists uniquely.
Specializing $\lambda=\LL$ defines $\Phi_{\LL}$ and $\Psi_{\LL}$ in
$\eu K[[z]]$, without inverting any Lefschetz polynomial.
The closed expression
\begin{equation}
\Phi_\lambda(z)=
\frac{(e^{\lambda z}-1)/\lambda-\lambda(e^z-1)}{1-\lambda}
\label{eq:phi-closed}
\end{equation}
is understood coefficientwise: its apparent denominators are removable.
In particular,
\[
\Phi_0(z)=z,\qquad \Phi_1(z)=(2-z)e^z-2.
\]
The Euler coordinates denoted by $\Phi$ and $\Psi$ in
\cite[Section~2.4]{CottiPaperI} are precisely $\Phi_1$ and $\Psi_1$.
We retain the Lefschetz subscript here to distinguish the coordinates
of the different realizations.

\begin{theorem}[Universal motivic heavy/light transformation]
\label{thm:pde}
For every $g\ge1$,
\begin{equation}
\mathcal H_g(y,z)=\mathcal G_g\bigl(y+\Phi_{\LL}(z)\bigr)
\quad\text{in }\eu K[[y,z]].
\label{eq:motivic-composition}
\end{equation}
Equivalently, $\mathcal H_g$ is the unique formal solution of
\begin{equation}
\partial_z\mathcal H_g=\Phi_{\LL}'(z)\partial_y\mathcal H_g,
\qquad \mathcal H_g(y,0)=\mathcal G_g(y).
\label{eq:motivic-pde}
\end{equation}
In particular,
\begin{equation}
\mathcal A_g=\mathcal G_g\circ\Phi_{\LL},\qquad
\mathcal G_g=\mathcal A_g\circ\Psi_{\LL}.
\label{eq:universal-chamber-composition}
\end{equation}
\end{theorem}

\begin{proof}
Multiply \eqref{eq:motivic-onestep} by $z^s/s!$ and sum over $s$.
Binomial convolution gives
\[
\left(1-\sum_{m\ge2}p_m(\LL)\frac{z^m}{m!}\right)
\mathcal H_{g,r}(z)=\mathcal H_{g,r-1}'(z).
\]
The factor in parentheses is $\Phi_{\LL}'(z)$. Summing against
$y^{r-1}/(r-1)!$ gives \eqref{eq:motivic-pde}.
Set $J(y,x)=\mathcal H_g(y,\Psi_{\LL}(x))$. The chain rule gives
$\partial_xJ=\partial_yJ$ and $J(y,0)=\mathcal G_g(y)$.
Over the $\Q$-algebra $\eu K$, the unique solution is
$J(y,x)=\mathcal G_g(y+x)$: if $J=\sum_{n\ge0}J_n(y)x^n$, then
$(n+1)J_{n+1}=J_n'$ with prescribed $J_0$.
Substitution proves \eqref{eq:motivic-composition}; setting $y=0$
and composing with the formal inverse gives the last identities.
\end{proof}

For a motivic measure $\mu:K_0(\Var_{\C})\to\mathsf R$, put
$\lambda=\mu(\LL)$. Coefficientwise realization gives the series
$H_g^\mu,A_g^\mu,G_g^\mu$; explicitly,
\[
G_g^\mu(t)=\sum_{n\ge0}\mu(\overline M_{g,n})\frac{t^n}{n!},
\]
with the same unstable genus-one constant convention. Since the
coordinate coefficients are polynomials in $\LL$, realization
commutes with forming their formal inverse. Therefore
\begin{corollary}
\label{cor:realized-composition}
For every such measure and $g\ge1$,
\begin{equation}
H_g^\mu(y,z)=G_g^\mu(y+\Phi_\lambda(z)),\qquad
G_g^\mu=A_g^\mu\circ\Psi_\lambda.
\label{eq:realized-composition}
\end{equation}
\end{corollary}

For Hodge--Deligne realization $\lambda=uv$; its virtual Poincar\'e
specialization $u=v=-t$ has $\lambda=t^2$. The Euler and virtual
total-Betti specializations both have Lefschetz value one and hence
use the same coordinate $\Phi_1$. The cohomological interpretation
on the proper Hassett spaces will be recalled with the Hodge
applications.
The coefficient transformation, its forest interpretation, and the
associated recurrences are developed in Appendix~\ref{app:triangular}.

\begin{remark}[Relation to the equivariant Hodge formula]
The scalar Hodge--Deligne specialization agrees with the
change-of-variables identity of Kannan--Serpente--Yun
\cite[Corollary~B]{KSY}. Their Theorem~A retains symmetric-group
representations. The theorem here instead establishes the scalar
identity for coarse classes in $K_0(\Var_{\C})\otimes\Q$, before
choosing a realization.
\qrem
\end{remark}

\subsection{Genus zero and the inverse coordinate}
\label{sec:genuszero-coordinate}

There is no stable minimal all-light chamber in genus zero. Nevertheless,
the ordinary genus-zero spaces determine the same inverse coordinate.
Set
\begin{equation}
\mathcal G_0(z)=\sum_{n\ge3}[\overline M_{0,n}]\frac{z^n}{n!},
\qquad
\mathcal G_0^\bullet(z)=1+z+\mathcal G_0'(z).
\label{eq:G0bullet}
\end{equation}
The terms $1+z$ are formal unstable initial terms. The classical
genus-zero recursion for Grothendieck classes is equivalent to
\begin{equation}
(\mathcal G_0^\bullet)'(z)=
\frac{\mathcal G_0^\bullet(z)}
{1+\LL(1+z)-\LL\mathcal G_0^\bullet(z)},
\qquad \mathcal G_0^\bullet(0)=1;
\label{eq:G0bullet-ode}
\end{equation}
see \cite{Manin}, \cite[Proposition~3.2]{ManinMarcolli}, and
\cite{AluffiMarcolliNascimento}. The denominator has
constant term one, so this is a formal identity in $\eu K[[z]]$.

\begin{proposition}[Genus-zero interpretation of the inverse]
\label{prop:Psi-genuszero}
One has
\begin{equation}
\Psi_{\LL}(z)=\log\mathcal G_0^\bullet(z),\qquad
\mathcal G_0^\bullet(z)=e^{\Psi_{\LL}(z)}.
\label{eq:Psi-log-G0}
\end{equation}
\end{proposition}
\begin{proof}
Set $W=\log\mathcal G_0^\bullet$. By \eqref{eq:G0bullet-ode},
\[
W'(z)=\frac1{1+\LL(1+z)-\LL e^{W(z)}},\qquad W(0)=0.
\]
The definition of the coordinate gives the formal identity
\begin{equation}
\Phi_{\LL}'(w)=1+\LL(1+\Phi_{\LL}(w))-\LL e^w.
\label{eq:Phi-genuszero-identity}
\end{equation}
Substituting $w=\Psi_{\LL}(z)$ and differentiating
$\Phi_{\LL}(\Psi_{\LL}(z))=z$ shows that $\Psi_{\LL}$ satisfies
the same differential equation and initial condition as $W$.
Formal coefficient recursion over $\eu K$ gives uniqueness.
\end{proof}

\begin{corollary}[Genus-zero form of the universal transformation]
\label{cor:genuszero-heavy-light}
For every $g\ge1$, one has
\begin{equation}
\mathcal G_g(z)
=\mathcal A_g\!\left(\log\mathcal G_0^\bullet(z)\right)
\quad\text{in }\eu K[[z]].
\label{eq:genuszero-heavy-light}
\end{equation}
Thus the ordinary genus-zero series determines the change of
coordinates from all-light to ordinary markings in every positive genus.
\end{corollary}
\begin{proof}
Substitute Proposition~\ref{prop:Psi-genuszero} into
$\mathcal G_g=\mathcal A_g\circ\Psi_{\LL}$ from
Theorem~\ref{thm:pde}.
\end{proof}

After realization, write
$G_0^{\mu,\bullet}=1+z+(G_0^\mu)'$ for the corresponding pointed
series. Then $G_0^{\mu,\bullet}=e^{\Psi_\lambda}$ and
\begin{equation}
G_g^\mu(z)=A_g^\mu\bigl(\log G_0^{\mu,\bullet}(z)\bigr),
\qquad g\ge1.
\label{eq:realized-genuszero-heavy-light}
\end{equation}
Thus genus zero determines the coordinate that relates the two
positive-genus chambers, without requiring a genus-zero all-light
series.

\subsection{The coordinate and the edge propagator}
\label{sec:propagator-coordinate}

The propagator from Section~\ref{sec:edge-propagator} and the wall coordinate satisfy a
single formal identity. Recall
$R_\lambda(z)=\sum_{m\ge1}(\lambda-1)^{m-1}z^m/m!$ and
$P_\lambda=(1-R_\lambda)^{-1}$.
\begin{proposition}
\label{prop:propagator-coordinate}
In $\Q[\lambda][[z]]$,
\begin{equation}
\Phi_\lambda'(z)=e^z(1-R_\lambda(z)),\qquad
P_\lambda(z)=\frac{e^z}{\Phi_\lambda'(z)}.
\label{eq:PandPhi}
\end{equation}
The same identities hold universally at $\lambda=\LL$ and after
every motivic realization.
\end{proposition}
\begin{proof}
Differentiating \eqref{eq:phi-closed} yields
\[
\Phi_\lambda'(z)=\frac{e^{\lambda z}-\lambda e^z}{1-\lambda}
=e^z\left(1-\frac{e^{(\lambda-1)z}-1}{\lambda-1}\right).
\]
Both sides are understood coefficientwise as polynomials in
$\lambda$. The first identity follows, including at $\lambda=1$;
inverting the unit $1-R_\lambda$ gives the second.
\end{proof}

The edge series therefore links the chain geometry to the heavy/light
wall geometry. For a scalar specialization, wherever these functions
are interpreted analytically, the critical points of $\Phi_\lambda$
are precisely the poles of $P_\lambda$, since the exponential has
no zeros. This supplies the common starting point for the all-light
and ordinary asymptotic analyses of Section~\ref{sec:motivic-asymptotics}. The exceptional
genus-one all-light geometry is treated next.

\section{Genus one: motivic decomposition and Hodge formulas}
\label{sec:g1}

The genus-one case requires a separate construction because forgetting
all markings no longer produces a stable unmarked curve. The all-light
locus instead consists of smooth elliptic curves and rational necklaces,
whose length can grow with the number of labels. We first compute the
necklace contribution in the Grothendieck ring, keeping the coarse
quotients that occur for lengths one and two. Summing over all lengths
then produces a logarithmic term in the edge propagator. The remaining
smooth elliptic contribution can be studied through elliptic local
systems; after Hodge--Deligne realization, their cohomology gives an
explicit series which we resum in terms of Bessel functions.

\subsection{The elliptic locus and rational necklaces}
\label{sec:genus-one-loci}

Let $n\ge1$. As in Section~\ref{sec:all-light-decomposition}, an all-light stable curve has no rational
tails, and every rational two-valent component carries light markings.
If a genus-one curve has a component of positive geometric genus, the
genus formula forces its dual graph, with the marking legs omitted,
to be a tree. Any additional rational
components would produce a rational tail. Thus the curve is smooth
elliptic. If all components are rational, this connected graph has
first Betti number one and no vertices of valence less than two. It is
therefore a circuit, with a nonempty light label set on every component.
This includes a nodal irreducible rational curve and a pair of rational
components meeting at two nodes.

Let $\mathcal I_n$ be the coarse smooth elliptic locus and
$\mathcal B_{n,j}$ the coarse locus of rational necklaces with exactly
$j$ components. We have finite locally closed decompositions
\[
\overline M_{1,(0^+)^n}=\mathcal I_n\sqcup\mathcal B_n,
\qquad \mathcal B_n=\bigsqcup_{j=1}^n\mathcal B_{n,j}.
\]
Their universal series are
\[
\mathcal I(z)=\sum_{n\ge1}[\mathcal I_n]\frac{z^n}{n!},\qquad
\mathcal B(z)=\sum_{n\ge1}[\mathcal B_n]\frac{z^n}{n!},\qquad
\mathcal A_1=\mathcal I+\mathcal B
\]
in $\eu K[[z]]$.
For a measure $\mu$ with values in $\mathsf R$, their realizations are
\[
I^\mu(z):=\sum_{n\ge1}\mu(\mathcal I_n)\frac{z^n}{n!},
\qquad
B^\mu(z):=\sum_{n\ge1}\mu(\mathcal B_n)\frac{z^n}{n!},
\qquad A_1^\mu=I^\mu+B^\mu.
\]

\subsection{Collision partitions on the elliptic locus}
\label{sec:elliptic-motivic-collisions}

Let $M_{1,k}$ be the coarse space of smooth elliptic curves with $k$
distinct labeled markings, and set
$\mathcal F_1(x)=\sum_{k\ge1}[M_{1,k}]x^k/k!$.
\begin{proposition}[Elliptic collision identity]
\label{prop:motivic-elliptic-collisions}
For $n\ge1$,
\begin{equation}
[\mathcal I_n]=\sum_{k=1}^n \genfrac{\{}{\}}{0pt}{}{n}{k}[M_{1,k}],\qquad
\mathcal I(z)=\mathcal F_1(e^z-1).
\label{eq:motivic-elliptic-collisions}
\end{equation}
\end{proposition}
\begin{proof}
The coincidence relation among the light labels gives a set partition
of $[n]$. For a partition with $k$ blocks, order the blocks by their
least labels. The corresponding locus is isomorphic to $M_{1,k}$;
there is no further quotient by $\mathfrak S_k$. Additivity gives the
coefficient formula, and the Stirling exponential generating function
gives the composition.
\end{proof}

There is also a useful relation with the vertex series of Section~\ref{sec:all-light-decomposition}.
Use the first marking as a distinguished point. Among the other $m$
markings, choose the $k$ that do not coincide with it. These are light
markings on the punctured elliptic curve, allowed to coincide among
themselves, and are parametrized by $\mathcal V_{1,1,k}$. Therefore
\begin{equation}
\mu(\mathcal I_{m+1})=\sum_{k=0}^m\binom{m}{k}
\mu(\mathcal V_{1,1,k}),\qquad
(I^\mu)'(z)=e^z V_{1,1}^\mu(z),\quad I^\mu(0)=0.
\label{eq:I-vertex-series}
\end{equation}
In Section~\ref{sec:motivic-asymptotics} this identity will control the regular genus-one
contribution without requiring an explicit Hodge formula.

\subsection{Necklace quotients and the motivic decomposition}
\label{sec:motivic-genus-one}

The short necklaces require an actual coarse quotient calculation.
\begin{lemma}[Simultaneous inversion]
\label{lem:torus-inversion-quotient}
Let $X_m=\GG^m/\langle\iota_m\rangle$, where
\[
\iota_m(x_1,\ldots,x_m)=(x_1^{-1},\ldots,x_m^{-1}),
\]
and set
$X_0=\operatorname{pt}$. Then
\begin{equation}
[X_m]=\frac{(\LL+1)^m+(\LL-1)^m}{2}
\quad\text{in }K_0(\Var_{\C}).
\label{eq:torus-inversion-class}
\end{equation}
The right-hand side belongs to $\Z[\LL]$, so no division by two in
the Grothendieck ring is required.
\end{lemma}
\begin{proof}
Compactify to $\overline X_m=(\mathbb P^1)^m/\langle\iota_m\rangle$.
For $m\ge1$, projection to the first factor induces a map to
$\mathbb P^1/\langle x\mapsto x^{-1}\rangle\cong\mathbb P^1$.
Away from its two branch values, corresponding to $x=\pm1$, the
quotient is trivial with fiber $(\mathbb P^1)^{m-1}$. Indeed, for
$x=x_1$ and $y=x_j$, the projective coordinate
\[
u_j=\frac{y-x}{xy-1}
\]
is invariant under $(x,y)\mapsto(x^{-1},y^{-1})$ and gives a
M\"obius coordinate whenever $x^2\ne1$. These coordinates give the
trivialization on the quotient; their inverses are regular in the
corresponding projective charts. The open base has class $\LL-1$,
and the reduced fiber over either branch value is $\overline X_{m-1}$.
Consequently
\[
[\overline X_m]=(\LL-1)(\LL+1)^{m-1}+2[\overline X_{m-1}]
=(\LL+1)^m
\]
by induction from $\overline X_0=\operatorname{pt}$.

To remove the boundary, choose the $k\ge1$ coordinates lying in
$\{0,\infty\}$. Their $2^k$ value assignments are paired freely by
simultaneous inversion; each pair contributes $\GG^{m-k}$.
Thus
\[
[\overline X_m]-[X_m]
=\sum_{k=1}^m\binom{m}{k}2^{k-1}(\LL-1)^{m-k}
=\frac{(\LL+1)^m-(\LL-1)^m}{2}.
\]
Subtracting proves the claim, including $m=0$ directly. The resulting
integer polynomial is
\[
\sum_{\substack{0\le r\le m\\r\ \mathrm{even}}}
\binom{m}{r}\LL^{m-r}.
\]
\end{proof}

\begin{figure}[htbp]
\centering
\begin{tikzpicture}[>=Stealth,line width=.65pt,
 every node/.style={font=\small},
 graph/.style={circle,draw,line width=.65pt,minimum size=4mm,
 inner sep=0pt,fill=white,font=\scriptsize}]
\begin{scope}
\node at (1.4,3) {\textbf{(a) One component}};
\node[graph] (v) at (1.4,0) {$0$};
\draw (v) .. controls (-.1,2.3) and (2.9,2.3) .. (v);
\draw[<->,densely dashed] (.98,.46) to[bend left=20] (1.82,.46);
\node at (1.4,-.55) {$B_1=[n]$};
\node[font=\footnotesize] at (1.4,-1.15) {branches exchanged};
\node at (1.4,-1.75) {$\mathbf x\longmapsto\mathbf x^{-1}$};
\node at (1.4,-2.35) {$\GG^{n-1}/C_2$};
\end{scope}
\begin{scope}[xshift=5cm]
\node at (1.4,3) {\textbf{(b) Two components}};
\node[graph] (a) at (0,.8) {$0$};
\node[graph] (b) at (2.8,.8) {$0$};
\draw (a) to[bend left=65] node[above] {$q_1$} (b);
\draw (a) to[bend right=65] node[below] {$q_2$} (b);
\draw[<->,densely dashed] (1.4,1.48)--(1.4,.12);
\node at (0,.25) {$B_1$};
\node at (2.8,.25) {$B_2$};
\node[font=\footnotesize] at (1.4,-1.15) {nodes exchanged};
\node at (1.4,-1.75) {$\mathbf x\longmapsto\mathbf x^{-1}$};
\node at (1.4,-2.35) {$\GG^{n-2}/C_2$};
\end{scope}
\begin{scope}[xshift=10cm]
\node at (1.4,3) {\textbf{(c) At least three}};
\node[graph] (a) at (0,0) {$0$};
\node[graph] (b) at (2.8,0) {$0$};
\node[graph] (c) at (1.4,1.8) {$0$};
\draw (a)--(b)--(c)--(a);
\node at (0,-.55) {$B_1$};
\node at (2.8,-.55) {$B_2$};
\node[above=3pt] at (c) {$B_3$};
\node[font=\footnotesize] at (1.4,-1.15) {label blocks fixed};
\node at (1.4,-1.75) {trivial stabilizer};
\node at (1.4,-2.35) {$\GG^{n-j}$};
\end{scope}
\end{tikzpicture}
\caption{The three rational-necklace regimes, shown as dual graphs.
White vertices are rational components; edges represent nodes, and
$B_i$ are nonempty disjoint light-label blocks. For one component,
exchanging the two branches of the self-node induces simultaneous
inversion on the component coordinates. For two components,
exchanging the two nodes induces simultaneous inversion on their
product, while fixing each label block. For $j\ge3$ (a three-cycle
is drawn), fixing every light label fixes the entire graph.
The spaces below the panels are the coarse parameter spaces for a
fixed label allocation and cyclic arrangement; the multiplicities
from choosing these data are counted in
\eqref{eq:necklace-one-class}--\eqref{eq:motivic-necklace-class}.
The arrows in (a) and (b) indicate the involutions of the graphs.}
\label{fig:necklace-short-cycle-symmetries}
\end{figure}
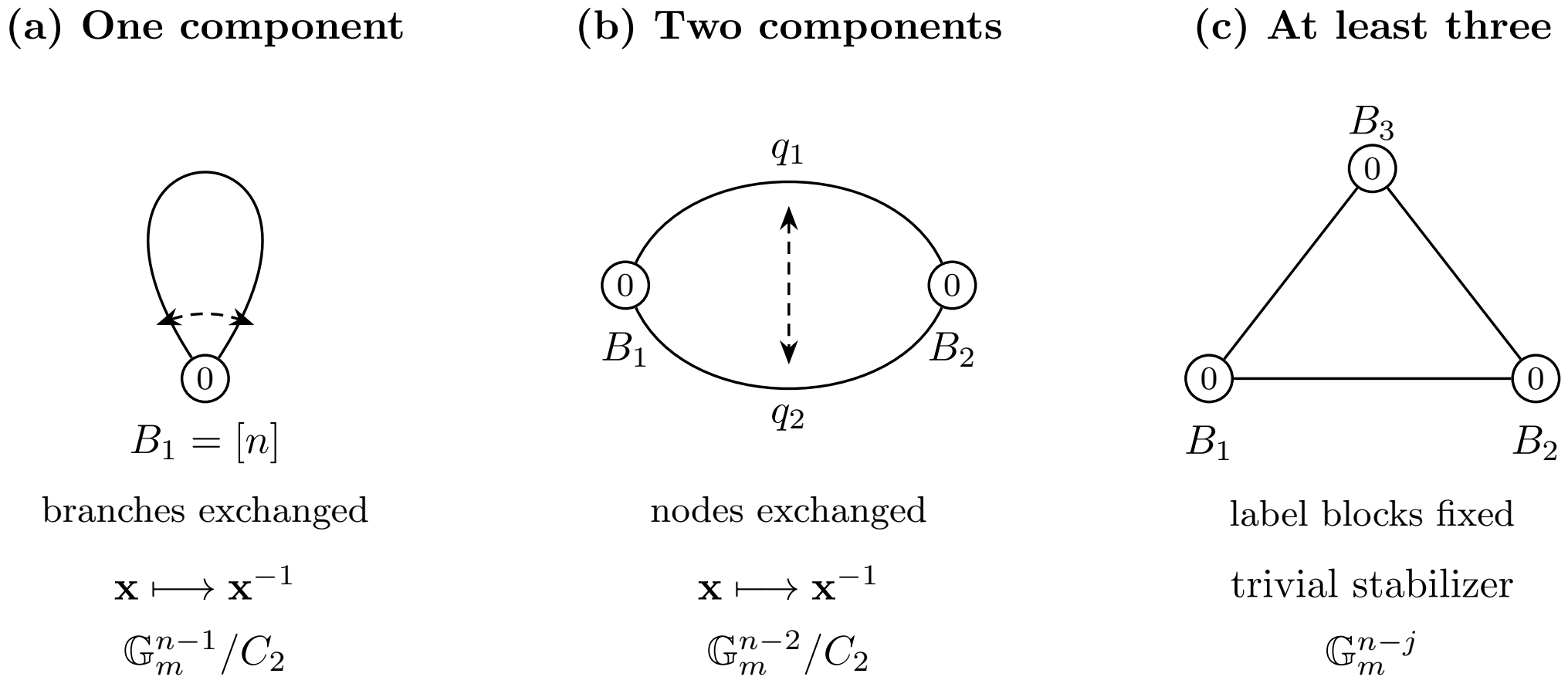

For one component, normalize the two branches of the node to $0,\infty$
and one light marking to $1$. Exchanging the branches acts by simultaneous
inversion on the other $n-1$ coordinates, so $\mathcal B_{n,1}\cong X_{n-1}$.
For two components, an unordered bipartition of the light labels determines
which component carries each nonempty block; there are $2^{n-1}-1$ choices.
Exchanging the two nodes induces simultaneous inversion on the product
of the two component tori. Each stratum is therefore $X_{n-2}$.
For $j\ge3$, the nonempty disjoint label blocks prevent nonidentity graph
automorphisms from fixing the decoration. A cyclic order up to reversal
gives $(j-1)!/2$ arrangements per set partition. Hence
\begin{align}
[\mathcal B_{n,1}]&=
\frac{(\LL+1)^{n-1}+(\LL-1)^{n-1}}2,
\label{eq:necklace-one-class}\\
[\mathcal B_{n,2}]&=(2^{n-1}-1)
\frac{(\LL+1)^{n-2}+(\LL-1)^{n-2}}2\quad(n\ge2),
\label{eq:necklace-two-class}\\
[\mathcal B_{n,j}]&=\genfrac{\{}{\}}{0pt}{}{n}{j}\frac{(j-1)!}{2}(\LL-1)^{n-j}
\quad(3\le j\le n).
\label{eq:motivic-necklace-class}
\end{align}
The $j=2$ term is absent for $n=1$. The factors $1/2$ here represent
integer polynomials or integer counts, not a motivic fixed-locus average.
The three discrete symmetry types and their parameter spaces are
summarized in Figure~\ref{fig:necklace-short-cycle-symmetries}.

Define the companion series
\begin{equation}
U_\lambda(z)=\sum_{m\ge1}(\lambda+1)^{m-1}\frac{z^m}{m!}
=\frac{e^{(\lambda+1)z}-1}{\lambda+1},
\label{eq:U-motivic}
\end{equation}
where the series defines its value also at $\lambda=-1$. The three
necklace formulas give
\[
\mathcal B=
\frac{R_{\LL}+U_{\LL}}2+
\frac{R_{\LL}^2+U_{\LL}^2}4+
\sum_{j\ge3}\frac{R_{\LL}^j}{2j}.
\]
Since $\log P_{\LL}=\sum_{j\ge1}R_{\LL}^j/j$, this proves:
\begin{proposition}[Motivic genus-one decomposition]
\label{prop:motivic-genus-one}
In $\eu K[[z]]$,
\begin{equation}
\mathcal A_1=\mathcal I+\frac12\log P_{\LL}
+\frac12 U_{\LL}+\frac14 U_{\LL}^2.
\label{eq:A1-motivic-complete}
\end{equation}
For every motivic measure with Lefschetz value $\lambda$,
\begin{equation}
A_1^\mu(z)=\frac12\log P_\lambda(z)+S_1^\mu(z),\qquad
S_1^\mu=I^\mu+\frac12U_\lambda+\frac14U_\lambda^2.
\label{eq:A1mu-complete}
\end{equation}
\end{proposition}
The logarithmic term records the unbounded necklace length. At this
stage $S_1^\mu$ is a formal series; its analytic regularity is not being
assumed. It will follow from vertex bounds for the realizations used
in Section~\ref{sec:motivic-asymptotics}.

\subsection{Hodge--Deligne realization and Bessel resummation}
\label{sec:bessel}

We use the compactly supported Hodge--Deligne polynomial
\[
E_c(X;u,v)=\sum_{k,p,r}(-1)^k
\dim\operatorname{Gr}_F^p\operatorname{Gr}_{p+r}^W H_c^k(X;\C)
\,u^p v^r,
\]
and put $q=uv$. It is additive and multiplicative, and sends $\LL$ to
$q$~\cite{DeligneII,DeligneIII}. Rational compactly supported cohomology
of a separated finite-type Deligne--Mumford stack agrees with that of
its coarse space, compatibly with mixed Hodge structures; we use this
comparison for the component stacks~\cite{Behrend}.

Specializing the series of Section~\ref{sec:genus-one-loci} to
$\mu=E_c(-;u,v)$ gives
\[
I(z;u,v):=I^{E_c}(z)=\sum_{n\ge1}E_c(\mathcal I_n;u,v)\frac{z^n}{n!},
\qquad
B(z;q):=B^{E_c}(z)=\sum_{n\ge1}E_c(\mathcal B_n;u,v)\frac{z^n}{n!}.
\]
Here $\mathcal I_n$ is the smooth elliptic locus and $\mathcal B_n$
the rational-necklace boundary. The boundary realization depends on
$u,v$ only through $q=uv$. Specializing the boundary contribution
$A_1^\mu-I^\mu$ in \eqref{eq:A1mu-complete} gives
\begin{equation}
B(z;q)=\frac12\log P_q(z)+\frac12U_q(z)+\frac14U_q(z)^2.
\label{eq:B}
\end{equation}
This agrees with the Hodge boundary formula of
\cite[Proposition~4.3]{Kiem}.

For the interior, we proceed as in \cite{Kiem}. Use the first marking
as origin. The stack
$\mathscr I_n$ with coarse space $\mathcal I_n$ is the $(n-1)$-fold
fiber product of the universal elliptic curve over $\mathcal M_{1,1}$.
Let $\mathbb V=R^1p_*\Q$ be its rank-two local system.
Odd tensor powers have no cohomology over the base because the
generic elliptic involution acts by minus one. For the even powers,
the Clebsch--Gordan rule gives
\[
\mathbb V^{\otimes2j}\cong\bigoplus_{i=0}^j
\alpha_{ij}\,\Sym^{2j-2i}\mathbb V\otimes\Q(-i),\qquad
\alpha_{ij}=\binom{2j}{i}-\binom{2j}{i-1}.
\]
The Tate twist contributes $q^i$. Set
$\beta_m=\dim_\C S_m(\mathrm{SL}_2(\Z))$; in particular $\beta_2=0$.
Their generating function is recalled in
Lemma~\ref{lem:cusp-dimensions} of
Appendix~\ref{app:poincare-resummation}.
Kiem's local-system calculation, using Eichler--Shimura, yields
\begin{align}
E_c(\mathcal I_n;u,v)&=
\sum_{j=0}^{\lfloor(n-1)/2\rfloor}\binom{n-1}{2j}
(1+q)^{n-1-2j}\mathcal C_j(u,v),
\label{eq:Kiem-interior}\\
\mathcal C_j(u,v)&=\alpha_{jj}q^{j+1}
-\sum_{i=0}^{j-1}\alpha_{ij}q^i
\left[1+\beta_{2j-2i+2}
(u^{2j-2i+1}+v^{2j-2i+1})\right].
\label{eq:Cj}
\end{align}
Here the binomial coefficient chooses the $2j$ factors contributing
$H^1$, and the remaining elliptic factors contribute $(1+q)^{n-1-2j}$;
the full cohomological input is \cite[Proposition~4.1]{Kiem}.

Define
\begin{equation}
\mathcal J_k(z;q)=\sum_{i\ge0}
\frac{(2k+1)q^i z^{2i+2k}}{i!(i+2k+1)!}
=\frac{2k+1}{zq^{k+1/2}}\mathsf I_{2k+1}(2z\sqrt q),
\label{eq:J-def}
\end{equation}
where $\mathsf I_\nu$ is the modified Bessel function of the first
kind~\cite{Watson,NIST}. The power series defines the expression also
at $q=0$ and $z=0$, independently of the choice of square root.

\begin{theorem}[Genus-one Hodge--Deligne series]
\label{thm:A1explicit}
For $q=uv$,
\begin{equation}
A_1^{E_c}(z;u,v)=B(z;q)+I(z;u,v),
\label{eq:A1-Hodge}
\end{equation}
with $B$ given by \eqref{eq:B} and
\begin{equation}
I(z;u,v)=\int_0^z e^{(1+q)\zeta}
\left(q\mathcal J_0(\zeta;q)-
\sum_{k\ge1}[1+\beta_{2k+2}(u^{2k+1}+v^{2k+1})]
\mathcal J_k(\zeta;q)\right)d\zeta.
\label{eq:bessel}
\end{equation}
\end{theorem}
\begin{proof}
Differentiate the interior EGF and use \eqref{eq:Kiem-interior}:
\[
I'(z;u,v)=e^{(1+q)z}\sum_{j\ge0}
\mathcal C_j(u,v)\frac{z^{2j}}{(2j)!}.
\]
The diagonal term uses
$\alpha_{jj}/(2j)!=1/(j!(j+1)!)$, giving $q\mathcal J_0$.
For the other terms put $j=i+k$, $k\ge1$, and use
\[
\frac{\alpha_{i,i+k}}{(2i+2k)!}
=\frac{2k+1}{i!(i+2k+1)!}.
\]
The sum over $i$ is $\mathcal J_k$, while the bracket in
\eqref{eq:Cj} depends only on $k$. Integration with $I(0)=0$
gives \eqref{eq:bessel}. All rearrangements are coefficientwise
finite since $\mathcal J_k$ starts in degree $2k$.
\end{proof}

The universal collision identity and boundary calculation precede
realization; the explicit cuspidal coefficients above come from the
Hodge calculation and are not asserted to define classes in
$K_0(\Var_{\C})$.

\subsection{Poincar\'e, Euler, and total Betti specializations}
\label{sec:poincare1}

We now specialize the Hodge series.

For a general variety $X$, $E_c(X;-t,-t)$ is a virtual Poincar\'e
polynomial. For the coarse space of a smooth proper Hassett stack,
rational cohomology agrees with that of the stack and is pure.
Thus in this case
\[
P(X;t)=E_c(X;-t,-t)=\sum_k b_k(X)t^k.
\]
In particular $P(X;-1)=\chi(X)$ and $P(X;1)=\sum_k b_k(X)$.
These interpretations apply to the complete spaces below, not to
each nonproper stratum separately.

Set $A_1^P(z;t)=A_1^{E_c}(z;-t,-t)$ and $I_P(z;t)=I(z;-t,-t)$.
Then
\begin{equation}
A_1^P(z;t)=B(z;t^2)+I_P(z;t),\qquad
I_P'(z;t)=e^{(1+t^2)z}\mathscr S(z;t),\quad I_P(0;t)=0,
\label{eq:A1P-decomposition}
\end{equation}
where
\begin{equation}
\mathscr S(z;t)=t^2\mathcal J_0(z;t^2)
-\sum_{k\ge1}(1-2\beta_{2k+2}t^{2k+1})\mathcal J_k(z;t^2).
\label{eq:SP-original}
\end{equation}
Define the unilateral even-Bessel series by
\begin{equation}
\mathcal N(z;t)=\sum_{k\ge1}t^{-2k}\mathsf I_{2k}(2tz)
=\sum_{i\ge0,k\ge1}\frac{t^{2i}z^{2i+2k}}{i!(i+2k)!}
\in\Q[t^2][[z]].
\label{eq:N-def}
\end{equation}
The second expression shows that the apparent negative powers of $t$
are removable. An integral representation, including the values at
$t=\pm1$, is given in Proposition~\ref{prop:N-integral} in
Appendix~\ref{app:poincare-resummation}.

\begin{theorem}[Poincar\'e resummation]
\label{thm:poincare-resummation}
The even and odd parts of $\mathscr S$ in $t$ are
\begin{align}
\mathscr S_{\mathrm{ev}}(z;t)&=t^2\mathsf I_0(2tz)
-(1+t^2)\mathsf I_2(2tz)-(1-t^2)\mathcal N(z;t),
\label{eq:SP-even}\\
\mathscr S_{\mathrm{odd}}(z;t)&=
t(\mathsf I_0(2tz)-2\mathsf I_2(2tz))
+\frac{t}{6}(\cosh(2tz)-4\cosh(tz)-3).
\label{eq:SP-odd}
\end{align}
Consequently,
\begin{equation}
\begin{aligned}
A_1^P(z;t)={}&B(z;t^2)+\int_0^z e^{(1+t^2)\zeta}
\Bigl[t(t+1)\mathsf I_0(2t\zeta)-(t+1)^2\mathsf I_2(2t\zeta)\\
&\hspace{15mm}-(1-t^2)\mathcal N(\zeta;t)
+\frac{t}{6}(\cosh(2t\zeta)-4\cosh(t\zeta)-3)\Bigr]d\zeta.
\end{aligned}
\label{eq:A1P-resummed}
\end{equation}
\end{theorem}
The proof is given in Appendix~\ref{app:poincare-resummation}.
The boundary series is even in $t$, so
\[
A_1^{P,\mathrm{ev}}=B(z;t^2)+\int_0^z
e^{(1+t^2)\zeta}\mathscr S_{\mathrm{ev}}(\zeta;t)d\zeta,
\qquad
A_1^{P,\mathrm{odd}}=\int_0^z
e^{(1+t^2)\zeta}\mathscr S_{\mathrm{odd}}(\zeta;t)d\zeta.
\]
In particular, the odd cohomology is entirely supplied by the cuspidal
terms, although the proper space and its interior have different
cohomological interpretations.

At the Euler specialization $t=-1$, all Bessel terms cancel and
\begin{equation}
\mathscr S(z;-1)=\frac12+\frac23\cosh z-\frac16\cosh(2z),
\qquad
(I^\chi)'(z)=e^{2z}\mathscr S(z;-1).
\label{eq:SP-Euler}
\end{equation}
This recovers the Euler formula and the starting point of the full
asymptotic expansion in \cite[Sections~3.3--3.4]{CottiPaperI}.
At $t=1$, the total-Betti series is
\begin{equation}
\begin{aligned}
A_1^b(z)={}&-\frac12\log(1-z)-\frac{e^z-1}{3}
-\frac{e^{2z}-1}{8}-\frac{e^{3z}-1}{9}
+\frac{e^{4z}-1}{12}+\frac z{12}\\
&+2\int_0^z e^{2\zeta}
\bigl(\mathsf I_0(2\zeta)-2\mathsf I_2(2\zeta)\bigr)d\zeta.
\end{aligned}
\label{eq:A1-total-Betti}
\end{equation}
Both specializations share the logarithmic boundary term. Their
leading asymptotics will follow in Section~\ref{sec:motivic-asymptotics} from control of the
remaining series. For either realization, and for the full
Hodge--Deligne series, the ordinary chamber is obtained using the
heavy/light composition of Section~\ref{sec:motivic-heavy-light}.

\section{Analytic asymptotics and Hodge tameness}
\label{sec:motivic-asymptotics}

The decompositions of Sections~\ref{sec:motivicgraphs} and~\ref{sec:g1} hold as formal identities, so
an asymptotic analysis requires additional control of the series that
occur in their remainders. We begin with the singularities determined
by the edge propagator and the heavy/light coordinate when the
Lefschetz value is positive. Once the necessary continuation and
remainder bounds have been stated for a scalar realization, these
singularities yield the coefficient asymptotics in both chambers.
The final part of the section gives geometric sufficient conditions
for the bounds and verifies them for every fixed Hodge--Deligne
evaluation.

A \emph{scalar motivic realization} is a unital ring homomorphism
$\mu:K_0(\Var_{\C})\to\C$. Throughout the main section we assume
$\lambda=\mu(\LL)>0$. The cases $\lambda=0$, $\lambda<0$ and complex
$\lambda$ require different singularity arguments and are treated
separately in Appendix~\ref{app:nonpositive-lambda}.

\subsection{The universal singularities and analytic composition}
\label{sec:universal-analytic-geometry}

For $r>0$, $\eta>0$ and $0<\phi<\pi/2$, set
\begin{equation}
\Delta(r;\eta,\phi)=
\{z\in\C:|z|<r+\eta,\ z\ne r,\ |\arg(z-r)|>\phi\}.
\label{eq:Delta-domain}
\end{equation}
Arguments are principal, and continuation always means continuation
of the germ represented by the series at the origin. Define
\begin{equation}
\sigma_\lambda=\frac{\log\lambda}{\lambda-1},\quad \sigma_1=1,
\qquad \rho_\lambda=\Phi_\lambda(\sigma_\lambda),\qquad
a_\lambda=\sqrt{\frac{2\rho_\lambda}{\lambda e^{\sigma_\lambda}}}.
\label{eq:singularity-constants}
\end{equation}

\begin{proposition}[Universal singularity pair]
\label{prop:universal-singularity-pair}
For $\lambda>0$, the propagator has a unique pole of minimal modulus,
at $z=\sigma_\lambda$. It is simple, with local expansion
\begin{equation}
P_\lambda(z)=\frac1{\lambda\sigma_\lambda}
(1-z/\sigma_\lambda)^{-1}+O(1).
\label{eq:Plambda-local}
\end{equation}
Moreover,
\begin{equation}
\begin{aligned}
\Phi_\lambda'(\sigma_\lambda)=0,\quad
&\Phi_\lambda''(\sigma_\lambda)=-\lambda e^{\sigma_\lambda},\\
\rho_\lambda=e^{\sigma_\lambda}-1-\lambda^{-1},\quad
&0<\rho_\lambda<\sigma_\lambda.
\end{aligned}
\label{eq:critical-data}
\end{equation}
The inverse germ $\Psi_\lambda(0)=0$ has radius $\rho_\lambda$ and
unique dominant singularity $t=\rho_\lambda$. It extends to a
$\Delta$-domain there and has a convergent local Puiseux expansion
beginning with
\begin{equation}
\Psi_\lambda(t)=\sigma_\lambda-a_\lambda
(1-t/\rho_\lambda)^{1/2}+O(1-t/\rho_\lambda).
\label{eq:Psi-lambda-local}
\end{equation}
\end{proposition}
\begin{proof}
For $\lambda\ne1$, the propagator poles are
\[
z_k=\frac{\log\lambda+2\pi i k}{\lambda-1},\qquad k\in\Z.
\]
The pole $z_0=\sigma_\lambda$ is positive, and every other one
has strictly larger modulus. Expanding the denominator of
$P_\lambda=(\lambda-1)/(\lambda-e^{(\lambda-1)z})$ gives
\eqref{eq:Plambda-local}. The case $\lambda=1$ follows from
$P_1=(1-z)^{-1}$. The identity $P_\lambda=e^z/\Phi_\lambda'$ and
direct differentiation of the coordinate give the first three
relations in \eqref{eq:critical-data}, including their limits at one.

Write
\[
\Phi_\lambda(w)=w-H_\lambda(w),\qquad
H_\lambda(w)=\sum_{m\ge2}p_m(\lambda)\frac{w^{m+1}}{(m+1)!}.
\]
Here $H_\lambda$ is entire and has positive coefficients in every
degree at least three. Since $\sigma_\lambda$ is the first positive
zero of $1-H_\lambda'$, the coordinate is strictly increasing on
$[0,\sigma_\lambda]$. Thus $\rho_\lambda>0$, whereas
$\rho_\lambda=\sigma_\lambda-H_\lambda(\sigma_\lambda)<\sigma_\lambda$.

To solve $w=t+H_\lambda(w)$, set
\[
w_0(t)=0,\qquad w_{j+1}(t)=t+H_\lambda(w_j(t)).
\]
The iterates increase coefficientwise and stabilize in each degree
because $H_\lambda$ starts in degree three. Their formal limit is
$\Psi_\lambda(t)=\sum_{k\ge1}a_kt^k$ with $a_k\ge0$. For
$0\le t\le\rho_\lambda$, induction gives
\[
0\le w_j(t)\le w_{j+1}(t)
\le\rho_\lambda+H_\lambda(\sigma_\lambda)=\sigma_\lambda.
\]
By monotone convergence of the nonnegative coefficient sums,
\[
\sum_{k\ge1}a_k\rho_\lambda^k
=\lim_{j\to\infty}w_j(\rho_\lambda)\le\sigma_\lambda.
\]
Thus the inverse series converges absolutely on $|t|\le\rho_\lambda$.
Passing to the limit in the iteration at $\rho_\lambda$ gives
$\Phi_\lambda(\Psi_\lambda(\rho_\lambda))=\rho_\lambda$.
Strict increase of $\Phi_\lambda$ on $[0,\sigma_\lambda]$ then gives
$\Psi_\lambda(\rho_\lambda)=\sigma_\lambda$.

For $F(t,w)=w-t-H_\lambda(w)$ we have
\[
\partial_w F(t,w)=1-H_\lambda'(w)\ne0\qquad(|w|<\sigma_\lambda),
\]
since $|H_\lambda'(w)|\le H_\lambda'(|w|)<1$.
On $|t|=\rho_\lambda$ the inequality
$|\Psi_\lambda(t)|\le\sigma_\lambda$ is strict unless $t=\rho_\lambda$:
equality would align the positive terms in degrees $1,3,4$, forcing
both $(t/\rho_\lambda)^2$ and $(t/\rho_\lambda)^3$ to equal one.
The analytic implicit-function theorem consequently extends the germ
across every other point of the circle.

At $t=\rho_\lambda$ the critical point is simple. Taylor expansion
gives $t-\rho_\lambda=-\tfrac12\lambda e^{\sigma_\lambda}
(w-\sigma_\lambda)^2+O((w-\sigma_\lambda)^3)$.
Local inversion gives a convergent square-root expansion, with the
negative root chosen because $\Psi_\lambda(t)<\sigma_\lambda$ on
$0<t<\rho_\lambda$. Its coefficient is $a_\lambda$.
This is a genuine singularity, and the local expansion together with
continuation across the rest of the circle gives a $\Delta$-domain.
\end{proof}

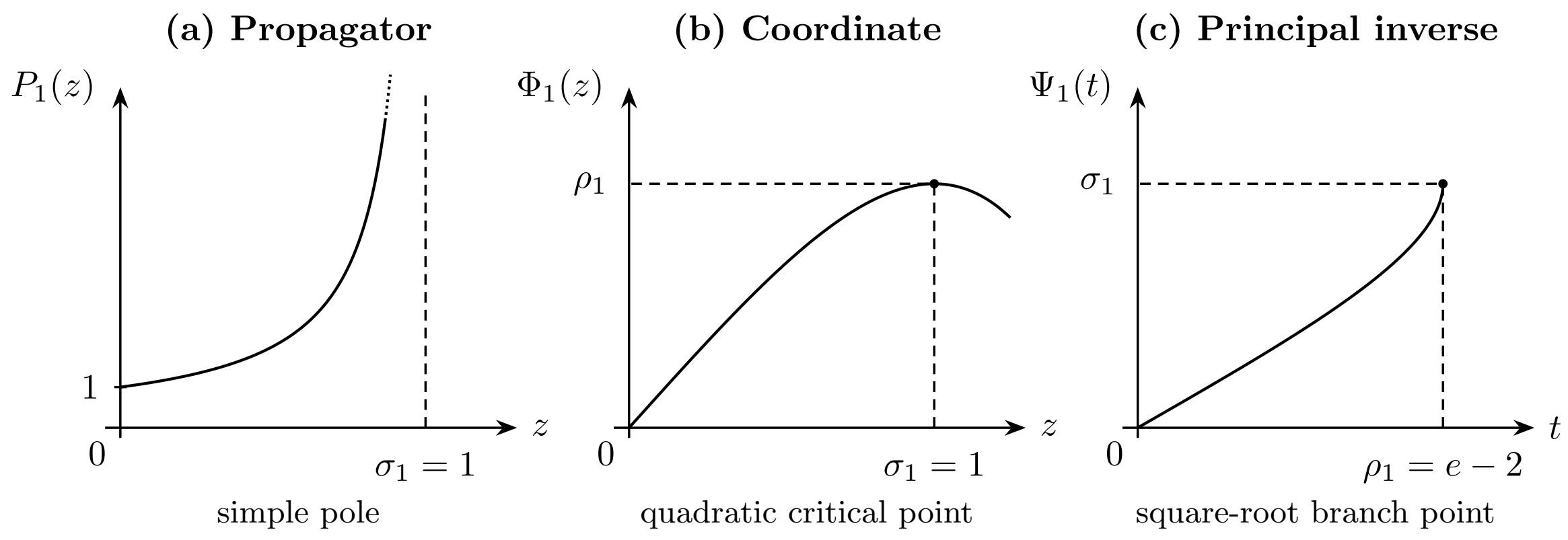
\begin{figure}[htbp]
\centering
\begin{tikzpicture}[>=Stealth,line width=.65pt,
 every node/.style={font=\small}]
\pgfmathsetmacro{\rhone}{exp(1)-2}
\begin{scope}
\node at (1.75,3.9) {\textbf{(a) Propagator}};
\draw[->] (-.15,0)--(3.9,0) node[right] {$z$};
\draw[->] (0,-.1)--(0,3.35) node[left=3pt] {$P_1(z)$};
\node[below left] at (0,0) {$0$};
\draw[densely dashed] (3,0)--(3,3.3);
\node[below=3pt] at (3,0) {$\sigma_1=1$};
\draw (-.06,.4)--(.06,.4);
\node[left=2pt] at (0,.4) {$1$};
\draw[thick,domain=0:.87,samples=100,variable=\x]
 plot ({3*\x},{.4/(1-\x)});
\draw[thick,densely dotted,domain=.87:.885,samples=40,variable=\x]
 plot ({3*\x},{.4/(1-\x)});
\node[font=\footnotesize] at (1.75,-.85) {simple pole};
\end{scope}
\begin{scope}[xshift=5cm]
\node at (1.75,3.9) {\textbf{(b) Coordinate}};
\draw[->] (-.15,0)--(3.9,0) node[right] {$z$};
\draw[->] (0,-.1)--(0,3.35) node[left=3pt] {$\Phi_1(z)$};
\node[below left] at (0,0) {$0$};
\draw[densely dashed] (3,0)--(3,2.4)--(0,2.4);
\node[below=3pt] at (3,0) {$\sigma_1=1$};
\node[left=2pt] at (0,2.4) {$\rho_1$};
\draw[thick,domain=0:1.25,samples=100,variable=\x]
 plot ({3*\x},{2.4*((2-\x)*exp(\x)-2)/\rhone});
\fill (3,2.4) circle (1.3pt);
\node[font=\footnotesize] at (1.75,-.85) {quadratic critical point};
\end{scope}
\begin{scope}[xshift=10cm]
\node at (1.75,3.9) {\textbf{(c) Principal inverse}};
\draw[->] (-.15,0)--(3.9,0) node[right] {$t$};
\draw[->] (0,-.1)--(0,3.35) node[left=3pt] {$\Psi_1(t)$};
\node[below left] at (0,0) {$0$};
\draw[densely dashed] (3,0)--(3,2.4)--(0,2.4);
\node[below=3pt] at (3,0) {$\rho_1=e-2$};
\node[left=2pt] at (0,2.4) {$\sigma_1$};
\draw[thick,domain=0:1,samples=130,variable=\x]
 plot ({3*((2-\x)*exp(\x)-2)/\rhone},{2.4*\x});
\fill (3,2.4) circle (1.3pt);
\node[font=\footnotesize] at (1.75,-.85) {square-root branch point};
\end{scope}
\end{tikzpicture}
\caption{The universal singularity mechanism, illustrated by the actual
real curves for $\lambda=1$. The propagator $P_1(z)=(1-z)^{-1}$
has a simple pole at $\sigma_1=1$; the coordinate
$\Phi_1(z)=(2-z)e^z-2$ is regular there, with a simple critical
point and critical value $\rho_1=e-2$. Its principal inverse has a
square-root singularity at $\rho_1$, with finite limiting value
$\sigma_1$. Axis scales differ between panels. For every
$\lambda>0$, the same local mechanism holds with the constants of
\eqref{eq:singularity-constants}.}
\label{fig:universal-real-geometry}
\end{figure}

Figure~\ref{fig:universal-real-geometry} illustrates the real
singularity mechanism in the Euler case $\lambda=1$.

\begin{lemma}[Analyticity under composition]
\label{lem:analytic-composition}
If $F$ is analytic in a domain $\Delta(\sigma_\lambda;\eta,\phi)$,
then $F\circ\Psi_\lambda$ is analytic in a suitable $\Delta$-domain
at $\rho_\lambda$.
\end{lemma}
\begin{proof}
For $|t|<\rho_\lambda$, positivity gives
$|\Psi_\lambda(t)|\le\Psi_\lambda(|t|)<\sigma_\lambda$.
The preceding proof gives strict inequality also on the boundary
circle away from $\rho_\lambda$. Continuation across these points
therefore keeps the image inside the domain of $F$.
Near $\rho_\lambda$, the expansion
$\Psi_\lambda(t)-\sigma_\lambda\sim
-a_\lambda(1-t/\rho_\lambda)^{1/2}$ maps a sufficiently small
dented neighborhood into the allowed sector for $F$: the argument of
this difference stays away from the excluded positive-real cone,
whose half-angle is less than $\pi/2$. Compactness of the rest of
the boundary circle and this local inclusion give the claimed domain.
\end{proof}

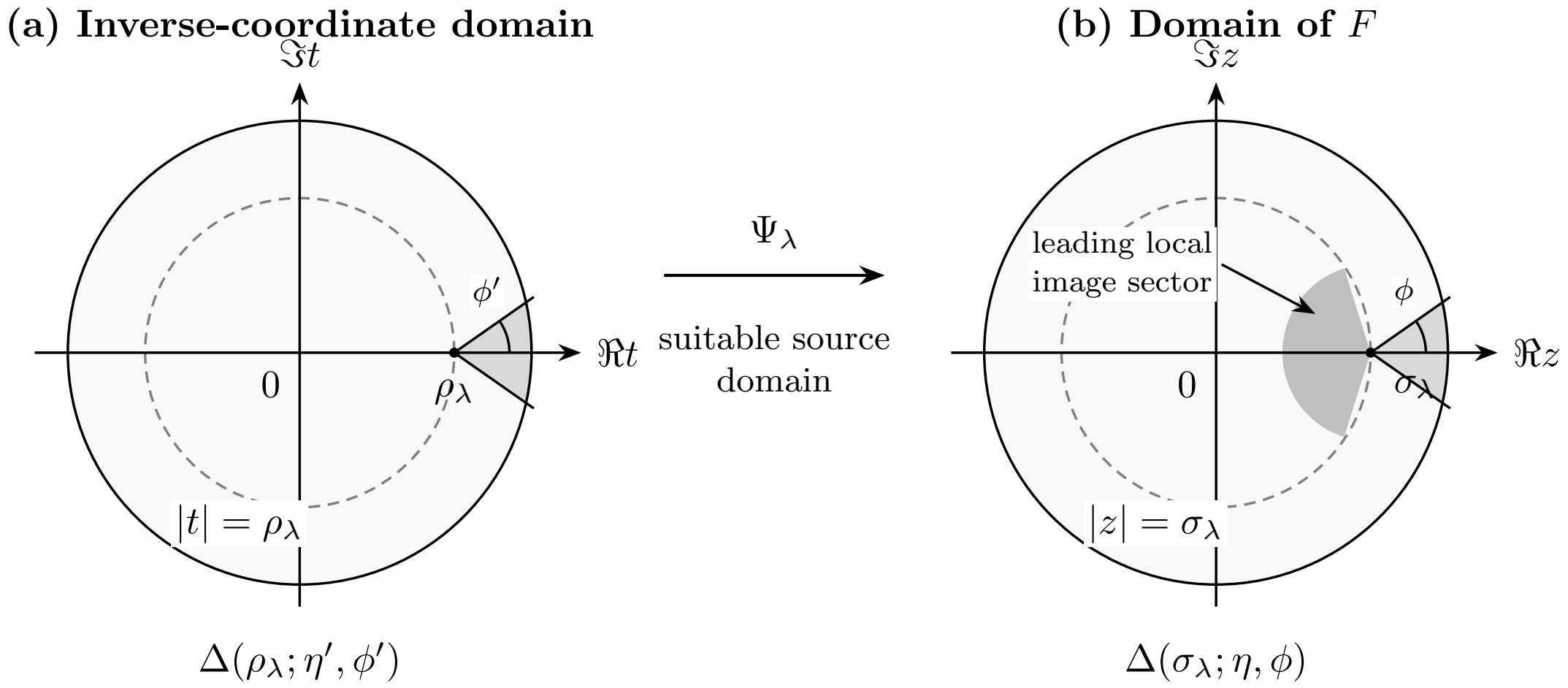
\begin{figure}[htbp]
\centering
\begin{tikzpicture}[>=Stealth,line width=.65pt,
 every node/.style={font=\small}]
\begin{scope}
\node at (0,2.95) {\textbf{(a) Inverse-coordinate domain}};
\fill[gray!5] (0,0) circle (2.1);
\begin{scope}
\clip (0,0) circle (2.1);
\fill[gray!30] (1.4,0)--(6,3.221)--(6,-3.221)--cycle;
\end{scope}
\draw (0,0) circle (2.1);
\draw[densely dashed,gray] (0,0) circle (1.4);
\draw[->] (-2.4,0)--(2.55,0) node[right] {$\Re t$};
\draw[->] (0,-2.3)--(0,2.45) node[above] {$\Im t$};
\draw (1.4,0)--++(35:.88);
\draw (1.4,0)--++(-35:.88);
\draw (1.9,0) arc[start angle=0,end angle=35,radius=.5];
\node[font=\scriptsize,fill=white,inner sep=.5pt] at (1.7,.55) {$\phi'$};
\fill (1.4,0) circle (1.3pt);
\node[below=3pt] at (1.4,0) {$\rho_\lambda$};
\node[below left=1pt] at (0,0) {$0$};
\node[fill=white,inner sep=1pt] at (-.55,-1.55) {$|t|=\rho_\lambda$};
\node at (0,-2.8) {$\Delta(\rho_\lambda;\eta',\phi')$};
\end{scope}
\draw[->,thick] (3.3,.7)--(5.3,.7);
\node[above=3pt] at (4.3,.7) {$\Psi_\lambda$};
\node[align=center,font=\footnotesize] at (4.3,-.05) {suitable source\\domain};
\begin{scope}[xshift=8.3cm]
\node at (0,2.95) {\textbf{(b) Domain of $F$}};
\fill[gray!5] (0,0) circle (2.1);
\begin{scope}
\clip (0,0) circle (2.1);
\fill[gray!30] (1.4,0)--(6,3.221)--(6,-3.221)--cycle;
\end{scope}
\fill[gray!50] (1.4,0)--++(107.5:.8)
 arc[start angle=107.5,end angle=252.5,radius=.8]--cycle;
\draw (0,0) circle (2.1);
\draw[densely dashed,gray] (0,0) circle (1.4);
\draw[->] (-2.4,0)--(2.55,0) node[right] {$\Re z$};
\draw[->] (0,-2.3)--(0,2.45) node[above] {$\Im z$};
\draw (1.4,0)--++(35:.88);
\draw (1.4,0)--++(-35:.88);
\draw (1.9,0) arc[start angle=0,end angle=35,radius=.5];
\node[font=\scriptsize,fill=white,inner sep=.5pt] at (1.7,.55) {$\phi$};
\fill (1.4,0) circle (1.3pt);
\node[below right=2pt] at (1.4,0) {$\sigma_\lambda$};
\node[below left=1pt] at (0,0) {$0$};
\node[fill=white,inner sep=1pt] at (-.55,-1.55) {$|z|=\sigma_\lambda$};
\node[align=center,font=\scriptsize,fill=white,inner sep=1pt] at (-.85,.8) {leading local\\image sector};
\draw[->] (.05,.8)--(.9,.35);
\node at (0,-2.8) {$\Delta(\sigma_\lambda;\eta,\phi)$};
\end{scope}
\node at (4.15,-3.55)
 {$\Psi_\lambda(t)-\sigma_\lambda
 \sim-a_\lambda(1-t/\rho_\lambda)^{1/2}$};
\end{tikzpicture}
\caption{Dented domains for analytic composition. The light grey cones
are excluded; the dashed circles mark the radii $\rho_\lambda$ and
$\sigma_\lambda$. The outer circles have radii $\rho_\lambda+\eta'$
and $\sigma_\lambda+\eta$, respectively. After choosing a suitable
source domain, $\Psi_\lambda$ takes it into the domain of $F$.
The darker sector illustrates only the leading local square-root
model near $\sigma_\lambda$, directed away from the excluded cone;
it is not the exact image of the whole source domain. The drawing
is schematic and uses independent scales in the two planes.}
\label{fig:delta-composition}
\end{figure}

The domain geometry and the leading local behavior of the inverse
are illustrated in Figure~\ref{fig:delta-composition}.

For $\lambda=1$ the constants are
\[
\sigma_1=1,\qquad\rho_1=e-2,\qquad
a_1=\sqrt{2(e-2)/e},
\]
as in \cite[Section~3.4]{CottiPaperI}.

\subsection{Minimal analytic hypotheses}
\label{sec:minimal-analytic-hypothesis}

Recall the exact formal decompositions from Sections~\ref{sec:motivicgraphs} and~\ref{sec:g1}:
\begin{equation}
\begin{aligned}
A_1^\mu&=\tfrac12\log P_\lambda+S_1^\mu,\\
A_g^\mu&=\gamma_g e^{2(g-1)(\lambda-2)z}P_\lambda^{3g-3}
+\mathcal R_g^\mu\quad(g\ge2).
\end{aligned}
\label{eq:motivicgraphexpansion-analytic}
\end{equation}
We impose the following conditions on the corresponding germs.

\smallskip\noindent\textbf{Hypothesis $(\mathrm H_1)$.}
The series $S_1^\mu$ converges near zero, continues analytically to
some $\Delta(\sigma_\lambda;\eta,\phi)$, and extends holomorphically
across $\sigma_\lambda$.

\smallskip\noindent\textbf{Hypothesis $(\mathrm H_g)$, $g\ge2$.}
The series $\mathcal R_g^\mu$ converges near zero and continues
analytically to some $\Delta(\sigma_\lambda;\eta,\phi)$, where
\begin{equation}
\mathcal R_g^\mu(z)=O\bigl(|1-z/\sigma_\lambda|^{-(3g-4)}\bigr)
\quad\text{as }z\to\sigma_\lambda.
\label{eq:Hg}
\end{equation}
The constants and domains may depend on the measure and genus.
These hypotheses concern only the relevant total remainder; they do
not require separate bounds on every vertex or residual family.
The tameness conditions below will be sufficient to establish them.

\subsection{Coefficient asymptotics in fixed genus}
\label{sec:motivic-fixed-genus}

All limits below keep the genus and realization fixed. We use
singularity transfer in $\Delta$-domains, including logarithmic
transfer and transfer of local Puiseux expansions
\cite[Chapter~VI]{FS}. The conversion from an ordinary series
coefficient to a moduli-space measure includes the factor $n!$.

\smallskip\noindent\emph{Genus zero.}
\begin{theorem}
\label{thm:motivic-genus-zero-asymptotics}
For a scalar motivic realization with $\lambda>0$, no additional
analytic hypothesis is required in genus zero, and
\begin{equation}
\mu(\overline M_{0,n})\sim
\frac{\rho_\lambda e^{\sigma_\lambda}a_\lambda}{2\sqrt\pi}
\,n!\,\rho_\lambda^{-n}n^{-5/2}.
\label{eq:motivic-g0-asymptotic}
\end{equation}
More precisely, there are coefficients $\kappa_j(\lambda)$ such that,
for every integer $N\ge1$,
\[
\mu(\overline M_{0,n})=
\frac{\rho_\lambda e^{\sigma_\lambda}a_\lambda}{2\sqrt\pi}
\,n!\,\rho_\lambda^{-n}n^{-5/2}
\left(1+\sum_{j=1}^{N-1}\frac{\kappa_j(\lambda)}{n^j}
+O(n^{-N})\right).
\]
These coefficients depend only on $\lambda$ and are obtained from
the local inverse of $\Phi_\lambda$ as described below.
\end{theorem}
\begin{proof}
Section~\ref{sec:genuszero-coordinate} gives
$1+t+(G_0^\mu)'(t)=e^{\Psi_\lambda(t)}$.
Put $X=(1-t/\rho_\lambda)^{1/2}$. By
Proposition~\ref{prop:universal-singularity-pair}, the inverse has a
convergent expansion
\[
w(X):=\Psi_\lambda\bigl(\rho_\lambda(1-X^2)\bigr)
=\sigma_\lambda+\sum_{j\ge1}b_j(\lambda)X^j,
\qquad b_1(\lambda)=-a_\lambda.
\]
Its coefficients are determined successively by
\[
\sum_{r\ge2}\frac{\Phi_\lambda^{(r)}(\sigma_\lambda)}{r!}
\left(\sum_{j\ge1}b_j(\lambda)X^j\right)^r
=-\rho_\lambda X^2.
\]
The choice $b_1=-a_\lambda$ fixes the branch, and for $j\ge2$
the coefficient of $X^{j+1}$ determines $b_j$ linearly, since
$\Phi_\lambda''(\sigma_\lambda)b_1\ne0$.

Integrating $(G_0^\mu)'=e^{\Psi_\lambda}-1-t$, with
$dt=-2\rho_\lambda X\,dX$, gives
\[
G_0^\mu(t)=H_{0,\lambda}(t)
+\sum_{j\ge0}c_j(\lambda)(1-t/\rho_\lambda)^{j+3/2},
\qquad
c_j(\lambda)=-\frac{2\rho_\lambda}{2j+3}
[X^{2j+1}]e^{w(X)},
\]
where $H_{0,\lambda}$ is holomorphic in a neighborhood of
$t=\rho_\lambda$. Both parts come from convergent local series.
In particular,
\[
c_0(\lambda)=\frac23\rho_\lambda e^{\sigma_\lambda}a_\lambda.
\]
The nonzero first nonanalytic term and continuation across all other
points of the boundary identify the unique dominant singularity.
For each fixed $N\ge1$, termwise singularity transfer
\cite[Chapter~VI]{FS} yields
\[
\mu(\overline M_{0,n})=
\rho_\lambda^{-n}\sum_{j=0}^{N-1}
c_j(\lambda)\frac{\Gamma(n-j-3/2)}{\Gamma(-j-3/2)}
+O\bigl(n!\,\rho_\lambda^{-n}n^{-N-5/2}\bigr).
\]
Expanding the gamma quotients in powers of $n^{-1}$ determines the
coefficients $\kappa_j(\lambda)$. The leading term follows from
$\Gamma(-3/2)=4\sqrt\pi/3$.
\end{proof}

For $\lambda=1$, the inverse can be written
\[
\Psi_1(t)=2+W_{-1}\bigl(-(t+2)e^{-2}\bigr).
\]
Thus the local inversion above recovers the expansion obtained from
the branch $W_{-1}$ of the Lambert function in
\cite[Theorem~2.18; see also Section~3.8]{Cotti}.
For general $\lambda>0$,
the same procedure gives every order directly from $\Phi_\lambda$.

\Needspace{6\baselineskip}
\smallskip\noindent\emph{Genus one.}
\begin{theorem}
\label{thm:motivic-genus-one-asymptotics}
If $\lambda>0$ and $(\mathrm H_1)$ holds, then
\begin{align}
\mu(\overline M_{1,(0^+)^n})&\sim\tfrac12 n!\,
\sigma_\lambda^{-n}n^{-1},
\label{eq:g1-alllight-asymptotic}\\
\mu(\overline M_{1,n})&\sim\tfrac14 n!\,
\rho_\lambda^{-n}n^{-1}.
\label{eq:g1-ordinary-asymptotic}
\end{align}
\end{theorem}
\begin{proof}
The simple propagator pole and $(\mathrm H_1)$ give, locally,
\[
A_1^\mu(z)=-\tfrac12\log(1-z/\sigma_\lambda)
+C^A_{1,\mu}+O(1-z/\sigma_\lambda),
\]
with a remainder holomorphic at $\sigma_\lambda$. Logarithmic
transfer gives the first equivalent. The heavy/light identity gives
$G_1^\mu=A_1^\mu\circ\Psi_\lambda$, and
\[
1-\Psi_\lambda(t)/\sigma_\lambda
=\frac{a_\lambda}{\sigma_\lambda}(1-t/\rho_\lambda)^{1/2}
\bigl(1+O((1-t/\rho_\lambda)^{1/2})\bigr).
\]
Consequently
$G_1^\mu(t)=-\tfrac14\log(1-t/\rho_\lambda)+C^G_{1,\mu}
+O((1-t/\rho_\lambda)^{1/2})$.
Lemma~\ref{lem:analytic-composition} gives the continuation needed
for logarithmic transfer and the second equivalent.
\end{proof}

\smallskip\noindent\emph{Higher genus.}
Fix $g\ge2$ and set
\[
m=3g-3,\qquad \alpha_{g,\lambda}=2(g-1)(\lambda-2).
\]
\begin{theorem}[All-light asymptotics]
\label{thm:motivic-alllight-asymptotics}
If $\lambda>0$ and $(\mathrm H_g)$ holds, then
\begin{equation}
\begin{aligned}
A_g^\mu(z)&=K^A_{g,\lambda}(1-z/\sigma_\lambda)^{-m}
+O(|1-z/\sigma_\lambda|^{-m+1}),\\
K^A_{g,\lambda}&=\frac{\gamma_g e^{\alpha_{g,\lambda}\sigma_\lambda}}
{(\lambda\sigma_\lambda)^m}.
\end{aligned}
\label{eq:motivic-alllight-local}
\end{equation}
In particular,
\begin{equation}
\mu(\overline M_{g,(0^+)^n})\sim
\frac{\gamma_g e^{\alpha_{g,\lambda}\sigma_\lambda}}
{\Gamma(m)(\lambda\sigma_\lambda)^m}
\,n!\,\sigma_\lambda^{-n}n^{m-1}.
\label{eq:motivic-alllight-asymptotic}
\end{equation}
\end{theorem}
\begin{proof}
Insert \eqref{eq:Plambda-local} into the trivalent term of the
stabilization-graph expansion. Its leading coefficient is $K^A_{g,\lambda}$,
and its next order is at most $m-1$. Hypothesis $(\mathrm H_g)$
controls the rest at that order in a $\Delta$-domain. Transfer of
$(1-z/\sigma_\lambda)^{-m}$ proves the result.
\end{proof}

\begin{theorem}[Ordinary asymptotics]
\label{thm:motivic-ordinary-asymptotics}
Under the same assumptions,
\begin{equation}
\begin{aligned}
G_g^\mu(t)&=K^G_{g,\lambda}(1-t/\rho_\lambda)^{-m/2}
+O(|1-t/\rho_\lambda|^{-(m-1)/2}),\\
K^G_{g,\lambda}&=\frac{\gamma_g e^{\alpha_{g,\lambda}\sigma_\lambda}}
{(\lambda a_\lambda)^m}.
\end{aligned}
\label{eq:motivic-ordinary-local}
\end{equation}
Hence
\begin{equation}
\mu(\overline M_{g,n})\sim
\frac{\gamma_g e^{\alpha_{g,\lambda}\sigma_\lambda}}
{\Gamma(m/2)(\lambda a_\lambda)^m}
\,n!\,\rho_\lambda^{-n}n^{m/2-1}.
\label{eq:motivic-ordinary-asymptotic}
\end{equation}
\end{theorem}
\begin{proof}
Substitute \eqref{eq:Psi-lambda-local} into the preceding local
expansion. The factor $1-\Psi_\lambda(t)/\sigma_\lambda$ begins
with $(a_\lambda/\sigma_\lambda)(1-t/\rho_\lambda)^{1/2}$.
Thus the leading constant is
$K^A_{g,\lambda}(\sigma_\lambda/a_\lambda)^m=K^G_{g,\lambda}$,
and the error has the stated order. Analytic continuation follows
from Lemma~\ref{lem:analytic-composition}; singularity transfer
gives the coefficient equivalent.
\end{proof}

\begin{corollary}[Heavy/light ratios]
\label{cor:motivic-ratio}
Under the corresponding hypotheses, the all-light to ordinary ratio is
\begin{equation}
\frac{\mu(\overline M_{g,(0^+)^n})}{\mu(\overline M_{g,n})}
\sim\begin{cases}
2(\rho_\lambda/\sigma_\lambda)^n,&g=1,\\[1mm]
\displaystyle
\frac{\Gamma(m/2)}{\Gamma(m)}
\left(\frac{a_\lambda}{\sigma_\lambda}\right)^m
\left(\frac{\rho_\lambda}{\sigma_\lambda}\right)^n n^{m/2},&g\ge2,
\end{cases}
\label{eq:motivic-ratio}
\end{equation}
In every case the ratio tends to zero.
\end{corollary}
\begin{proof}
Divide the equivalents. The common trivalent factor cancels for
$g\ge2$, and $0<\rho_\lambda/\sigma_\lambda<1$ gives the limit.
\end{proof}

The leading terms depend on the realization only through $\lambda$,
subject to the stated hypotheses. At $\lambda=1$ these formulas
recover the Euler equivalents in
\cite[Sections~2.5, 3.4 and~6]{CottiPaperI}. The complete series,
however, need not be determined by their Lefschetz value.

\subsection{Vertex and residual quotient tameness}
\label{sec:vertex-tameness}

\begin{definition}
\label{def:vertex-tame}
A scalar realization is \emph{vertex-tame} if for each stable pair
$(h,d)$ there are constants $B,C>0$ such that
\[
|\mu(\mathcal V_{h,d,n})|\le BC^n\qquad(n\ge0).
\]
It is \emph{residual-quotient tame} if for every residual family
from Section~\ref{sec:motivicgraphs} there are constants $B,C>0$ such that
\[
|\mu(Y_{\Gamma,\tau,n})|+
|\mu(Y_{\Gamma,\tau,n}/G_{\Gamma,\tau})|\le BC^n.
\]
It is \emph{graph-tame} if it has both properties.
\label{def:graph-tame}
\end{definition}
The constants may depend on the realization and the indicated geometric
type, but not on $n$. To see the analytic consequence of these
bounds, let $F(z)=\sum_{n\ge0}a_nz^n/n!$ be any one of the vertex,
residual, or quotient series above. An estimate $|a_n|\le BC^n$ gives
\[
\sum_{n\ge0}\left|a_n\frac{z^n}{n!}\right|
\le B\sum_{n\ge0}\frac{(C|z|)^n}{n!}=Be^{C|z|}
\qquad(z\in\C).
\]
Thus the series converges absolutely everywhere and defines an entire
function of exponential type. In particular, graph tameness makes
every $V_{h,d}^\mu$, $Y_{\Gamma,\tau}^\mu$, and
$Q_{\Gamma,\tau}^\mu$ occurring in the graph expansion entire.

\begin{unnumberedremark}[Vertex tameness and residual varieties]
Vertex tameness already implies the estimate
$|\mu(Y_{\Gamma,\tau,n})|\le BC^n$ for the residual varieties before
quotienting. Indeed, fix $\tau=(p,J)$ and use the construction of
\eqref{eq:residual-varieties}. Its number of component positions is
\[
N=|\Vertices(\Gamma)|+|\{e:a_p(e)=1\}|+|L_p\setminus J|
\ge |\Vertices(\Gamma)|\ge1.
\]
Each admissible allocation $[n]=B_1\sqcup\cdots\sqcup B_N$ contributes
a product $\prod_{i=1}^N F_{i,B_i}$ of vertex or one-component chain
factors, with a point allowed at an empty optional position.
Vertex tameness and
$\mu(\mathcal C_B)=(\lambda-1)^{|B|-1}$ for $B\ne\varnothing$
allow constants $b,C\ge1$, independent of $n$ and the allocation,
such that $|\mu(F_{i,B_i})|\le bC^{|B_i|}$ for every allowed factor.
There are at most $N^n$ allocations. Hence
\[
|\mu(Y_{\Gamma,\tau,n})|
\le\sum_{\text{admissible allocations}}
\prod_{i=1}^N|\mu(F_{i,B_i})|
\le N^n b^N C^n=b^N(NC)^n.
\]
Thus, under vertex tameness, the additional requirement for graph
tameness is the bound on the actual quotient:
\[
|\mu(Y_{\Gamma,\tau,n}/G_{\Gamma,\tau})|\le B'C'^{\,n}.
\]
\qrem
\end{unnumberedremark}

\begin{proposition}[Tameness implies the analytic hypotheses]
\label{prop:vertex-tame-Hg}
For $\lambda>0$, vertex tameness implies $(\mathrm H_1)$ and graph
tameness implies $(\mathrm H_g)$ for every $g\ge2$.
\end{proposition}
\begin{proof}
The identity $(I^\mu)'=e^zV_{1,1}^\mu$ from Section~\ref{sec:g1} makes $I^\mu$
entire under vertex tameness. Since $U_\lambda$ is entire, so is
$S_1^\mu$, proving $(\mathrm H_1)$.

For $g\ge2$, the explicit global remainder of Section~\ref{sec:motivicgraphs} is a finite
sum of terms $F(z)P_\lambda(z)^r$, with $0\le r\le3g-4$.
Here each $F$ is a product of vertex series or a constant multiple
of $Q_{\Gamma,\tau}^\mu-Y_{\Gamma,\tau}^\mu/|G_{\Gamma,\tau}|$.
Graph tameness makes every $F$ entire. Choose $\eta>0$ so that the
disk $|z|<\sigma_\lambda+\eta$ contains no other propagator pole.
The finite sum is meromorphic in that disk with pole order at most
$3g-4$, proving $(\mathrm H_g)$.
\end{proof}
The tameness conditions are sufficient, not necessary. At any fixed
genus it suffices to control the finitely many vertex pairs and
residual families that actually occur.

\subsection{Cohomological bounds and Hodge realizations}
\label{sec:hodge-vertex-tame}

All cohomology groups in this subsection are taken in the ordinary
analytic topology, with rational coefficients and with compact supports
where indicated. Write $b_c(Y)=\sum_k\dim_\Q H_c^k(Y;\Q)$.
The following bounds will control the residual geometry without
requiring an explicit computation of the quotient Hodge polynomials.

\begin{lemma}[Finite quotients]
\label{lem:finite-quotient-bound}
For a finite group $G$ acting algebraically on a complex quasi-projective
variety $Y$,
\begin{equation}
H_c^\bullet(Y/G;\Q)\cong H_c^\bullet(Y;\Q)^G
\label{eq:finite-quotient-invariants}
\end{equation}
as mixed Hodge structures. In particular, $b_c(Y/G)\le b_c(Y)$.
\end{lemma}
\begin{proof}
The quotient map $q$ is finite and proper. Proper base change gives
$R^i q_*\Q_Y=0$ for $i>0$ and $(q_*\Q_Y)^G=\Q_{Y/G}$, since
each fiber is a single finite orbit. For any rational representation
$W$ of $G$, the averaging operator
\[
P_G:=\frac1{|G|}\sum_{g\in G}g:W\longrightarrow W
\qquad\text{satisfies}\qquad
P_G^2=P_G,\quad\operatorname{im}P_G=W^G.
\]
Thus taking invariants is exact: if an equivariant map is surjective,
applying $P_G$ to any lift of an invariant vector gives an invariant
lift. The same operator acts on sheaves with a $G$-action and commutes
with their equivariant differentials. Applying compactly supported
cohomology therefore gives the stated isomorphism. It is induced by
proper algebraic pullback and is compatible with mixed Hodge structures;
the invariant subspace has no larger dimension. For the smooth case see also
\cite[Section~4.2, Corollary~4.4]{FlorentinoSilva}.
\end{proof}

\begin{lemma}[Quotient presentation of the vertex spaces]
\label{lem:vertex-quotient}
For each stable pair $(h,d)$ there is a finite group $K$ such that, for
every $n\ge0$, there exists a smooth quasi-projective variety
$\widetilde{\mathcal V}_{h,d,n}$ equipped with a $K$-action and satisfying
\[
\mathcal V_{h,d,n}\cong\widetilde{\mathcal V}_{h,d,n}/K.
\]
\end{lemma}
\begin{proof}
If $h=0$, stability gives $d\ge3$. Set $U=M_{0,d}$ and $K=\{1\}$;
then $U$ is a fine moduli variety with its universal smooth $d$-pointed
rational curve. If $h>0$, fix an integer $m\ge3$ and add a full
symplectic level-$m$ structure. For $h\ge2$, let
\[
U=M_h[m]\times_{\mathcal M_h}\mathcal M_{h,d},
\]
where $M_h[m]$ is the fine moduli variety of smooth genus-$h$ curves
with this level structure. The existence of the universal family
and the level cover is recalled in
\cite[Remark~2.4 and Chapter~III, pp.~249--251]{Popp}.
For $h=1$, stability gives $d\ge1$; using the first heavy marking
as origin, take instead
\[
U=M_{1,1}[m]\times_{\mathcal M_{1,1}}\mathcal M_{1,d},
\]
with $M_{1,1}[m]$ the fine modular curve of full level $m$;
see \cite{DS}. In either case the remaining heavy markings are
pairwise distinct points of the universal curve. The resulting $U$
is smooth and quasi-projective, and forgetting the level gives a
finite \'{e}tale $K$-torsor $U\to\mathcal M_{h,d}$, where $K$ is the
group of changes of symplectic level structure.

Over $U$, let $\mathcal C^\circ\to U$ be the universal smooth curve
with its $d$ heavy sections removed, and set
\[
f_n:\widetilde{\mathcal V}_{h,d,n}
:=(\mathcal C^\circ)^{\times_U n}\longrightarrow U.
\]
The $n$ light markings may coincide, so no diagonals are removed
from this fiber product. The $K$-action lifts to the universal family
and fixes all marking labels. Forgetting the level structure gives
\[
[\widetilde{\mathcal V}_{h,d,n}/K]
\cong\mathscr V_{h,[d],[n]},\qquad
\mathcal V_{h,d,n}\cong\widetilde{\mathcal V}_{h,d,n}/K,
\]
with the coarse-quotient identification of
Section~\ref{sec:parameter-and-quotient}. This also covers $h=0$
with the trivial group.

\end{proof}

\begin{lemma}[Vertex Betti bounds]
\label{lem:vertex-betti-bound}
For every stable pair $(h,d)$ there are constants $B,C>0$ such that
\[
b_c(\mathcal V_{h,d,n})\le BC^n\qquad(n\ge0).
\]
\end{lemma}
\begin{proof}
Use the construction in the proof of Lemma~\ref{lem:vertex-quotient} and
write $f_n:\widetilde{\mathcal V}_{h,d,n}\to U$ for the projection.
For the cohomological estimate we now use the analytifications of
these varieties, retaining the same notation. Before removing the
sections, the universal curve $\pi:\mathcal C\to U$ is a smooth
proper submersion. C.\,Ehresmann's fibration theorem
\cite[Section~1]{Ehresmann} therefore gives local differentiable
trivializations. These can be chosen to preserve the $d$ disjoint
marked sections: choose horizontal lifts of base vector fields that
are tangent to each section, extend them by a partition of unity,
and use properness to integrate them over a sufficiently small base
neighborhood. Parallel transport then preserves each section.
Removing the sections and taking fiber powers consequently makes
$f_n$ a locally trivial fibration in the ordinary analytic topology,
with fiber $F_{h,d}^{\,n}$, where $F_{h,d}$ is a smooth genus-$h$
curve with $d$ punctures. Put
\[
b:=\sum_q\dim_\Q H^q(F_{h,d};\Q).
\]

The quasi-projective variety $U$ has the homotopy type of a finite
CW complex. Indeed, choose a projective compactification and a finite
semialgebraic triangulation compatible with its boundary; the
complement of the boundary subcomplex deformation retracts onto a
finite complex after barycentric subdivision. For the triangulation
theorem see \cite[Section~9.2]{BochnakCosteRoy}.
Choose such a CW model $T$ with $c_p$ cells in dimension $p$, and put
$c_U=\sum_p c_p$. For a local system $\mathcal L$ of rank $r$, its
pullback to this model has cellular cochains with one $r$-dimensional
summand for each cell. Consequently,
\[
\dim_\Q C^p(T;\mathcal L)=c_p r,\qquad
\sum_p\dim_\Q H^p(U;\mathcal L)
\le\sum_p\dim_\Q C^p(T;\mathcal L)=c_Ur,
\]
where $\mathcal L$ in the cellular complex denotes its pullback to $T$.

Since $f_n$ is locally trivial, the sheaves
$\mathcal L_n^q:=R^q(f_n)_*\Q$ are local systems with fibers
$H^q(F_{h,d}^{\,n};\Q)$. K\"unneth gives
\[
\sum_q\operatorname{rank}\mathcal L_n^q
=\sum_q\dim_\Q H^q(F_{h,d}^{\,n};\Q)=b^n.
\]
The Leray spectral sequence is
\[
E_2^{p,q}=H^p(U;\mathcal L_n^q)
\quad\Longrightarrow\quad
H^{p+q}(\widetilde{\mathcal V}_{h,d,n};\Q).
\]
Its limiting terms are the graded pieces of a filtration on the
cohomology of the total space, and each page is obtained by taking
cohomology of the preceding one. Hence
\begin{align*}
\sum_k\dim_\Q H^k(\widetilde{\mathcal V}_{h,d,n};\Q)
&=\sum_{p,q}\dim_\Q E_\infty^{p,q}
\le\sum_{p,q}\dim_\Q E_2^{p,q}\\
&\le c_U\sum_q\operatorname{rank}\mathcal L_n^q
=c_Ub^n.
\end{align*}
Finally, the total space is smooth of complex dimension
$\dim U+n$. Poincar\'e duality identifies the dimensions of
$H_c^k$ and $H^{2(\dim U+n)-k}$, so it preserves the total Betti
number. Applying Lemma~\ref{lem:finite-quotient-bound} to the
$K$-quotient therefore gives
\[
b_c(\mathcal V_{h,d,n})
\le b_c(\widetilde{\mathcal V}_{h,d,n})\le c_Ub^n.
\]
The constants $c_U$ and $b$ depend only on $(h,d)$, as required.
\end{proof}

\begin{proposition}[Hodge--Deligne realizations are graph-tame]
\label{prop:hodge-vertex-tame}
For every fixed $u,v\in\C$, the realization
$\mu_{u,v}(X)=E_c(X;u,v)$ is graph-tame. If $uv$ is a positive real
number, it satisfies $(\mathrm H_1)$ and $(\mathrm H_g)$ for all $g\ge2$.
\end{proposition}
\begin{proof}
For a variety of dimension at most $D$, the compactly supported Hodge
indices lie in $[0,D]$. With $M=\max\{1,|u|,|v|\}$ this gives
\begin{equation}
|E_c(Y;u,v)|\le M^{2D}b_c(Y).
\label{eq:scalar-Hodge-bound}
\end{equation}
Since $\dim\mathcal V_{h,d,n}=3h-3+d+n$, the vertex Betti bound
proves vertex tameness for every fixed pair of Hodge variables.

Fix $\Gamma,\tau$. Recall from \eqref{eq:residual-varieties} that the
residual varieties are disjoint unions of products of vertex spaces
and empty or one-component chain spaces. We may take
\[
d_\Gamma:=\sum_{v\in\Vertices(\Gamma)}(3h_v-3+d_v).
\]
This is the total dimension of the vertex factors before light
markings are added; it depends only on $\Gamma$, not on $\tau$ or $n$.
Since a vertex factor with $s$ light
labels has dimension $3h_v-3+d_v+s$, and a one-component chain with
$s\ge1$ labels has dimension $s-1\le s$, every summand has dimension
at most $d_\Gamma+n$.

Vertex
factors have exponential compact-Betti bounds; a one-component
chain with $s\ge1$ labels has $b_c(\GG^{s-1})=2^{s-1}$, and an
empty optional position contributes a point. K\"unneth and the
exponential bound on label assignments give
\[
b_c(Y_{\Gamma,\tau,n})\le A B^n,\qquad
\dim Y_{\Gamma,\tau,n}\le d_\Gamma+n
\]
for constants $A,B>0$ independent of $n$. The finite quotient has no larger
total compact Betti number and has the same dimension bound.
Applying \eqref{eq:scalar-Hodge-bound} to both gives
\[
|E_c(Y_{\Gamma,\tau,n};u,v)|+
|E_c(Y_{\Gamma,\tau,n}/G_{\Gamma,\tau};u,v)|
\le2A M^{2d_\Gamma}(M^2B)^n.
\]
This proves residual quotient tameness. The final assertion follows
from Proposition~\ref{prop:vertex-tame-Hg} with $\lambda=uv>0$.
\end{proof}

The quotient estimates use invariant cohomology, not the Hodge
polynomials of fixed loci. For each fixed pair $(u,v)$ with positive
product, they verify the hypotheses of all the preceding positive
Lefschetz asymptotic theorems. Their consequences for fixed-product
universality and concentration of Hodge mass are derived in Section~\ref{sec:asymptotic-hodge-tate}.

\section{Hodge universality and asymptotic diagonal concentration}
\label{sec:asymptotic-hodge-tate}

The Hodge-tameness result of Section~\ref{sec:hodge-vertex-tame} verifies the analytic hypotheses
for every fixed scalar Hodge--Deligne evaluation with positive product.
The leading asymptotics therefore depend on the two Hodge variables
only through that product. We now turn this universality into a
statement about the Hodge numbers themselves: the mass outside the
diagonal is negligible at fixed genus. At product one this also
identifies the leading growth of the Euler characteristic and the
total Betti number.

\subsection{Universality at fixed Hodge product}
\label{sec:fixed-product-universality}

Fix $q>0$. We use the constants of Section~\ref{sec:universal-analytic-geometry} with $\lambda=q$:
\[
\sigma_q=\frac{\log q}{q-1}\quad(q\ne1),\quad\sigma_1=1,
\qquad \rho_q=e^{\sigma_q}-1-q^{-1},\qquad
a_q=\sqrt{\frac{2\rho_q}{q e^{\sigma_q}}}.
\]
In particular, $0<\rho_q<\sigma_q$.

\begin{theorem}[Fixed-product Hodge universality]
\label{thm:fixed-q-hodge-universality}
Fix $q>0$ and $g\ge1$. For every fixed pair $u,v\in\C$ with $uv=q$,
the leading asymptotics are independent of the individual values of
$u$ and $v$. In genus one,
\begin{align}
E_c(\overline M_{1,(0^+)^n};u,v)
&\sim\tfrac12 n!\,\sigma_q^{-n}n^{-1},
\label{eq:hodge-g1-alllight-fixed-q}\\
E_c(\overline M_{1,n};u,v)
&\sim\tfrac14 n!\,\rho_q^{-n}n^{-1}.
\label{eq:hodge-g1-ordinary-fixed-q}
\end{align}
For $g\ge2$, set $m=3g-3$. Then
\begin{align}
E_c(\overline M_{g,(0^+)^n};u,v)
&\sim\mathcal A_{g,q}\,n!\,\sigma_q^{-n}n^{m-1},
\label{eq:hodge-alllight-fixed-q}\\
E_c(\overline M_{g,n};u,v)
&\sim\mathcal G_{g,q}\,n!\,\rho_q^{-n}n^{m/2-1},
\label{eq:hodge-ordinary-fixed-q}
\end{align}
where the positive leading constants are
\begin{equation}
\mathcal A_{g,q}=
\frac{\gamma_g e^{2(g-1)(q-2)\sigma_q}}{\Gamma(m)(q\sigma_q)^m},
\qquad
\mathcal G_{g,q}=
\frac{\gamma_g e^{2(g-1)(q-2)\sigma_q}}{\Gamma(m/2)(q a_q)^m}.
\label{eq:hodge-leading-constants}
\end{equation}
The limits keep $g,q,u,v$ fixed. No uniformity over all pairs with
product $q$, or over increasing genus, is asserted.
\end{theorem}
\begin{proof}
The measure $\mu_{u,v}(X)=E_c(X;u,v)$ is graph-tame by Proposition~\ref{prop:hodge-vertex-tame}.
Its Lefschetz value is $q>0$, so it satisfies the analytic hypotheses
for genus one and for each higher genus. Applying the corresponding
scalar asymptotic theorems gives the formulas above. All factors in
their leading terms depend only on $g$ and $q$.
\end{proof}

The full Hodge polynomials may still depend separately on $u$ and $v$;
the theorem concerns their leading fixed-genus growth. This distinction
is visible already in the cuspidal terms of the genus-one formula.

\subsection{Concentration of Hodge mass on the diagonal}
\label{sec:diagonal-concentration}

For the remainder of the section let $X_n$ be either of the coarse
families $\overline M_{g,(0^+)^n}$ or $\overline M_{g,n}$, with
$g\ge1$ fixed. As recalled in Section~\ref{sec:g1}, its rational cohomology
agrees with that of the smooth proper Hassett stack and is pure.
Write $h^{p,r}(X_n)$ for the dimension of the $(p,r)$-summand of
$H^{p+r}(X_n;\C)$.
Thus the coefficients are nonnegative and
\begin{equation}
E_c(X_n;-x,-y)=\sum_{p,r}h^{p,r}(X_n)x^p y^r
\qquad(x,y>0).
\label{eq:positive-Hodge-evaluation}
\end{equation}
Define the total and off-diagonal $q$-weighted masses by
\begin{equation}
T_q(X_n)=\sum_{p,r}h^{p,r}(X_n)q^{(p+r)/2},\qquad
N_q(X_n)=\sum_{p\ne r}h^{p,r}(X_n)q^{(p+r)/2}.
\label{eq:hodge-masses}
\end{equation}
The total mass is positive, since constant cohomology classes
contribute. For a fixed real $\theta$ set
\begin{equation}
H_{n,q}(\theta)=
E_c(X_n;-\sqrt q\,e^\theta,-\sqrt q\,e^{-\theta})
=\sum_{p,r}h^{p,r}(X_n)q^{(p+r)/2}e^{\theta(p-r)}.
\label{eq:hodge-moment}
\end{equation}
The two Hodge arguments have product $q$, and $H_{n,q}(0)=T_q(X_n)$.

\begin{theorem}[Asymptotic diagonal concentration]
\label{thm:asymptotic-hodge-tate}
For every fixed $q>0$ and $g\ge1$, both families satisfy
\begin{equation}
\frac{N_q(X_n)}{T_q(X_n)}\longrightarrow0.
\label{eq:asymptotic-hodge-tate}
\end{equation}
Equivalently,
\begin{equation}
\sum_p h^{p,p}(X_n)q^p\sim T_q(X_n).
\label{eq:diagonal-hodge-dominance}
\end{equation}
\end{theorem}
\begin{proof}
Choose one fixed $\theta_0>0$. Fixed-product universality gives
$H_{n,q}(\theta_0)\sim H_{n,q}(0)$ and
$H_{n,q}(-\theta_0)\sim H_{n,q}(0)$. Averaging the two evaluations
and subtracting the value at zero yields the nonnegative expression
\begin{align*}
D_n&:=\frac{H_{n,q}(\theta_0)+H_{n,q}(-\theta_0)}2-H_{n,q}(0)\\
&=\sum_{p,r}h^{p,r}(X_n)q^{(p+r)/2}
\bigl(\cosh(\theta_0(p-r))-1\bigr)
=o(T_q(X_n)).
\end{align*}
For $p\ne r$, the integer $p-r$ has absolute value at least one,
so the summand's last factor is at least $\cosh\theta_0-1>0$.
Consequently
\[
0\le(\cosh\theta_0-1)N_q(X_n)\le D_n=o(T_q(X_n)).
\]
Dividing by $T_q(X_n)$ proves the limit; subtracting $N_q$ from
$T_q$ gives the equivalent formulation.
Only the three fixed evaluations $0,\theta_0,-\theta_0$ are used,
so no uniformity in $\theta$ is needed.
\end{proof}

\begin{remark}
The description as asymptotic Hodge--Tate concentration refers to
the proportion of the total weighted Hodge mass on the diagonal.
It does not assert that any individual space eventually has only
diagonal Hodge components, or that each off-diagonal Hodge number
vanishes as $n$ grows.
\qrem
\end{remark}

\subsection{Euler characteristic and total Betti number}
\label{sec:euler-betti-consequence}

Let $b(X_n)=\sum_k b_k(X_n)$ and
$b_{\mathrm{odd}}(X_n)=\sum_{k\text{ odd}}b_k(X_n)$.
Purity identifies $b(X_n)=T_1(X_n)$. Every odd-degree Hodge
component is off-diagonal, since $p+r$ odd implies $p\ne r$.

\begin{corollary}
\label{cor:euler-total-betti}
For either family at fixed genus $g\ge1$,
\begin{equation}
\frac{b_{\mathrm{odd}}(X_n)}{b(X_n)}\longrightarrow0,
\qquad \chi(X_n)\sim b(X_n).
\label{eq:euler-total-betti-asymptotic}
\end{equation}
\end{corollary}
\begin{proof}
The exact identity and bounds are
\[
0\le b(X_n)-\chi(X_n)=2b_{\mathrm{odd}}(X_n)
\le2N_1(X_n).
\]
Divide by $b(X_n)=T_1(X_n)$ and apply diagonal concentration at $q=1$.
\end{proof}

Thus the Euler equivalents in \cite[Sections~3.4 and~6]{CottiPaperI}
also give the leading growth of the total Betti numbers in both chambers.
The odd cohomology can be nonzero; its total dimension is negligible
relative to the full cohomology.

\subsection{Hodge--Deligne heavy/light comparison}
\label{sec:hodge-chamber-comparison}

\begin{corollary}
\label{cor:hodge-heavy-light-ratio}
Fix $q>0$, $g\ge1$, and $u,v\in\C$ with $uv=q$. Then
\begin{equation}
\frac{E_c(\overline M_{g,(0^+)^n};u,v)}
{E_c(\overline M_{g,n};u,v)}
\sim\begin{cases}
2(\rho_q/\sigma_q)^n,&g=1,\\[1mm]
\displaystyle
\frac{\Gamma((3g-3)/2)}{\Gamma(3g-3)}
\left(\frac{a_q}{\sigma_q}\right)^{3g-3}
\left(\frac{\rho_q}{\sigma_q}\right)^n n^{(3g-3)/2},&g\ge2,
\end{cases}
\label{eq:hodge-heavy-light-ratio}
\end{equation}
In particular the ratio tends to zero.
\end{corollary}
\begin{proof}
Divide the two positive leading equivalents in
Theorem~\ref{thm:fixed-q-hodge-universality}. The nonzero leading
constant in the denominator ensures that the ratio is defined for
all sufficiently large $n$, including for complex Hodge arguments.
The inequality $0<\rho_q/\sigma_q<1$ gives exponential decay.
\end{proof}

\begin{remark}[Relation to the question of Kannan--Serpente--Yun]
\label{rem:hodge-KSY}
At $u=v=1$, we recover the Euler comparison of
\cite[Section~6.2]{CottiPaperI}, with $\sigma_1=1$ and $\rho_1=e-2$.
For every fixed $g\ge2$, the resulting equivalent for the all-light
to ordinary Euler ratio, and its limit zero, give the affirmative
answer to \cite[Question~5.3]{KSY} obtained there.
At $u=v=-1$, the same rates and constants describe the ratio of
total Betti numbers. More generally, for any fixed $s>0$, the
Poincar\'e evaluation $P(X_n;s)$ uses $u=v=-s$ and $q=s^2$
and obeys the same comparison.
\qrem
\end{remark}

\begin{remark}[Genus zero]
There is no stable minimal all-light family in genus zero. The
ordinary spaces $\overline M_{0,n}$ are Hodge--Tate, so their Hodge
polynomials depend exactly on $uv$ and $N_q$ vanishes for every $n$.
Their positive-product coefficient asymptotics are the genus-zero
formula of Section~\ref{sec:motivic-asymptotics}. The new concentration assertion concerns the
positive-genus families, whose full Hodge polynomials need not be
diagonal.
\qrem
\end{remark}

\appendix

\section{Universal triangular transformations and recurrences}
\label{app:triangular}

We give the coefficient form of the heavy/light transformation of
Section~\ref{sec:motivic-heavy-light}. The transformation coefficients
are independent of the genus and of the chosen motivic measure.

Fix $g\ge1$ and a motivic measure
$\mu:K_0(\Var_{\C})\to\mathsf R$, and put $\lambda=\mu(\LL)$.
Write
\[
 a_{g,n}^{\mu}:=\mu(\overline M_{g,(0^+)^n}),
 \qquad
 b_{g,n}^{\mu}:=\mu(\overline M_{g,n}),
\]
with $a_{1,0}^{\mu}=b_{1,0}^{\mu}=0$.
For $g\ge2$, both constant terms equal $\mu(\overline M_g)$.
By Corollary~\ref{cor:realized-composition},
\[
 G_g^\mu(x)=A_g^\mu(\Psi_\lambda(x)),
 \qquad \Psi_\lambda=\Phi_\lambda^{\langle-1\rangle}.
\]
Below, the universal coefficients are defined with $\lambda$
an indeterminate and evaluated at $\mu(\LL)$ when used
with the invariants above.

\subsection{The universal triangular coefficients}

Define
\begin{equation}\label{eq:Cnm-def}
 C_{n,m}(\lambda)
 :=\frac{n!}{m!}[x^n]\Psi_\lambda(x)^m,
 \qquad 0\le m\le n,
\end{equation}
and set $C_{n,m}=0$ outside this range.
In particular, $C_{0,0}=1$ and $C_{n,0}=0$ for $n>0$.

\begin{proposition}[Universal triangular transformation]
\label{prop:Cnm-transform}
For every $n\ge0$,
\begin{equation}\label{eq:triangular-transform}
 b_{g,n}^{\mu}
 =\sum_{m=0}^{n}
 C_{n,m}(\lambda)a_{g,m}^{\mu}.
\end{equation}
For $1\le m\le n$,
\begin{equation}\label{eq:Cnm-Lagrange}
 C_{n,m}(\lambda)
 =\frac{(n-1)!}{(m-1)!}
 [z^{n-m}]
 \left(\frac{z}{\Phi_\lambda(z)}\right)^n.
\end{equation}
Moreover, $C_{n,n}=1$ for $n\ge0$, and
$C_{n,n-1}=0$ for $n\ge1$.
\end{proposition}

\begin{proof}
Expand $G_g^\mu=A_g^\mu\circ\Psi_\lambda$ and compare
coefficients of $x^n/n!$. The second formula follows from
Lagrange--B\"urmann inversion, and the last assertions from
$\Psi_\lambda(x)=x+O(x^3)$.
\end{proof}

\subsection{Derangements and the universal coefficients}

The coefficients $C_{n,m}(\lambda)$ satisfy a recurrence
involving derangement polynomials. We first recall the relevant
terminology.

A permutation $\sigma\in\mathfrak S_r$ is a
\emph{derangement} if it has no fixed points, that is,
$\sigma(i)\ne i$ for every $i\in[r]$. Denote by
$\mathfrak D_r\subseteq\mathfrak S_r$ the set of derangements.
An index $i\in[r]$ is an \emph{excedance} of $\sigma$ if
$\sigma(i)>i$, and set
\[
 \operatorname{exc}(\sigma)
 :=\#\{i\in[r]:\sigma(i)>i\}.
\]
The derangement polynomials are
\begin{equation}\label{eq:q-derangement-def}
 d_r(\lambda)
 :=\sum_{\sigma\in\mathfrak D_r}
 \lambda^{\operatorname{exc}(\sigma)}.
\end{equation}
We include the empty permutation, so that $d_0=1$.
These polynomials were studied by F.\,Brenti
\cite{BrentiDer}; see also \cite{PetersenEulerian}.
The coefficient of $\lambda^j$ counts derangements
with exactly $j$ excedances.

Their classical exponential generating function is
\begin{equation}\label{eq:derangement-egf-classical}
 \sum_{r\ge0}d_r(\lambda)\frac{z^r}{r!}
 =
 \frac{1}{
  1-\displaystyle\sum_{r\ge2}
  (\lambda+\cdots+\lambda^{r-1})\frac{z^r}{r!}}.
\end{equation}
Since the denominator is $\Phi_\lambda'(z)$, this gives
\begin{equation}\label{eq:derangement-egf}
 \frac1{\Phi_\lambda'(z)}
 =
 \sum_{r\ge0}d_r(\lambda)\frac{z^r}{r!}.
\end{equation}
This identity also follows directly from the propagator.
The Eulerian polynomial $\mathrm A_r(\lambda)$ of
Section~\ref{sec:eulerian-propagator} counts permutations by excedances
as well as by descents; see \cite{StanleyEC1,PetersenEulerian}.
Deleting the fixed points of a permutation and relabeling the remaining
entries in increasing order preserves its number of excedances.
Choosing the fixed-point set therefore gives
\[
\mathrm A_r(\lambda)
=\sum_{j=0}^r\binom{r}{j}d_j(\lambda),
\qquad \mathrm A_0(\lambda):=1.
\]
By the exponential product rule and
\eqref{eq:Eulerian-chain-expansion}, this implies
\[
\sum_{r\ge0}d_r(\lambda)\frac{z^r}{r!}
=e^{-z}P_\lambda(z)=\frac1{\Phi_\lambda'(z)},
\]
where the last equality is the propagator identity of
Section~\ref{sec:propagator-coordinate}.
Thus the reciprocal derivative of the universal coordinate
enumerates derangements by their number of excedances.
The first terms are
\[
 d_0=1,\qquad d_1=0,\qquad d_2=\lambda,\qquad
 d_3=\lambda(1+\lambda),\qquad
 d_4=\lambda(1+7\lambda+\lambda^2).
\]

\begin{proposition}[Derangement recurrence]
\label{prop:Cnm-derangement}
For $n\ge0$ and $1\le m\le n+1$,
\begin{equation}\label{eq:Cnm-derangement}
 C_{n+1,m}(\lambda)
 =
 \sum_{k=m-1}^{n}
 \binom{k}{m-1}
 d_{k-m+1}(\lambda)C_{n,k}(\lambda).
\end{equation}
Together with $C_{0,0}=1$ and $C_{n,0}=0$ for $n>0$,
this determines the triangular coefficients. In particular,
\[
 C_{n,m}(\lambda)\in\Z_{\ge0}[\lambda].
\]
\end{proposition}

\begin{proof}
Differentiating $\Phi_\lambda(\Psi_\lambda(x))=x$ and using
\eqref{eq:derangement-egf} gives
\[
 \Psi_\lambda'(x)
 =
 \frac1{\Phi_\lambda'(\Psi_\lambda(x))}
 =
 \sum_{r\ge0}d_r(\lambda)
 \frac{\Psi_\lambda(x)^r}{r!}.
\]
Hence
\[
 \frac{d}{dx}\frac{\Psi_\lambda(x)^m}{m!}
 =
 \sum_{r\ge0}
 \binom{m+r-1}{m-1}d_r(\lambda)
 \frac{\Psi_\lambda(x)^{m+r-1}}{(m+r-1)!}.
\]
Comparing coefficients of $x^n/n!$ proves the recurrence.
The positivity assertion follows by induction.
\end{proof}

\begin{remark}[A forest interpretation]
The coefficients $C_{n,m}(\lambda)$ also admit a direct
combinatorial interpretation. The inverse-coordinate identity
can be written as
\[
 \Psi_\lambda(x)
 =
 x+\sum_{k\ge3}
 (\lambda+\cdots+\lambda^{k-2})
 \frac{\Psi_\lambda(x)^k}{k!}.
\]
This is the recursive generating-function equation for rooted
trees with labeled leaves and unordered children. Each internal
vertex has $k\ge3$ children and carries a decoration
$j\in\{1,\ldots,k-2\}$ of weight $\lambda^j$.
The weight of a tree is the product of its vertex weights;
a single leaf has weight $1$.

By \eqref{eq:Cnm-def}, $C_{n,m}(\lambda)$ is therefore the
weight enumerator of unordered forests of $m$ such trees
whose leaves are labeled by $[n]$.
\qrem
\end{remark}

\subsection{A bivariate refinement}

Define
\[
 F_{g,n}^{\mu}(s)
 :=
 \sum_{m=0}^{n}
 s^m C_{n,m}(\lambda)a_{g,m}^{\mu}.
\]
The variable $s$ records the all-light index $m$, and
$F_{g,n}^{\mu}(1)=b_{g,n}^{\mu}$.
The associated generating series is
\begin{equation}\label{eq:universal-bivariate-egf}
 \mathscr F_g^\mu(x,s)
 :=\sum_{n\ge0}F_{g,n}^{\mu}(s)\frac{x^n}{n!}
 =A_g^\mu(s\Psi_\lambda(x)).
\end{equation}

Define universal polynomials $\eta_j(\lambda)\in\Q[\lambda]$ by
\begin{equation}\label{eq:eta-def}
 \frac{x\Psi_\lambda'(x)}{\Psi_\lambda(x)}
 =
 \sum_{j\ge0}\eta_j(\lambda)\frac{x^j}{j!},
 \qquad \eta_0(\lambda)=1.
\end{equation}

\begin{proposition}[Universal bivariate recurrence]
\label{prop:universal-bivariate}
Put $\Delta_s=s\partial_s$. For every $n\ge1$,
\begin{equation}\label{eq:universal-bivariate-rec}
 (n-\Delta_s)F_{g,n}^{\mu}(s)
 =
 \sum_{j=1}^{n-1}
 \binom{n}{j}\eta_j(\lambda)
 \Delta_s F_{g,n-j}^{\mu}(s).
\end{equation}
Together with
\[
 [s^n]F_{g,n}^{\mu}(s)=a_{g,n}^{\mu},
 \qquad F_{g,0}^{\mu}=a_{g,0}^{\mu},
\]
this determines the refinement.
\end{proposition}

\begin{proof}
Differentiating \eqref{eq:universal-bivariate-egf} gives
\[
 x\partial_x\mathscr F_g^\mu
 =
 \frac{x\Psi_\lambda'(x)}{\Psi_\lambda(x)}
 \Delta_s\mathscr F_g^\mu.
\]
Using \eqref{eq:eta-def} and comparing coefficients proves
the recurrence. The term involving $F_{g,0}^{\mu}$ vanishes
because it is constant in $s$.

On the coefficient of $s^m$, the operator $n-\Delta_s$
acts by $n-m$. Thus the recurrence determines all coefficients
with $m<n$, while the coefficient of $s^n$ is the prescribed
diagonal datum $a_{g,n}^{\mu}$.
\end{proof}

\begin{remark}[Comparison with genus zero]
For the Poincar\'e realization, $\lambda=t^2$, this construction
is analogous to the genus-zero refinement of
\cite[Theorem~2.8]{Cotti}: setting $s=1$ recovers the ordinary
Poincar\'e polynomial, and the refinement satisfies a linear
recurrence expressed through $\Delta_s$.
The genus-zero construction uses the Bell-polynomial
decomposition of the Getzler--Manin formula, whereas the present
one uses the universal heavy/light coordinate.
\qrem
\end{remark}

\section{Poincar\'e resummation and the unilateral Bessel series}
\label{app:poincare-resummation}
\subsection{Proof of the Poincar\'e resummation}
\begin{lemma}\label{lem:cusp-dimensions}
Let
\(
\beta_{2k+2}
:=
\dim_{\C}S_{2k+2}\bigl(\mathrm{SL}_2(\Z)\bigr),
\) for \(k\geq0,\)
where $S_m(\mathrm{SL}_2(\Z))$ denotes the space of cusp forms of weight
$m$ for the full modular group. Then
\begin{equation}
\sum_{k\geq0}\beta_{2k+2}y^k
=
\frac{y^5}{(1-y^2)(1-y^3)}.
\label{eq:cusp-dim-series-P}
\end{equation}
\end{lemma}

\begin{proof}
Recall that the graded ring of holomorphic modular forms for the full
modular group is
\(
M_\bullet\bigl(\mathrm{SL}_2(\Z)\bigr)
=
\C[E_4,E_6],
\)
where $E_4$ and $E_6$ are the normalized Eisenstein series of weights
$4$ and $6$, respectively. Moreover, the graded ideal of cusp forms is
generated by the modular discriminant $\Delta$, the cusp form of weight
$12$ given, with $Q=e^{2\pi i\tau}$, by
\[
\Delta(\tau)
=
Q\prod_{n\geq1}(1-Q^n)^{24}
=
\frac{E_4(\tau)^3-E_6(\tau)^2}{1728}.
\]
Thus
\(
S_\bullet\bigl(\mathrm{SL}_2(\Z)\bigr)
=
\Delta\,\C[E_4,E_6].
\)
See, for instance,
\cite[Ch.~VII, \S\S3.2, 4.4--4.5]{Serre};
see also \cite[Ch.~3]{DS}.
If $y$ records half the modular weight, then $E_4$, $E_6$, and
$\Delta$ have degrees $2$, $3$, and $6$, respectively. Hence
\[
\sum_{m\geq0}
\dim_{\C}S_{2m}\bigl(\mathrm{SL}_2(\Z)\bigr)y^m
=
\frac{y^6}{(1-y^2)(1-y^3)}.
\]
Since \(S_0(\mathrm{SL}_2(\Z))=0\), dividing by \(y\) and setting \(m=k+1\) gives \eqref{eq:cusp-dim-series-P}.
\end{proof}

\begin{proof}[Proof of Theorem \ref{thm:poincare-resummation}]
We use standard identities for modified Bessel functions; see
\cite[Eqs.~(10.29.1), (10.35.1)--(10.35.2)]{NIST} and, classically,
\cite[Chs.~III and VI]{Watson}. The recurrence
\[
\mathsf I_{\nu-1}(x)-\mathsf I_{\nu+1}(x)
=
\frac{2\nu}{x}\mathsf I_\nu(x),
\]
applied with $x=2tz$ and $\nu=2k+1$, gives
\begin{equation}
\mathcal J_k(z;t^2)
=
t^{-2k}
\bigl(
\mathsf I_{2k}(2tz)-\mathsf I_{2k+2}(2tz)
\bigr).
\label{eq:J-even-Bessel}
\end{equation}
As in \eqref{eq:N-def}, the right-hand side is regular at $t=0$
when interpreted coefficientwise.

We first resum the terms in \eqref{eq:SP-original} which are independent
of the cusp-form dimensions $\beta_{2k+2}$. By
\eqref{eq:J-even-Bessel},
\[
t^2\mathcal J_0(z;t^2)
=
t^2
\bigl(
\mathsf I_0(2tz)-\mathsf I_2(2tz)
\bigr),
\]
while
\[
\begin{aligned}
\sum_{k\geq1}\mathcal J_k(z;t^2)
&=
\sum_{k\geq1}
t^{-2k}
\bigl(
\mathsf I_{2k}(2tz)-\mathsf I_{2k+2}(2tz)
\bigr)\\
&=
(1-t^2)\mathcal N(z;t)+\mathsf I_2(2tz).
\end{aligned}
\]
Combining these two expressions gives precisely
\eqref{eq:SP-even}.

We now turn to the cuspidal contribution. Set
\[
b_k:=\beta_{2k+2},\qquad k\geq0.
\]
Notice that $b_0=\beta_2=0$. By Lemma~\ref{lem:cusp-dimensions},
\[
\sum_{k\geq0}b_k y^k
=
\frac{y^5}{(1-y^2)(1-y^3)}.
\]
Using \eqref{eq:J-even-Bessel}, the cuspidal part of
\eqref{eq:SP-original} becomes
\[
2t\sum_{k\geq1}
b_k
\bigl(
\mathsf I_{2k}(2tz)-\mathsf I_{2k+2}(2tz)\bigr)=2t\sum_{k\geq1}\varepsilon_k\mathsf I_{2k}(2tz),
\]
where we introduced the differences
\[
\varepsilon_k:=b_k-b_{k-1},
\qquad k\geq1.
\]
Their generating function is
\begin{equation}
\sum_{k\geq1}\varepsilon_k y^k
=
(1-y)\frac{y^5}{(1-y^2)(1-y^3)}
=
-y
+\frac{y}{6(1-y)}
+\frac{y+2y^2}{3(1+y+y^2)}
+\frac{y}{2(1+y)}.
\label{eq:epsilon-gf}
\end{equation}
Using
\[
\sum_{k\geq1}\cos\frac{2\pi k}{3}\,y^k
=
-\frac{y+2y^2}{2(1+y+y^2)},
\]
we may equivalently write
\[
\sum_{k\geq1}\varepsilon_k y^k
=
-y
+\frac16\sum_{k\geq1}y^k
-\frac23\sum_{k\geq1}
\cos\frac{2\pi k}{3}\,y^k
-\frac12\sum_{k\geq1}(-1)^k y^k.
\]
Comparison of coefficients therefore gives
\begin{equation}
\varepsilon_k
=
\frac16
-\frac23\cos\frac{2\pi k}{3}
-\frac12(-1)^k
-\delta_{k,1},
\qquad k\geq1.
\label{eq:epsilon-explicit}
\end{equation}
Consequently,
\begin{equation}
\sum_{k\geq1}\varepsilon_k\mathsf I_{2k}(x)
={}
-\mathsf I_2(x)
+\frac16\sum_{k\geq1}\mathsf I_{2k}(x)
-\frac23\sum_{k\geq1}
\cos\frac{2\pi k}{3}\,\mathsf I_{2k}(x)
-\frac12\sum_{k\geq1}
(-1)^k\mathsf I_{2k}(x).
\label{eq:epsilon-Bessel}
\end{equation}

To evaluate the three sums, we use the even part of the
Jacobi--Anger expansion,
\begin{equation}
\cosh(x\cos\theta)
=
\mathsf I_0(x)
+
2\sum_{k\geq1}
\mathsf I_{2k}(x)\cos(2k\theta),
\label{eq:JA-even}
\end{equation}
which follows from \cite[Eq.~(10.35.2)]{NIST}. Specializing respectively
to $\theta=0,\pi/3,\pi/2$, we obtain
\begin{align*}
\cosh x
&=
\mathsf I_0(x)
+
2\sum_{k\geq1}\mathsf I_{2k}(x),\\
\cosh(x/2)
&=
\mathsf I_0(x)
+
2\sum_{k\geq1}
\cos\frac{2\pi k}{3}\,\mathsf I_{2k}(x),\\
1
&=
\mathsf I_0(x)
+
2\sum_{k\geq1}
(-1)^k\mathsf I_{2k}(x).
\end{align*}
Substitution into \eqref{eq:epsilon-Bessel} yields
\begin{equation}
\sum_{k\geq1}\varepsilon_k\mathsf I_{2k}(x)
=
\frac12\mathsf I_0(x)
-\mathsf I_2(x)
+\frac1{12}\cosh x
-\frac13\cosh(x/2)
-\frac14.
\label{eq:epsilon-resummed}
\end{equation}
Taking $x=2tz$ and multiplying by $2t$ proves
\eqref{eq:SP-odd}. Finally, adding \eqref{eq:SP-even} and
\eqref{eq:SP-odd} gives the integrand in
\eqref{eq:A1P-resummed}; integration from $0$ to $z$ completes the
proof.
\end{proof}

\subsection{Integral representation of the unilateral Bessel series}
\label{app:N-integral}

We record here an integral representation of the unilateral Bessel
series $\mathcal N(z;t)$ introduced in \eqref{eq:N-def}. Although this
representation is not needed for the main resummation theorem, it gives
a complementary analytic description of the only non-elementary
infinite sum appearing in \eqref{eq:A1P-resummed}.

\begin{proposition}[Integral representation of the unilateral Bessel series]
\label{prop:N-integral}
For $|t|>1$,
\begin{equation}
\mathcal N(z;t)
=
\frac1\pi
\int_0^\pi
e^{2tz\cos\theta}
\frac{
t^2\cos(2\theta)-1
}{
t^4-2t^2\cos(2\theta)+1
}\,d\theta .
\label{eq:N-integral-outside}
\end{equation}
For $|t|<1$,
\begin{equation}
\mathcal N(z;t)
=
\cosh((1+t^2)z)
-
\frac1\pi
\int_0^\pi
e^{2tz\cos\theta}
\frac{
1-t^2\cos(2\theta)
}{
1-2t^2\cos(2\theta)+t^4
}\,d\theta .
\label{eq:N-integral-inside}
\end{equation}
At the boundary values $t=\pm1$,
\begin{equation}
\mathcal N(z;\pm1)
=
\frac12
\left(
\cosh(2z)-\mathsf I_0(2z)
\right).
\label{eq:N-boundary}
\end{equation}
\end{proposition}

\begin{proof}
For integer $m\ge0$, the standard integral representation
\begin{equation}
\mathsf I_m(x)
=
\frac1\pi
\int_0^\pi
e^{x\cos\theta}\cos(m\theta)\,d\theta
\label{eq:I-integral}
\end{equation}
is \cite[Eq.~(10.32.3)]{NIST}; see also
\cite[Ch.~VI]{Watson}.

Assume first that $|t|>1$. The series in \eqref{eq:N-def} converges
uniformly under the integral in \eqref{eq:I-integral}, hence
\[
\mathcal N(z;t)
=
\frac1\pi
\int_0^\pi
e^{2tz\cos\theta}
\sum_{k\ge1}t^{-2k}\cos(2k\theta)\,d\theta.
\]
For $|r|<1$,
\begin{equation}
\sum_{k\ge1}r^k\cos(2k\theta)
=
\frac{r\cos(2\theta)-r^2}
{1-2r\cos(2\theta)+r^2}.
\label{eq:geometric-cos}
\end{equation}
Taking $r=t^{-2}$ proves \eqref{eq:N-integral-outside}.

For $|t|<1$, introduce
\[
\eu M(z;t)
:=
\sum_{k\ge1}
t^{2k}\mathsf I_{2k}(2tz).
\]
The bilateral generating function
\begin{equation}
e^{\frac{x}{2}(w+w^{-1})}
=
\sum_{m\in\Z}
w^m\mathsf I_m(x),
\label{eq:Bessel-bilateral}
\end{equation}
see \cite[Eq.~(10.35.1)]{NIST}, gives after taking its even part
\begin{equation}
\mathcal N(z;t)+\eu M(z;t)
=
\cosh((1+t^2)z)-\mathsf I_0(2tz).
\label{eq:N-M}
\end{equation}
On the other hand, \eqref{eq:I-integral} and
\eqref{eq:geometric-cos}, now with $r=t^2$, yield
\[
\eu M(z;t)
=
\frac1\pi
\int_0^\pi
e^{2tz\cos\theta}
\frac{
t^2\cos(2\theta)-t^4
}{
1-2t^2\cos(2\theta)+t^4
}\,d\theta.
\]
Using also
\[
\mathsf I_0(2tz)
=
\frac1\pi
\int_0^\pi
e^{2tz\cos\theta}\,d\theta
\]
in \eqref{eq:N-M} proves \eqref{eq:N-integral-inside}.

Finally, for $t=\pm1$ the even part of the Jacobi--Anger expansion
gives
\[
\cosh(2z)
=
\mathsf I_0(2z)
+
2\sum_{k\ge1}\mathsf I_{2k}(2z),
\]
which proves \eqref{eq:N-boundary}.
\end{proof}

\begin{remark}
The two integral representations in
Proposition~\ref{prop:N-integral} should not be identified by naive
analytic continuation across $|t|=1$: the poles of the rational kernel
meet the integration contour there. Formula \eqref{eq:N-boundary}
instead gives the boundary values $t=\pm1$ directly.
\qrem
\end{remark}

\section{Beyond the positive real Lefschetz regime}
\label{app:nonpositive-lambda}

Theorems~\ref{thm:motivic-genus-zero-asymptotics}--\ref{thm:motivic-ordinary-asymptotics}
and Corollary~\ref{cor:motivic-ratio} assume a positive real Lefschetz
value $\lambda$. The same assumption underlies the singularity analysis
in Proposition~\ref{prop:universal-singularity-pair} and the Hodge
specializations in Theorem~\ref{thm:fixed-q-hodge-universality}
and Corollary~\ref{cor:hodge-heavy-light-ratio}. This assumption is not merely technical: outside the positive real regime, the singularity geometry of the propagator and of the inverse heavy/light coordinate can change qualitatively. In this appendix we discuss the cases of zero, negative, and complex Lefschetz values, which require separate consideration and are not covered directly by the preceding transfer statements.

\phantomsection
\label{rem:lambda-zero}
The value $\lambda=0$ is degenerate. Indeed,
\[
P_0(z)=e^z,
\qquad
\Phi_0(z)=z,
\qquad
\Psi_0(t)=t.
\]
Thus the propagator has no finite pole and the inverse coordinate has
no critical singularity. The mechanism responsible for the
fixed-genus asymptotics of the positive-real analysis therefore
disappears completely.

For a nonzero complex point $\zeta$, use the rotated domain
\[
\Delta_\zeta(\eta,\phi)=\{z:|z|<|\zeta|+\eta,\ z\ne\zeta,
\ |\arg(z/\zeta-1)|>\phi\}.
\]
A multiple-$\Delta$ domain is the intersection of these domains at the
finitely many dominant singularities under consideration. As in
Section~\ref{sec:motivic-asymptotics}, continuation refers to the germ at the origin.

We first describe the pole geometry away from the negative real axis.
For $\lambda\in\C\setminus(-\infty,0]$, $\lambda\ne1$, choose the
principal logarithm and put
\[
\sigma_{\lambda,k}
:=\frac{\Log\lambda+2\pi i k}{\lambda-1},\qquad k\in\mathbb Z.
\]
Writing $\Log\lambda=x+iy$, with $|y|<\pi$, gives
\[
|\sigma_{\lambda,k}|^2
=\frac{x^2+(y+2\pi k)^2}{|\lambda-1|^2}.
\]
Consequently, $\sigma_{\lambda,0}$ is the unique pole of smallest
modulus. The value $\lambda=1$ is already covered by $P_1(z)=(1-z)^{-1}$.
For a genuinely complex Lefschetz value, the all-light transfer argument
therefore remains valid, with a complex dominant singularity and a
complex leading constant, provided the remainder is analytic in a
rotated $\Delta$ domain at $\sigma=\sigma_{\lambda,0}$ and satisfies
$\mathcal R_g^\mu(z)=O(|1-z/\sigma|^{-3g+4})$ there. Explicitly,
for $g\ge2$ and $m=3g-3$,
\[
\mu\bigl(\overline M_{g,(0^+)^n}\bigr)
=\frac{\gamma_g e^{2(g-1)(\lambda-2)\sigma}}
{\Gamma(m)(\lambda\sigma)^m}
\,n!\,\sigma^{-n}n^{m-1}
+O\bigl(n!\,|\sigma|^{-n}n^{m-2}\bigr).
\]
In genus one, if $S_1^\mu$ is analytic in the corresponding rotated
$\Delta$ domain and holomorphic across $\sigma$, logarithmic transfer
instead gives
\[
\mu\bigl(\overline M_{1,(0^+)^n}\bigr)
\sim\frac12 n!\,\sigma^{-n}n^{-1}.
\]
These are statements about the all-light series. For the ordinary series,
the location of the singularities of the inverse germ must still be
established independently. The negative real axis is the boundary case
of this pole geometry: the principal argument reaches $\pm\pi$, and
two poles have the same smallest modulus.

We next consider a negative real Lefschetz value. Write
\[
\lambda=-a,
\qquad a>0.
\]
The poles of the propagator
\[
P_\lambda(z)
=
\frac{\lambda-1}
{\lambda-e^{(\lambda-1)z}}
\]
are
\begin{equation}
\sigma_{\lambda,k}
=
-\frac{\log a+(2k+1)\pi i}{a+1},
\qquad
k\in\mathbb Z.
\label{eq:negative-lambda-poles}
\end{equation}
There is no longer a positive real dominant pole. Instead, the two
poles of smallest modulus form the conjugate pair
\begin{equation}
\sigma_{\lambda,\pm}
=
-\frac{\log a}{a+1}
\pm
\frac{\pi i}{a+1}.
\label{eq:negative-dominant-poles}
\end{equation}
Set
\begin{equation}
\delta_\lambda
:=
|\sigma_{\lambda,\pm}|
=
\frac{\sqrt{(\log a)^2+\pi^2}}{a+1}.
\label{eq:Rnegative}
\end{equation}
All other poles have modulus strictly larger than $\delta_\lambda$.

For a fixed genus $g\ge2$, put
\[
m:=3g-3,
\qquad
\alpha_{g,\lambda}:=2(g-1)(\lambda-2),
\]
and define
\begin{equation}
K_{g,\lambda}^{\pm}
:=
\frac{
\gamma_g e^{\alpha_{g,\lambda}\sigma_{\lambda,\pm}}
}{
(\lambda\sigma_{\lambda,\pm})^m
}.
\label{eq:negative-K}
\end{equation}

The analogue of Hypothesis $(\mathrm H_g)$ in this situation is a
multiple-singularity condition: we require
$\mathcal R_g^\mu$ to admit analytic continuation to a common
multiple-$\Delta$ domain at $\sigma_{\lambda,+}$ and
$\sigma_{\lambda,-}$, with outer radius greater than
$\delta_\lambda$, and to satisfy
\begin{equation*}
\mathcal R_g^\mu(z)
=
O\left(
\left|1-\frac z{\sigma_{\lambda,\pm}}\right|^{-m+1}
\right)
\qquad
(z\to\sigma_{\lambda,\pm}).
\tag{$\mathrm H_g^-$}
\label{eq:Hgminus}
\end{equation*}

\begin{proposition}[Oscillatory all-light asymptotics]
\label{prop:negative-alllight}
Let $\mu$ be a scalar motivic realization with
$\lambda=\mu(\mathbb L)<0$, and assume $(\mathrm H_g^-)$. Then
\begin{align}
\mu\bigl(\overline M_{g,(0^+)^n}\bigr)
={}&
\frac{n!\,n^{m-1}}{\Gamma(m)}
\left(
K_{g,\lambda}^+
\sigma_{\lambda,+}^{-n}
+
K_{g,\lambda}^-
\sigma_{\lambda,-}^{-n}
\right)
\nonumber\\
&\quad
+
O\left(
n!\delta_\lambda^{-n}n^{m-2}
\right).
\label{eq:negative-alllight-asymptotic}
\end{align}
Equivalently, writing
\[
\sigma_{\lambda,+}
=
\delta_\lambda e^{i\theta_\lambda},
\qquad
K_{g,\lambda}^+
=
|K_{g,\lambda}^+|e^{i\varphi_{g,\lambda}},
\]
one has
\begin{align}
\mu\bigl(\overline M_{g,(0^+)^n}\bigr)
={}&
\frac{
2|K_{g,\lambda}^+|
}{
\Gamma(m)
}
\,n!\delta_\lambda^{-n}n^{m-1}
\cos\bigl(
\varphi_{g,\lambda}-n\theta_\lambda
\bigr)
\nonumber\\
&\quad
+
O\left(
n!\delta_\lambda^{-n}n^{m-2}
\right).
\label{eq:negative-alllight-cosine}
\end{align}
\end{proposition}

\begin{proof}
At every pole $\sigma$ of $P_\lambda$,
\[
P_\lambda(z)
=
\frac1{\lambda\sigma}
\left(1-\frac z\sigma\right)^{-1}
+O(1).
\]
Hence the maximal trivalent contribution in
\eqref{eq:motivicgraphexpansion-analytic} has, at
$\sigma_{\lambda,\pm}$, the local form
\[
K_{g,\lambda}^{\pm}
\left(
1-\frac z{\sigma_{\lambda,\pm}}
\right)^{-m}
+
O\left(
\left(
1-\frac z{\sigma_{\lambda,\pm}}
\right)^{-m+1}
\right).
\]
Hypothesis $(\mathrm H_g^-)$ gives the same subdominant estimate for
the remainder. Since
$\sigma_{\lambda,+}$ and $\sigma_{\lambda,-}$ are the only
singularities on $|z|=\delta_\lambda$, the multiple-singularity transfer
theorem yields \eqref{eq:negative-alllight-asymptotic}.

For real $\lambda<0$ the two dominant singularities are conjugate, and
so are the constants in \eqref{eq:negative-K}. Combining the two
contributions gives \eqref{eq:negative-alllight-cosine}.
\end{proof}

\begin{remark}
Unlike the positive-$\lambda$ case, the leading term in
\eqref{eq:negative-alllight-cosine} may vanish along subsequences
because of the oscillatory cosine factor. For this reason the natural
statement is the additive asymptotic expansion
\eqref{eq:negative-alllight-asymptotic}, rather than a single
asymptotic equivalence.
\qrem
\end{remark}

The ordinary-weight problem is subtler. The critical points of
$\Phi_\lambda$ are again the poles
$\sigma_{\lambda,k}$ of the propagator, and their critical values are
\begin{equation}
\rho_{\lambda,k}
:=
\Phi_\lambda(\sigma_{\lambda,k})
=
e^{\sigma_{\lambda,k}}
-1-\frac1\lambda.
\label{eq:negative-critical-values}
\end{equation}
For $\lambda=-a<0$, this becomes
\begin{equation}
\rho_{\lambda,k}
=
a^{-1/(a+1)}
\exp\left(
-\frac{(2k+1)\pi i}{a+1}
\right)
+
\frac{1-a}{a}.
\label{eq:negative-critical-circle}
\end{equation}
Thus all the critical values lie on the circle
\begin{equation}
\left|
\rho-\frac{1-a}{a}
\right|
=
a^{-1/(a+1)}.
\label{eq:critical-circle}
\end{equation}
If $1/(a+1)$ is rational, only finitely many distinct critical values
occur; if it is irrational, they are dense on this circle.

This shows why Proposition~\ref{prop:universal-singularity-pair}
does not extend directly to negative $\lambda$. Positivity of the
coefficients in the implicit equation
\[
w=t+H_\lambda(w)
\]
was essential in identifying a unique dominant singularity of the
inverse branch $\Psi_\lambda$. For $\lambda<0$, one must instead
determine which critical values in
\eqref{eq:negative-critical-values} are actually reached by the branch
of $\Phi_\lambda^{-1}$ issued from the origin, and whether additional
singular values arising from the global inversion problem are present.
Consequently, the ordinary-weight asymptotics cannot in general be
deduced from the propagator poles alone.

More generally, suppose that the principal branch $\Psi_\lambda$ has
finitely many dominant singularities
$\tau_1,\dots,\tau_N$, each arising from a simple critical point
$\zeta_j$ of $\Phi_\lambda$, so that
\[
\Phi_\lambda(\zeta_j)=\tau_j,
\qquad
\Phi_\lambda'(\zeta_j)=0,
\qquad
\Phi_\lambda''(\zeta_j)\ne0.
\]
Then locally
\[
\Psi_\lambda(t)
=
\zeta_j
-
b_j
\left(1-\frac t{\tau_j}\right)^{1/2}
+
O\left(1-\frac t{\tau_j}\right),
\qquad
b_j^2
=
-\frac{2\tau_j}{\Phi_\lambda''(\zeta_j)},
\]
with the square root chosen according to the corresponding inverse
branch. Assume also that $G_g^\mu$ has analytic continuation to a common
multiple-$\Delta$ domain at the $\tau_j$, that the inverse germ has
no other singularities on its circle of convergence, and that the
composed remainder is of order
$O(|1-t/\tau_j|^{-(m-1)/2})$ at each $\tau_j$. Then the argument of
Theorem~\ref{thm:motivic-ordinary-asymptotics}  gives a sum of
contributions of the form
\[
\frac{
\gamma_g e^{\alpha_{g,\lambda}\zeta_j}
}{
(\lambda b_j)^m
\Gamma(m/2)
}
\,n!\,
\tau_j^{-n}
n^{m/2-1}.
\]
Thus multiple dominant inverse singularities lead naturally to
oscillatory ordinary-weight asymptotics as well.

\begin{example}[The value $\lambda=-1$]
\label{ex:lambda-minus-one}
For $\lambda=-1$ the universal coordinate simplifies drastically:
\[
\Phi_{-1}(z)=\sinh z,
\qquad
\Psi_{-1}(t)=\operatorname{arsinh}t,
\qquad
P_{-1}(z)=\frac{e^z}{\cosh z}.
\]
The dominant propagator poles are
\[
\sigma_\pm=\pm\frac{\pi i}{2},
\]
while the principal inverse has exactly the two dominant branch
points
\[
\rho_\pm=\pm i.
\]
Under the corresponding remainder continuation and subdominance
hypotheses, both the all-light and ordinary-weight generating series are
governed by conjugate pairs of dominant singularities, and the
preceding multiple-singularity transfer mechanism yields oscillatory
fixed-genus expansions in both chambers.
\qetr
\end{example}

\begin{remark}
The positive real regime is therefore distinguished by the existence
of a unique dominant real propagator pole and a unique dominant
square-root singularity of the principal inverse coordinate. Negative
real Lefschetz values instead lead naturally to conjugate or, more
generally, multiple dominant singularities and hence to oscillatory
asymptotics. For genuinely complex values of $\lambda$, the same
principle applies, but the dominant poles and singular values must be
determined from their moduli and from the global analytic continuation
of the chosen inverse branch.
\qrem
\end{remark}

\end{document}